# The difference between two Stieltjes constants

Donal F. Connon

dconnon@btopenworld.com

1 June 2009

**Abstract**

The Stieltjes constants $\gamma_n(u)$ are the coefficients of the Laurent expansion of the Hurwitz zeta function $\varsigma(s,u)$ about $s=1$ and surprisingly little is known about them. In this paper we derive some relations for the difference between two Stieltjes constants, two examples of which are shown below.

$$\gamma_1\left(\frac{1}{4}\right) - \gamma_1\left(\frac{3}{4}\right) = -\pi\left[\gamma + 4\log 2 + 3\log \pi - 4\log \Gamma\left(\frac{1}{4}\right)\right]$$

$$\gamma_2\left(\frac{1}{4}\right) - \gamma_2\left(\frac{3}{4}\right) =$$

$$2\pi\left[\varsigma''\left(0,\frac{1}{4}\right) - \varsigma''\left(0,\frac{3}{4}\right)\right] - 4\pi[\gamma + \log(8\pi)]\left[2\log \Gamma\left(\frac{1}{4}\right) - \log \pi - \frac{1}{2}\log 2\right]$$

$$+ \frac{\pi}{2}\left(2[\gamma + \log(8\pi)]^2 + \frac{\pi^2}{6}\right)$$

Using a formula previously obtained by Coffey [11], we may also obtain some specific values of the Stieltjes constants function $\gamma_n(u)$, for example

$$\gamma_1\left(\frac{1}{4}\right) = \frac{1}{2}[2\gamma_1 - 7\log^2 2 - 6\gamma \log 2] - \frac{1}{2}\pi\left[\gamma + 4\log 2 + 3\log \pi - 4\log \Gamma\left(\frac{1}{4}\right)\right]$$

$$\gamma_1\left(\frac{3}{4}\right) = \frac{1}{2}[2\gamma_1 - 7\log^2 2 - 6\gamma \log 2] + \frac{1}{2}\pi\left[\gamma + 4\log 2 + 3\log \pi - 4\log \Gamma\left(\frac{1}{4}\right)\right]$$

## 1. INTRODUCTION

We recall Hasse's formula [32] for the Hurwitz zeta function which is valid for all $s \in \mathbf{C}$ except $s=1$

(1.1)
$$(s-1)\varsigma(s,u) = \sum_{n=0}^{\infty} \frac{1}{n+1} \sum_{k=0}^{n} \binom{n}{k} \frac{(-1)^k}{(u+k)^{s-1}}$$

where we have the limit

(1.2)
$$\lim_{s \to 1}(s-1)\varsigma(s,u) = \sum_{n=0}^{\infty} \frac{1}{n+1} \sum_{k=0}^{n} \binom{n}{k}(-1)^k = \sum_{n=0}^{\infty} \frac{1}{n+1} \delta_{n,0} = 1$$

This limit is independent of $u$ and we therefore have

$$\frac{\partial}{\partial u} \lim_{s \to 1}[(s-1)\varsigma(s,u)] = \lim_{s \to 1} \frac{\partial}{\partial u}[(s-1)\varsigma(s,u)] = 0$$

We also note that

$$\frac{\partial}{\partial u}\varsigma(s,u) = \frac{\partial}{\partial u} \sum_{n=0}^{\infty} \frac{1}{(n+u)^s} = -s\varsigma(s+1,u)$$

and hence we have

$$\lim_{s \to 1} \frac{\partial}{\partial u}[(s-1)\varsigma(s,u)] = -\lim_{s \to 1}[s(s-1)\varsigma(s+1,u)] = 0$$

As a result of (1.2) we may apply L'Hôpital's rule to the limit

$$\lim_{s \to 1}\left[\frac{(s-1)\varsigma(s,u)-1}{s-1}\right] = \lim_{s \to 1}\left[\frac{(s-1)\varsigma'(s,u)+\varsigma(s,u)}{1}\right]$$

Differentiating (1.1) with respect to $s$ we see that with the notation $\varsigma'(s,u) = \frac{\partial}{\partial s}\varsigma(s,u)$

(1.3)
$$(s-1)\varsigma'(s,u) + \varsigma(s,u) = -\sum_{n=0}^{\infty} \frac{1}{n+1} \sum_{k=0}^{n} \binom{n}{k}(-1)^k \frac{\log(u+k)}{(u+k)^{s-1}}$$

and hence we have

$$\lim_{s \to 1}[(s-1)\varsigma'(s,u) + \varsigma(s,u)] = -\sum_{n=0}^{\infty} \frac{1}{n+1} \sum_{k=0}^{n} \binom{n}{k}(-1)^k \log(u+k)$$

Guillera and Sondow [29] showed in 2005 that the digamma function $\psi(u)$ may be expressed as

$$\text{(1.4)} \qquad \psi(u) = \sum_{n=0}^{\infty} \frac{1}{n+1} \sum_{k=0}^{n} \binom{n}{k} (-1)^k \log(u+k)$$

where the digamma function $\psi(u)$ is the logarithmic derivative of the gamma function $\psi(u) = \frac{d}{du} \log \Gamma(u)$. This identity was also independently derived in [20].

Hence we may then deduce that

$$\text{(1.5)} \qquad \lim_{s \to 1} [(s-1)\varsigma'(s,u) + \varsigma(s,u)] = -\psi(u)$$

and therefore conclude that

$$\text{(1.6)} \qquad \lim_{s \to 1} \left[ \frac{(s-1)\varsigma(s,u) - 1}{s-1} \right] = \lim_{s \to 1} \left[ \varsigma(s,u) - \frac{1}{s-1} \right] = -\psi(u)$$

The above limit is derived in a more complex manner in [42, p.91] and [45, p.271] using Hermite's integral formula for the Hurwitz zeta function.

With $u = 1$ we immediately see that

$$\text{(1.7)} \qquad \lim_{s \to 1} \left[ \varsigma(s,1) - \frac{1}{s-1} \right] = \lim_{s \to 1} \left[ \varsigma(s) - \frac{1}{s-1} \right] = -\psi(1) = \gamma$$

Differentiating (1.6) with respect to $u$ gives us

$$\frac{\partial}{\partial u} \lim_{s \to 1} \left[ \varsigma(s,u) - \frac{1}{s-1} \right] = -\psi'(u)$$

and hence we get

$$\frac{\partial}{\partial u} \lim_{s \to 1} \varsigma(s,u) = -\psi'(u)$$

We have by interchanging the order of the derivative and the limit operations

$$\frac{\partial}{\partial u} \lim_{s \to 1} \varsigma(s,u) = \lim_{s \to 1} \frac{\partial}{\partial u} \varsigma(s,u)$$

$$= \lim_{s \to 1} [-s\varsigma(s+1,u)] = -\varsigma(2,u)$$

giving us the well-known result

(1.8) $\qquad \varsigma(2,u) = \psi'(u)$.

The above formula may also be obtained directly by differentiating (1.4) and comparing the result with (1.1).

Using (1.4) we see that

$$\psi(u) - \log u = \sum_{n=0}^{\infty} \frac{1}{n+1} \sum_{k=0}^{n} \binom{n}{k} (-1)^k \log(u+k) - \log u$$

$$= \sum_{n=0}^{\infty} \frac{1}{n+1} \sum_{k=0}^{n} \binom{n}{k} (-1)^k \log(u+k) - \sum_{n=0}^{\infty} \frac{1}{n+1} \sum_{k=0}^{n} \binom{n}{k} (-1)^k \log u$$

$$= \sum_{n=0}^{\infty} \frac{1}{n+1} \sum_{k=0}^{n} \binom{n}{k} (-1)^k \log\left(1 + \frac{k}{u}\right)$$

and hence we formally obtain the well-known limit

$$\lim_{u \to \infty} [\psi(u) - \log u] = 0$$

From (1.6) we see that

(1.9) $\qquad \displaystyle\lim_{s \to 1} \left[ \varsigma(s,u) - \frac{1}{s-1} + \psi(u) \right] = 0$

and we have the Laurent expansion of the Hurwitz zeta function $\varsigma(s,u)$ about $s=1$

(1.10) $\qquad \displaystyle \varsigma(s,u) = \frac{1}{s-1} + \sum_{p=0}^{\infty} \frac{(-1)^p}{p!} \gamma_p(u)(s-1)^p$

where $\gamma_p(u)$ are known as the generalised Stieltjes constants and we have [49]

(1.11) $\qquad \gamma_0(u) = -\psi(u)$

We have for $p \geq 0$

(1.12) $\qquad \displaystyle \gamma_p(u) = (-1)^p \lim_{s \to 1} \left[ \varsigma^{(p)}(s,u) - \frac{(-1)^p p!}{(s-1)^{p+1}} \right]$

A derivation of the logarithmic series for the Stieltjes constants is shown below. We have from (1.10)

$$(s-1)\varsigma(s,u) = 1 + \sum_{p=0}^{\infty} \frac{(-1)^p}{p!} \gamma_p(u)(s-1)^{p+1}$$

Differentiation with respect to $s$ gives us

$$\frac{\partial}{\partial s}[(s-1)\varsigma(s,u)] = \sum_{p=0}^{\infty} \frac{(-1)^p}{p!} \gamma_p(u)(p+1)(s-1)^p = \gamma_0(u) + \sum_{p=1}^{\infty} \frac{(-1)^p}{p!} \gamma_p(u)(p+1)(s-1)^p$$

and as $s \to 1$ we have using (1.5)

$$-\psi(u) = \gamma_0(u)$$

A further differentiation gives us

$$\frac{\partial^2}{\partial s^2}[(s-1)\varsigma(s,u)] = \sum_{p=0}^{\infty} \frac{(-1)^p}{p!} \gamma_p(u)(p+1)p(s-1)^{p-1}$$

$$= -2\gamma_1(u) + \sum_{p=2}^{\infty} \frac{(-1)^p}{p!} \gamma_p(u)(p+1)p(s-1)^{p-1}$$

As $s \to 1$ we have

(1.13) $$\lim_{s \to 1}[(s-1)\varsigma''(s,u) + 2\varsigma'(s,u)] = \sum_{n=0}^{\infty} \frac{1}{n+1} \sum_{k=0}^{n} \binom{n}{k} (-1)^k \log^2(k+u)$$

and therefore

(1.14) $$\sum_{n=0}^{\infty} \frac{1}{n+1} \sum_{k=0}^{n} \binom{n}{k} (-1)^k \log^2(k+u) = -2\gamma_1(u)$$

The following result is readily derived

$$\frac{d^{m+1}}{ds^{m+1}}[(s-1)f(s)] = (s-1)f^{(m+1)}(s) + (m+1)f^{(m)}(s)$$

and we have

$$\frac{d^{m+1}}{ds^{m+1}}(s-1)\left[\sum_{p=0}^{\infty} \frac{(-1)^p}{p!} \gamma_p(s-1)^p\right]\bigg|_{s=1} = (m+1)\frac{d^m}{ds^m}\left[\sum_{p=0}^{\infty} \frac{(-1)^p}{p!} \gamma_p(s-1)^p\right]\bigg|_{s=1}$$

$$= (m+1)\sum_{p=0}^{\infty} \frac{(-1)^p}{p!} \gamma_p(u) p(p-1)...(p+1-m)(s-1)^{p-m}\bigg|_{s=1}$$

When $s=1$ the only non-vanishing term is when $p=m$ and we therefore we see that

(1.15)  $$\lim_{s\to 1}\frac{\partial^{m+1}}{\partial s^{m+1}}[(s-1)\varsigma(s,u)] = (m+1)(-1)^m \gamma_m(u)$$

Using the Hasse identity (1.1) we have

$$\lim_{s\to 1}\frac{\partial^{m+1}}{\partial s^{m+1}}[(s-1)\varsigma(s,u)] = (-1)^{m+1}\sum_{n=0}^{\infty}\frac{1}{n+1}\sum_{k=0}^{n}\binom{n}{k}(-1)^k \log^{m+1}(u+k)$$

This then gives us for $m \geq 0$

(1.16)  $$\gamma_m(u) = -\frac{1}{m+1}\sum_{n=0}^{\infty}\frac{1}{n+1}\sum_{k=0}^{n}\binom{n}{k}(-1)^k \log^{m+1}(u+k)$$

We see that

$$\frac{d}{du}\gamma_0(u) = -\sum_{n=0}^{\infty}\frac{1}{n+1}\sum_{k=0}^{n}\binom{n}{k}\frac{(-1)^k}{u+k}$$

$$= -\varsigma(2,u) = -\sum_{n=0}^{\infty}\frac{1}{(n+u)^2}$$

where we have employed the Hasse formula (1.1). Integration then results in the familiar formula for the digamma function [42, p.14]

(1.17)  $$\gamma_0(x) - \gamma_0(1) = \sum_{n=0}^{\infty}\left[\frac{1}{n+x} - \frac{1}{n+1}\right] = -(x-1)\sum_{n=0}^{\infty}\frac{1}{(n+1)(n+x)} = -\psi(x) - \gamma$$

Similarly we find that

$$\frac{d}{du}\gamma_1(u) = -\sum_{n=0}^{\infty}\frac{1}{n+1}\sum_{k=0}^{n}\binom{n}{k}\frac{(-1)^k \log(u+k)}{u+k}$$

The Hasse formula (1.1) gives us

$$\frac{\partial}{\partial s}[(s-1)\varsigma(s,u)] = (s-1)\varsigma'(s,u) + \varsigma(s,u) = -\sum_{n=0}^{\infty}\frac{1}{n+1}\sum_{k=0}^{n}\binom{n}{k}\frac{(-1)^k \log(u+k)}{(u+k)^{s-1}}$$

and with $s = 2$ we have

$$\varsigma'(2,u) + \varsigma(2,u) = -\sum_{n=0}^{\infty} \frac{1}{n+1} \sum_{k=0}^{n} \binom{n}{k} \frac{(-1)^k \log(u+k)}{u+k}$$

$$\frac{d}{du}\gamma_1(u) = \varsigma'(2,u) + \varsigma(2,u) = -\sum_{n=0}^{\infty} \frac{\log(n+u)}{(n+u)^2} + \sum_{n=0}^{\infty} \frac{1}{(n+u)^2} = \sum_{n=0}^{\infty} \frac{1 - \log(n+u)}{(n+u)^2}$$

We may easily evaluate the following integral

$$\int_{1}^{x} \frac{1 - \log(n+u)}{(n+u)^2} du = \frac{\log(n+u)}{n+u} \bigg|_{1}^{x}$$

and hence we have upon integrating

(1.18) $$\gamma_1(x) - \gamma_1(1) = \sum_{n=0}^{\infty} \left[ \frac{\log(n+x)}{n+x} - \frac{\log(n+1)}{n+1} \right]$$

With $x = 2$, since the series telescopes, we see that

$$\gamma_1(2) - \gamma_1(1) = \sum_{n=0}^{\infty} \left[ \frac{\log(n+2)}{n+2} - \frac{\log(n+1)}{n+1} \right] = 0$$

and hence we have

$$\gamma_1(2) = \gamma_1(1)$$

which is a particular case of (2.17.1).

Ivić [33, p.41] reports that the coefficients $\gamma_p(x)$ for $x \in (0,1]$ may be expressed as

(1.19) $$\gamma_p(x) = \lim_{N \to \infty} \left( \sum_{n=0}^{N} \frac{\log^p(n+x)}{n+x} - \frac{\log^{p+1}(N+x)}{p+1} \right)$$

and in particular we have

(1.20) $$\gamma_0(x) = \lim_{N \to \infty} \left( \sum_{n=0}^{N} \frac{1}{n+x} - \log(N+x) \right)$$

(1.21) $$\gamma_1(x) = \lim_{N \to \infty} \left( \sum_{n=0}^{N} \frac{\log(n+x)}{n+x} - \frac{1}{2}\log^2(N+x) \right)$$

We also see that

$$\gamma_p(1) = \lim_{N\to\infty}\left(\sum_{n=0}^{N} \frac{\log^p(n+1)}{n+1} - \frac{\log^{p+1}(N+1)}{p+1}\right)$$

$$= \lim_{N\to\infty}\left(\sum_{n=1}^{N} \frac{\log^p(n+1)}{n+1} - \frac{\log^{p+1}(N+1)}{p+1}\right)$$

$$= \lim_{N\to\infty}\left(\sum_{k=1}^{N} \frac{\log^p k}{k} + \frac{\log^p(N+1)}{N+1} - \frac{\log^{p+1}(N+1)}{p+1}\right)$$

and since $\lim_{N\to\infty}\frac{\log^p(N+1)}{N+1} = 0$ this becomes

$$= \lim_{N\to\infty}\left(\sum_{k=1}^{N} \frac{\log^p k}{k} - \frac{\log^{p+1}(N+1)}{p+1}\right)$$

$$= \lim_{N\to\infty}\left(\sum_{k=1}^{N} \frac{\log^p k}{k} - \frac{\log^{p+1} N}{p+1} + \frac{\log^{p+1} N}{p+1} - \frac{\log^{p+1}(N+1)}{p+1}\right)$$

Therefore, assuming that $\lim_{N\to\infty}\left[\log^{p+1} N - \log^{p+1}(N+1)\right] = 0$ (as shown below), we have the well-known result (see for example the 1955 paper, "The power series coefficients of $\varsigma(s)$", by Briggs and Chowla [9]

(1.22) $\quad \gamma_p = \gamma_p(1) = \lim_{N\to\infty}\left(\sum_{k=1}^{N} \frac{\log^p k}{k} - \frac{\log^{p+1} N}{p+1}\right) = \lim_{N\to\infty}\left(\sum_{k=1}^{N} \frac{\log^p k}{k} - \int_{1}^{N} \frac{\log^p x}{x} dx\right)$

Using (1.22) gives us

$$\gamma_1(x) - \gamma_1(1) = \lim_{N\to\infty}\left(\sum_{n=0}^{N}\left[\frac{\log(n+x)}{n+x} - \frac{\log(n+1)}{n+1}\right] - \frac{1}{2}\log^2(N+x) + \frac{1}{2}\log^2(N+1)\right)$$

Define $f(x) = \log^{p+1}(N+x)$ and, by the mean value theorem of calculus, we have

$$\log^{p+1}(N+x) - \log^{p+1}(N+1) = (x-1)(p+1)\frac{\log^p(N+\alpha)}{N+\alpha}$$

where $1 < \alpha < x$. Successive applications of L'Hôpital's rule then readily shows that

(1.22.1) $$\lim_{N\to\infty} \frac{\log^p(N+\alpha)}{N+\alpha} = 0$$

Hence we have

(1.23) $$\gamma_p(x) - \gamma_p(1) = \sum_{n=0}^{\infty}\left[\frac{\log^p(n+x)}{n+x} - \frac{\log^p(n+1)}{n+1}\right]$$

and in particular we have (as previously seen above)

(1.24) $$\gamma_0(x) - \gamma_0(1) = \sum_{n=0}^{\infty}\left[\frac{1}{n+x} - \frac{1}{n+1}\right] = -(x-1)\sum_{n=0}^{\infty}\frac{1}{(n+1)(n+x)} = -\gamma - \psi(x)$$

where we have employed (1.17).

Differentiating (1.23) gives us

$$\gamma_p'(x) = \sum_{n=0}^{\infty}\left[\frac{p\log^{p-1}(n+x)}{(n+x)^2} - \frac{\log^p(n+x)}{(n+x)^2}\right]$$

and hence we have

(1.25) $$\gamma_p'(x) = (-1)^{p+1}\left[p\varsigma^{(p-1)}(2,x) + \varsigma^{(p)}(2,x)\right]$$

In particular, as noted above, we see that

(1.26) $$\gamma_1'(x) = \varsigma(2,x) + \varsigma'(2,x)$$

Higher derivatives may also be computed in a similar manner. This result was subsequently obtained by Coffey [17a] by a different method.

□

There are several proofs of (1.22). For example, Bohman and Fröberg [7c] noted that

$$(s-1)\varsigma(s) = \sum_{k=1}^{\infty}\frac{(s-1)}{k^s} \quad \text{and} \quad \sum_{k=1}^{\infty}\left[\frac{1}{k^{s-1}} - \frac{1}{(k+1)^{s-1}}\right] = 1$$

and, assuming that $s$ is real and greater than 1, the above two equations may be subtracted to give

(1.27) $$(s-1)\varsigma(s) = 1 + \sum_{k=1}^{\infty}\left[\frac{1}{(k+1)^{s-1}} - \frac{1}{k^{s-1}} + \frac{(s-1)}{k^s}\right]$$

The well-known limit may be immediately derived from the above

$$\lim_{s \to 1}(s-1)\varsigma(s) = 1$$

Equation (1.27) may be written as

$$(s-1)\varsigma(s) = 1 + \sum_{k=1}^{\infty}\left[\exp(-(s-1)\log(k+1)) - \exp(-(s-1)\log k) + (s-1)k^{-1}\exp(-(s-1)\log k)\right]$$

$$= 1 + \sum_{k=1}^{\infty}\left[\sum_{n=0}^{\infty}\frac{(-1)^n(s-1)^n}{n!}\left(\log^n(k+1) - \log^n k\right) + \frac{s-1}{k}\sum_{n=0}^{\infty}\frac{(-1)^n(s-1)^n}{n!}\log^n k\right]$$

Dividing by $(s-1)$ we get

$$\varsigma(s) = \frac{1}{s-1} + \sum_{n=0}^{\infty}\frac{(-1)^n}{n!}\gamma_n(s-1)^n$$

where

(1.28) $$\gamma_n = \sum_{k=1}^{\infty}\left[\frac{\log^n k}{k} - \frac{\log^{n+1}(k+1) - \log^{n+1} k}{n+1}\right] = \sum_{k=1}^{\infty}\left[\frac{\log^n k}{k} - \int_k^{k+1}\frac{\log^n t}{t}dt\right]$$

which is equivalent to (1.22). A different derivation is given in (11.21) below.

## 2. THE DIFFERENCE BETWEEN TWO STIELTJES CONSTANTS

We recall Rademacher's formula [5, p.261] for the Hurwitz zeta function where for all $s$ and $1 \leq p \leq q$ where $p$ and $q$ are positive integers

(2.1) $$\varsigma\left(s, \frac{p}{q}\right) = 2\Gamma(1-s)(2\pi q)^{s-1}\sum_{j=1}^{q}\sin\left(\frac{\pi s}{2} + \frac{2\pi jp}{q}\right)\varsigma\left(1-s, \frac{j}{q}\right)$$

$$= \frac{\pi}{\Gamma(s)}(2\pi q)^{s-1}\sum_{j=1}^{q}\left[\sec\left(\frac{\pi s}{2}\right)\cos\left(\frac{2\pi jp}{q}\right) + \csc\left(\frac{\pi s}{2}\right)\sin\left(\frac{2\pi jp}{q}\right)\right]\varsigma\left(1-s, \frac{j}{q}\right)$$

where letting $p = q$ results in the functional equation for the Riemann zeta function.

$$\varsigma(s,1) = \varsigma(s) = 2\Gamma(1-s)(2\pi)^{s-1} \sin\left(\frac{\pi s}{2}\right) \varsigma(1-s)$$

A straightforward subtraction gives us

$$\varsigma\left(s, \frac{p}{q}\right) - \varsigma\left(s, 1 - \frac{p}{q}\right) =$$

$$2\Gamma(1-s)(2\pi q)^{s-1} \sum_{j=1}^{q} \left[\sin\left(\frac{\pi s}{2} + \frac{2\pi jp}{q}\right) - \sin\left(\frac{\pi s}{2} + \frac{2\pi j(q-p)}{q}\right)\right] \varsigma\left(1-s, \frac{j}{q}\right)$$

The familiar trigonometric identity immediately gives us

$$\sin\left(\frac{\pi s}{2} + \frac{2\pi jp}{q}\right) - \sin\left(\frac{\pi s}{2} + \frac{2\pi j(q-p)}{q}\right) = 2\cos\left(\frac{\pi s}{2}\right)\sin\left(\frac{2\pi jp}{q}\right)$$

and hence we obtain

(2.1.1) $$\varsigma\left(s, \frac{p}{q}\right) - \varsigma\left(s, 1 - \frac{p}{q}\right) = 4\Gamma(1-s)(2\pi q)^{s-1} \cos\left(\frac{\pi s}{2}\right) \sum_{j=1}^{q} \sin\left(\frac{2\pi jp}{q}\right) \varsigma\left(1-s, \frac{j}{q}\right)$$

For convenience we designate

$$f(s) = \Gamma(1-s)(2\pi q)^{s-1} \cos\left(\frac{\pi s}{2}\right)$$

and logarithmic differentiation results in

$$\frac{f'(s)}{f(s)} = -\psi(1-s) + \log(2\pi q) - \frac{\pi}{2}\tan\left(\frac{\pi s}{2}\right)$$

With the definition $\varsigma'(s,a) = \frac{\partial}{\partial s}\varsigma(s,a)$, we then see that

(2.2) $$\varsigma'\left(s, \frac{p}{q}\right) - \varsigma'\left(s, 1 - \frac{p}{q}\right) = -4\Gamma(1-s)(2\pi q)^{s-1} \cos\left(\frac{\pi s}{2}\right) \sum_{j=1}^{q} \sin\left(\frac{2\pi jp}{q}\right) \varsigma'\left(1-s, \frac{j}{q}\right)$$

$$-4\Gamma(1-s)(2\pi q)^{s-1} \cos\left(\frac{\pi s}{2}\right)\left[\psi(1-s) - \log(2\pi q) + \frac{\pi}{2}\tan\left(\frac{\pi s}{2}\right)\right] \sum_{j=1}^{q} \sin\left(\frac{2\pi jp}{q}\right) \varsigma\left(1-s, \frac{j}{q}\right)$$

and we now wish to consider the limit as $s \to 1$. We designate

$$\varsigma'\left(1, \frac{p}{q}\right) - \varsigma'\left(1, 1 - \frac{p}{q}\right) = \lim_{s \to 1}\left[\varsigma'\left(s, \frac{p}{q}\right) - \varsigma'\left(s, 1 - \frac{p}{q}\right)\right]$$

and, with reference to (2.2), we first of all we consider

$$\lim_{s \to 1}\left[\Gamma(1-s)(2\pi q)^{s-1}\cos\left(\frac{\pi s}{2}\right)\sum_{j=1}^{q}\sin\left(\frac{2\pi jp}{q}\right)\varsigma'\left(1-s, \frac{j}{q}\right)\right]$$

$$= \lim_{s \to 1}\left[\Gamma(1-s)\cos\left(\frac{\pi s}{2}\right)\right]\sum_{j=1}^{q}\sin\left(\frac{2\pi jp}{q}\right)\varsigma'\left(0, \frac{j}{q}\right)$$

Using Euler's reflection formula for the gamma function

$$\Gamma(s)\Gamma(1-s) = \frac{\pi}{\sin \pi s}$$

we see that

(2.3) $$\Gamma(1-s)\cos\left(\frac{\pi s}{2}\right) = \cos\left(\frac{\pi s}{2}\right)\frac{\pi}{\Gamma(s)\sin \pi s} = \frac{\pi}{2}/\Gamma(s)\sin\left(\frac{\pi s}{2}\right)$$

and we therefore have

(2.4) $$\lim_{s \to 1}\left[\Gamma(1-s)\cos\left(\frac{\pi s}{2}\right)\right] = \frac{\pi}{2}$$

Using Lerch's identity [7] for $x > 0$

(2.5) $$\varsigma'(0, x) = \log \Gamma(x) - \frac{1}{2}\log(2\pi)$$

we see that

$$\sum_{j=1}^{q}\sin\left(\frac{2\pi jp}{q}\right)\varsigma'\left(0, \frac{j}{q}\right) = \sum_{j=1}^{q}\sin\left(\frac{2\pi jp}{q}\right)\log \Gamma\left(\frac{j}{q}\right) - \frac{1}{2}\log(2\pi)\sum_{j=1}^{q}\sin\left(\frac{2\pi jp}{q}\right)$$

and since (as proved in (2.43) below)

(2.6) $$\sum_{j=1}^{q}\sin\left(\frac{2\pi jp}{q}\right) = 0$$

this becomes

$$\sum_{j=1}^{q}\sin\left(\frac{2\pi jp}{q}\right)\varsigma'\left(0,\frac{j}{q}\right)=\sum_{j=1}^{q}\sin\left(\frac{2\pi jp}{q}\right)\log\Gamma\left(\frac{j}{q}\right)$$

Accordingly, the first term on the right-hand side of (2.2) may be written as

$$-2\pi\sum_{j=1}^{q}\sin\left(\frac{2\pi jp}{q}\right)\log\Gamma\left(\frac{j}{q}\right)$$

We now consider the second term on the right-hand side of (2.2). We need to determine the limit

$$\lim_{s\to 1}\left[\psi(1-s)+\frac{\pi}{2}\tan\left(\frac{\pi s}{2}\right)\right]$$

A straightforward derivation of this limit may be obtained as follows. Since [42, p.14]

$$\psi(s)-\psi(1-s)=-\pi\cot\pi s$$

and since $\cot x-\tan x=\dfrac{\cos x}{\sin x}-\dfrac{\sin x}{\cos x}=\dfrac{\cos^2 x-\sin^2 x}{\sin x\cos x}=2\cot 2x$ we have

$$\cot\pi s=\frac{1}{2}\cot(\pi s/2)-\frac{1}{2}\tan(\pi s/2)$$

Hence we may write

(2.7) $$\psi(1-s)+\frac{\pi}{2}\tan\left(\frac{\pi s}{2}\right)=\psi(s)+\frac{\pi}{2}\cot\left(\frac{\pi s}{2}\right)$$

and the limit becomes obvious.

(2.8) $$\lim_{s\to 1}\left[\psi(1-s)+\frac{\pi}{2}\tan\left(\frac{\pi s}{2}\right)\right]=-\gamma$$

This also conforms with the formula in [44, p.20] that in the neighbourhood of $s=1$

(2.9) $$\frac{\pi}{2}\tan\left(\frac{\pi s}{2}\right)=-\frac{1}{s-1}+O(|s-1|)$$

(2.10) $$\lim_{s\to 1}\left[\frac{\pi}{2}\tan\left(\frac{\pi s}{2}\right)-\frac{\varsigma'(s)}{\varsigma(s)}\right]=\gamma$$

We have the well-known identity [5, p.264]

$$\varsigma(0,x) = \frac{1}{2} - x$$

and thus we see that

$$\sum_{j=1}^{q} \sin\left(\frac{2\pi jp}{q}\right)\varsigma\left(0,\frac{j}{q}\right) = \frac{1}{2}\sum_{j=1}^{q}\sin\left(\frac{2\pi jp}{q}\right) - \frac{1}{q}\sum_{j=1}^{q} j\sin\left(\frac{2\pi jp}{q}\right)$$

$$= -\frac{1}{q}\sum_{j=1}^{q} j\sin\left(\frac{2\pi jp}{q}\right)$$

We also have (as proved in (2.44) below)

(2.11) $$\sum_{j=1}^{q} j\sin\left(\frac{2\pi jp}{q}\right) = -\frac{q}{2}\cot\left(\frac{\pi p}{q}\right) \text{ for } p < q$$

It may be noted that this is the only place in the proof where we require the strict inequality $p < q$.

Therefore, the second term on the right-hand side of (2.2) may be written as

$$\pi[\log(2\pi q) + \gamma]\cot\left(\frac{\pi p}{q}\right)$$

and hence, as previously reported by Adamchik [2], we obtain

(2.12) $$\varsigma'\left(1,\frac{p}{q}\right) - \varsigma'\left(1,1-\frac{p}{q}\right) = \pi[\log(2\pi q) + \gamma]\cot\left(\frac{\pi p}{q}\right) - 2\pi\sum_{j=1}^{q-1}\log\Gamma\left(\frac{j}{q}\right)\sin\left(\frac{2\pi jp}{q}\right)$$

where $p$ and $q$ are positive integers and $p < q$. Adamchik [2] notes that this formula was first proved by Almkvist and Meurman. Adamchik's derivation [2] is rather terse and it is hoped that this expanded exposition may be of interest.

In particular we have

(2.13) $$\varsigma'\left(1,\frac{1}{4}\right) - \varsigma'\left(1,\frac{3}{4}\right) = \pi\left[\gamma + 4\log 2 + 3\log\pi - 4\log\Gamma\left(\frac{1}{4}\right)\right]$$

(2.14) $$\varsigma'\left(1,\frac{1}{3}\right) - \varsigma'\left(1,\frac{2}{3}\right) = \frac{\pi}{2\sqrt{3}}\left[2\gamma - \log 3 + 8\log(2\pi) - 12\log\Gamma\left(\frac{1}{3}\right)\right]$$

We recall the Laurent expansion (1.10) of the Hurwitz zeta function

$$\varsigma(s,u) = \frac{1}{s-1} + \sum_{p=0}^{\infty} \frac{(-1)^p}{p!} \gamma_p(u)(s-1)^p$$

and we may write this as

$$\varsigma(s,u) - \varsigma(s,1-u) = \sum_{p=0}^{\infty} \frac{(-1)^p}{p!} [\gamma_p(u) - \gamma_p(1-u)](s-1)^p$$

Differentiation with respect to $s$ results in

(2.15) $\qquad \varsigma'(s,u) - \varsigma'(s,1-u) = \sum_{p=0}^{\infty} \frac{(-1)^p}{p!} p[\gamma_p(u) - \gamma_p(1-u)](s-1)^{p-1}$

and in the limit as $s \to 1$ we have

(2.16) $\qquad \varsigma'(1,u) - \varsigma'(1,1-u) = -[\gamma_1(u) - \gamma_1(1-u)]$

In fact, as previously noted in equation (4.3.228b) in [21], more generally we have

(2.17) $\qquad \gamma_p(x) - \gamma_p(y) = \lim_{s \to 1} (-1)^p \frac{\partial^p}{\partial s^p} [\varsigma(s,x) - \varsigma(s,y)]$

It is easily seen from the definition of the Hurwitz zeta function that

$$\varsigma(s, x+n) = \varsigma(s,x) - \sum_{k=0}^{n-1} \frac{1}{(k+x)^s}$$

and in particular

$$\varsigma(s, x+1) = \varsigma(s,x) - \frac{1}{x^s}$$

Differentiation gives us

$$\varsigma'(s, x+1) = \varsigma'(s,x) + \frac{\log x}{x^s}$$

and in particular

$$\varsigma'(s,2) = \varsigma'(s,1) = \varsigma'(s)$$

We then see from (2.17) that

(2.17.1) $\quad \gamma_p(2) = \gamma_p(1)$

We obtain from (2.13) and (2.17)

(2.18) $\quad \gamma_1\left(\dfrac{1}{4}\right) - \gamma_1\left(\dfrac{3}{4}\right) = -\pi\left[\gamma + 4\log 2 + 3\log \pi - 4\log \Gamma\left(\dfrac{1}{4}\right)\right]$

We showed in [23b] that

(2.18.1) $\quad \displaystyle\sum_{r=0}^{q-1} \gamma_p\left(\dfrac{r+x}{q}\right) = q(-1)^p \dfrac{\log^{p+1} q}{p+1} + q\sum_{k=0}^{p}\binom{p}{k}(-1)^k \gamma_{p-k}(x)\log^k q$

and noting that

$$\sum_{r=0}^{q-1} f\left(\dfrac{r+x}{q}\right) = \sum_{m=1}^{q} f\left(\dfrac{m-1+x}{q}\right) = \sum_{m=1}^{q-1} f\left(\dfrac{m-1+x}{q}\right) + f\left(\dfrac{q-1+x}{q}\right)$$

we see that for integers $q \geq 2$ and $x = 1$

(2.19) $\quad \displaystyle\sum_{r=1}^{q-1} \gamma_p\left(\dfrac{r}{q}\right) = -\gamma_p + q(-1)^p \dfrac{\log^{p+1} q}{p+1} + q\sum_{j=0}^{p}\binom{p}{j}(-1)^j \gamma_{p-j}\log^j q$

which was previously derived by Coffey [11] using the following relation in [30]

(2.20) $\quad \displaystyle\sum_{r=1}^{q-1} \varsigma\left(s,\dfrac{r}{q}\right) = (q^s - 1)\varsigma(s)$

The relationship (2.19) has some history in that the case $p = 1$ was originally discovered by Ramanujan [7a, p.198]; as regards this, see (11.20) below.

Taking the simplest case $q = 2$ we get

(2.21) $\quad \gamma_p\left(\dfrac{1}{2}\right) = -\gamma_p + 2(-1)^p \dfrac{\log^{p+1} 2}{p+1} + 2\sum_{j=0}^{p}\binom{p}{j}(-1)^j \gamma_{p-j}\log^j 2$

and with $p = 0$ we have

(2.22) $$\gamma_0\left(\frac{1}{2}\right) = \gamma + 2\log 2 = -\psi\left(\frac{1}{2}\right)$$

in agreement with (1.11).

With $p = 1$ we see that

(2.23) $$\gamma_1\left(\frac{1}{2}\right) = \gamma_1 - \log^2 2 - 2\gamma \log 2$$

Using Coffey's formula (2.19) we see that

(2.24) $$\sum_{r=1}^{3} \gamma_1\left(\frac{r}{4}\right) = 3\gamma_1 - 8\log^2 2 - 8\gamma \log 2$$

and by definition we also have

$$\sum_{r=1}^{3} \gamma_1\left(\frac{r}{4}\right) = \gamma_1\left(\frac{1}{4}\right) + \gamma_1\left(\frac{1}{2}\right) + \gamma_1\left(\frac{3}{4}\right)$$

(2.25) $$= \gamma_1\left(\frac{1}{4}\right) + \gamma_1\left(\frac{3}{4}\right) + \gamma_1 - \log^2 2 - 2\gamma \log 2$$

where we have used (2.23).

Therefore we get from (2.24) and (2.25)

(2.26) $$\gamma_1\left(\frac{1}{4}\right) + \gamma_1\left(\frac{3}{4}\right) = 2\gamma_1 - 7\log^2 2 - 6\gamma \log 2$$

and using (2.8) we easily see that

(2.27) $$\gamma_1\left(\frac{1}{4}\right) = \frac{1}{2}[2\gamma_1 - 7\log^2 2 - 6\gamma \log 2] - \frac{1}{2}\pi\left[\gamma + 4\log 2 + 3\log \pi - 4\log \Gamma\left(\frac{1}{4}\right)\right]$$

(2.28) $$\gamma_1\left(\frac{3}{4}\right) = \frac{1}{2}[2\gamma_1 - 7\log^2 2 - 6\gamma \log 2] + \frac{1}{2}\pi\left[\gamma + 4\log 2 + 3\log \pi - 4\log \Gamma\left(\frac{1}{4}\right)\right]$$

Assuming that each of the terms $\pi, \gamma, \log 2, \log \pi$ and $\log \Gamma\left(\frac{1}{4}\right)$ has weight equal to one and $\gamma_1$ has weight 2, then the above represents a homogenous equation of weight 2.

Since $\psi\left(\frac{1}{4}\right) = -\gamma - \frac{1}{2}\pi - 3\log 2$ we may write this as

(2.29) $\gamma_1\left(\frac{1}{4}\right) = \gamma\psi\left(\frac{1}{4}\right) + \gamma^2 + \frac{1}{2}[2\gamma_1 - 7\log^2 2] - \frac{1}{2}\pi\left[4\log 2 + 3\log\pi - 4\log\Gamma\left(\frac{1}{4}\right)\right]$

Similarly we have

$$\gamma_1\left(\frac{1}{3}\right) - \gamma_1\left(\frac{2}{3}\right) = -\left[\varsigma'\left(1,\frac{1}{3}\right) - \varsigma'\left(1,\frac{2}{3}\right)\right]$$

and using (2.12) this becomes

(2.30) $\gamma_1\left(\frac{1}{3}\right) - \gamma_1\left(\frac{2}{3}\right) = -\frac{\pi}{2\sqrt{3}}\left[2\gamma - \log 3 + 8\log(2\pi) - 12\log\Gamma\left(\frac{1}{3}\right)\right]$

An alternative proof of this formula is contained below in (4.3).

Equation (2.19) gives us

$$\gamma_1\left(\frac{1}{3}\right) + \gamma_1\left(\frac{2}{3}\right) = 2\gamma_1 - 3\gamma\log 3 - \frac{3}{2}\log^2 3$$

and this then gives us

(2.31) $\gamma_1\left(\frac{1}{3}\right) = \gamma_1 - \frac{3}{2}\gamma\log 3 - \frac{3}{4}\log^2 3 - \frac{\pi}{4\sqrt{3}}\left[2\gamma - \log 3 + 8\log(2\pi) - 12\log\Gamma\left(\frac{1}{3}\right)\right]$

with a similar expression resulting for $\gamma_1\left(\frac{2}{3}\right)$.

In [25] Dilcher defined a generalised polygamma function by

$$\psi_p(x) = -\gamma_p - \frac{\log^p x}{x} - \sum_{n=1}^{\infty}\left[\frac{\log^p(n+x)}{n+x} - \frac{\log^p n}{n}\right]$$

$$= -\gamma_p - \frac{\log^p x}{x} - \sum_{n=1}^{\infty}\left[\frac{\log^p(n+x)}{n+x} - \frac{\log^p(n+1)}{n+1} - \frac{\log^p n}{n} + \frac{\log^p(n+1)}{n+1}\right]$$

$$= -\gamma_p - \frac{\log^p x}{x} - \sum_{n=1}^{\infty}\left[\frac{\log^p(n+x)}{n+x} - \frac{\log^p(n+1)}{n+1}\right]$$

Hence, comparing this with (1.23), we see that

(2.32) $\qquad \psi_p(x) = -\gamma_p(x)$

This relationship was also noted by Coffey in [11].

Dilcher [25] also showed that

$$\psi_1\left(\frac{1}{3}\right) = -\gamma_1 + \frac{1}{2}\left[3\log 3 + \frac{\pi}{\sqrt{3}}\right]\gamma + \frac{3}{4}\log^3 3 + \pi\sqrt{3}\left[\frac{2}{3}\log(2\pi) - \frac{1}{12}\log 3 - \log\Gamma\left(\frac{1}{3}\right)\right]$$

and a little algebra shows that this is the same as (2.31) above.

It was also shown by Dilcher [25] that

$$\psi_k\left(\frac{1}{2}\right) = -\gamma_k\left(\frac{1}{2}\right) = (-1)^{k+1} 2 \frac{\log^{k+1} 2}{k+1} - \gamma_k - 2\sum_{j=0}^{k}\binom{k}{j}(-1)^j \gamma_{k-j}\log^j 2$$

which, as noted by Coffey [11], is the same as (2.11).

As shown by Dilcher in [25] we also have for $k \geq 0$

(2.33) $\qquad (-1)^{k+1}\varsigma_a^{(k)}(1) = \sum_{n=1}^{\infty}(-1)^n \frac{\log^k n}{n} = \sum_{j=0}^{k-1}\binom{k}{j}\gamma_j \log^{k-j} 2 - \frac{\log^{k+1} 2}{k+1}$

where the alternating Riemann zeta function $\varsigma_a(s)$ is defined by

$$\varsigma(s) = \sum_{n=1}^{\infty}\frac{1}{n^s} = \frac{1}{1-2^{-s}}\sum_{n=1}^{\infty}\frac{1}{(2n-1)^s} = \frac{1}{1-2^{-s}}\sum_{n=0}^{\infty}\frac{1}{(2n+1)^s} \qquad , (\text{Re}(s) > 1)$$

$$= \frac{1}{1-2^{1-s}}\sum_{n=1}^{\infty}\frac{(-1)^{n+1}}{n^s} = \frac{1}{1-2^{1-s}}\varsigma_a(s) \qquad , (\text{Re}(s) > 0;\ s \neq 1)$$

For completeness, a further derivation of Dilcher's formula (2.33) is set out below.

We have by direct substitution

$$(s-1)\varsigma_a(s) = 1 - 2^{1-s} + (s-1)\left[\left(1-2^{1-s}\right)\sum_{p=0}^{\infty}\frac{(-1)^p}{p!}\gamma_p(s-1)^p\right]$$

and we now proceed to differentiate this equation this $n+1$ times.

It is easily seen that

$$\frac{d^{n+1}}{ds^{n+1}}\left[1-2^{1-s}\right] = 2^{1-s}(-1)^n \log^{n+1} 2$$

The following result is readily derived

$$\frac{d^{n+1}}{ds^{n+1}}[(s-1)f(s)] = (s-1)f^{(n+1)}(s) + (n+1)f^{(n)}(s)$$

and with $f(s) = \varsigma_a(s)$ we immediately see that

$$\frac{d^{n+1}}{ds^{n+1}}[(s-1)\varsigma_a(s)]\bigg|_{s=1} = (n+1)\varsigma_a^{(n)}(1)$$

and similarly we also have

$$\frac{d^{n+1}}{ds^{n+1}}(s-1)\left[\left(1-2^{1-s}\right)\sum_{p=0}^{\infty}\frac{(-1)^p}{p!}\gamma_p(s-1)^p\right]\bigg|_{s=1} = (n+1)\frac{d^n}{ds^n}\left[\left(1-2^{1-s}\right)\sum_{p=0}^{\infty}\frac{(-1)^p}{p!}\gamma_p(s-1)^p\right]\bigg|_{s=1}$$

Using Leibniz's rule for the derivative of a product we see that

$$\frac{d^n}{ds^n}\left[\left(1-2^{1-s}\right)\sum_{p=0}^{\infty}\frac{(-1)^p}{p!}\gamma_p(s-1)^p\right] = \sum_{j=0}^{n}\binom{n}{j}\frac{d^j}{ds^j}\left[\sum_{p=0}^{\infty}\frac{(-1)^p}{p!}\gamma_p(s-1)^p\right]\frac{d^{n-j}}{ds^{n-j}}\left[1-2^{1-s}\right]$$

and we have

$$\frac{d^j}{ds^j}\left[\sum_{p=0}^{\infty}\frac{(-1)^p}{p!}\gamma_p(s-1)^p\right] = \sum_{p=0}^{\infty}\frac{(-1)^p}{p!}\gamma_p p(p-1)...(p+1-j)(s-1)^{p-j}$$

When $s=1$ the only non-vanishing term is when $p=j$ and hence we obtain

$$\frac{d^j}{ds^j}\left[\sum_{p=0}^{\infty}\frac{(-1)^p}{p!}\gamma_p(s-1)^p\right]\bigg|_{s=1} = (-1)^j \gamma_j$$

Using the fact that
$$\frac{d^{n-j}}{ds^{n-j}}\left[1-2^{1-s}\right]\bigg|_{s=1} = (-1)^{n-j+1} \log^{n-j} 2$$

$$= 0 \text{ when } j = n$$

we obtain

$$(n+1)\varsigma_a^{(n)}(1) = (-1)^n \log^{n+1} 2 + (n+1)\sum_{j=0}^{n-1}\binom{n}{j}(-1)^j \gamma_j (-1)^{n-j+1} \log^{n-j} 2$$

which may be written as

$$(-1)^{n+1}\varsigma_a^{(n)}(1) = \sum_{j=0}^{n-1}\binom{n}{j}\gamma_j \log^{n-j} 2 - \frac{\log^{n+1} 2}{n+1}$$

The $n$ th derivative of the alternating Riemann zeta function $\varsigma_a(s)$ is

$$(-1)^n \varsigma_a^{(n)}(1) = \sum_{k=1}^{\infty}(-1)^{k-1}\frac{\log^n k}{k}$$

and accordingly we obtain

(2.34) $$-\varsigma_a^{(1)}(1) = \sum_{k=1}^{\infty}(-1)^{k-1}\frac{\log k}{k} = \frac{1}{2}\log^2 2 - \gamma \log 2$$

(2.35) $$\varsigma_a^{(2)}(1) = \sum_{k=1}^{\infty}(-1)^{k-1}\frac{\log^2 k}{k} = \frac{1}{3}\log^3 2 - \gamma \log^2 2 - 2\gamma_1 \log 2$$

The above formula (2.33) was in fact first reported by Briggs and Chowla [9] in 1955 where they showed that

(2.36) $$\varsigma_a^{(k)}(1) = k!\sum_{r=1}^{k+1}\frac{(-1)^{r+1}\log^r 2}{r!}A_{k-r}$$

with $A_n = \frac{(-1)^n}{n!}\gamma_n$ and $A_{-1} = 1$.

It was rediscovered by Liang and Todd [36] in 1972 who used it to numerically compute the values of the first 20 Stieltjes constants (see Dilcher's paper [25] for more details).

Equation (2.34) was posed as a problem by Klamkin in 1954 and is closely related to an earlier problem posed by Sandham [40] in 1950.

Reference should also be made to the paper by Collins [18] entitled "The role of Bell polynomials in integration" where the integrals $\int_0^\infty \frac{\log^n x}{e^x+1} dx$ are evaluated in terms of the Stieltjes constants and the complete Bell polynomials. In this regard, it may of course be noted that

$$\int_0^\infty \frac{\log^n x}{e^x+1} dx = \frac{d^n}{ds^n} \int_0^\infty \frac{x^{s-1}}{e^x+1} dx \bigg|_{s=1} = \frac{d^n}{ds^n} [\Gamma(s)\varsigma_a(s)] \bigg|_{s=1}$$

$$= \sum_{k=0}^n \binom{n}{k} \varsigma_a^{(k)}(1) \Gamma^{(n-k)}(1)$$

and substituting Dilcher's formula (2.33)

$$(-1)^{k+1} \varsigma_a^{(k)}(1) = \sum_{j=0}^{k-1} \binom{k}{j} \gamma_j \log^{k-j} 2 - \frac{\log^{k+1} 2}{k+1}$$

we then determine that

(2.37) $$\int_0^\infty \frac{\log^n x}{e^x+1} dx = \sum_{k=0}^n \binom{n}{k} (-1)^k \left[ \frac{\log^{k+1} 2}{k+1} - \sum_{j=0}^{k-1} \binom{k}{j} \gamma_j \log^{k-j} 2 \right] \Gamma^{(n-k)}(1)$$

as previously determined by Collins [18].

Making the substitution $x = e^{-u}$ we get

$$\int_0^\infty \frac{\log u}{e^u+1} du = \int_0^1 \frac{\log \log(1/x)}{1+x} dx$$

Adamchik [2] has considered logarithmic integrals of this type and has shown that

$$\int_0^1 \frac{x^{p-1}}{1+x^n} \log\log\left(\frac{1}{x}\right) dx = \frac{\gamma + \log(2n)}{2n} \left[ \psi\left(\frac{p}{2n}\right) - \psi\left(\frac{n+p}{2n}\right) \right] + \frac{1}{2n}\left[ \varsigma'\left(1, \frac{p}{2n}\right) - \varsigma'\left(1, \frac{n+p}{2n}\right) \right]$$

and in particular we have another proof of the above integral

$$\int_0^1 \frac{1}{1+x} \log\log\left(\frac{1}{x}\right) dx = -\frac{\log^2 2}{2}$$

We may express $\Gamma^{(m)}(x)$ in terms of $\psi(x)$ and the Hurwitz zeta functions as arguments of the (exponential) complete Bell polynomials.

The (exponential) complete Bell polynomials may be defined by $Y_0 = 1$ and for $n \geq 1$

(2.38) $$Y_n(x_1,...,x_n) = \sum_{\pi(n)} \frac{n!}{k_1! k_2! ... k_n!} \left(\frac{x_1}{1!}\right)^{k_1} \left(\frac{x_2}{2!}\right)^{k_2} ... \left(\frac{x_n}{n!}\right)^{k_n}$$

where the sum is taken over all partitions $\pi(n)$ of $n$, i.e. over all sets of integers $k_j$ such that
$$k_1 + 2k_2 + 3k_3 + ... + nk_n = n$$

The complete Bell polynomials have integer coefficients and the first six are set out below [19, p.307]

(2.39) $$Y_1(x_1) = x_1$$

$$Y_2(x_1, x_2) = x_1^2 + x_2$$

$$Y_3(x_1, x_2, x_3) = x_1^3 + 3x_1 x_2 + x_3$$

$$Y_4(x_1, x_2, x_3, x_4) = x_1^4 + 6x_1^2 x_2 + 4x_1 x_3 + 3x_2^2 + x_4$$

$$Y_5(x_1, x_2, x_3, x_4, x_5) = x_1^5 + 10x_1^3 x_2 + 10x_1^2 x_3 + 15x_1 x_2^2 + 5x_1 x_4 + 10x_2 x_3 + x_5$$

$$Y_6(x_1, x_2, x_3, x_4, x_5, x_6) = x_1^6 + 6x_1 x_5 + 15x_2 x_4 + 10x_3^2 + 15x_1^2 x_4 + 15x_2^3 + 60x_1 x_2 x_3$$
$$+ 20x_1^3 x_3 + 45x_1^2 x_2^2 + 15x_1^4 x_1 + x_6$$

The complete Bell polynomials are also given by the exponential generating function reported in Comtet's book [19, p.134]

(2.40) $$\exp\left(\sum_{j=1}^{\infty} x_j \frac{t^j}{j!}\right) = 1 + \sum_{n=1}^{\infty} Y_n(x_1,...,x_n) \frac{t^n}{n!} = \sum_{n=0}^{\infty} Y_n(x_1,...,x_n) \frac{t^n}{n!}$$

In particular, Coffey [14a] notes that

(2.41) $$\Gamma^{(m)}(1) = Y_m(-\gamma, x_1,..., x_{m-1})$$

where $x_p = (-1)^{p+1} p! \varsigma(p+1)$.

Values of $\Gamma^{(m)}(1)$ are reported in [42, p.265] for $m \leq 10$ and the first three are

$$\Gamma^{(1)}(1) = -\gamma$$

$$\Gamma^{(2)}(1) = \varsigma(2) + \gamma^2$$

$$\Gamma^{(3)}(1) = -2\varsigma(3) - 3\gamma\varsigma(2) - \gamma^3$$

□

We note from [21, Eq. (4.3.231)] that for $t > 0$

$$\int_1^t \gamma_p(x) dx = \frac{(-1)^{p+1}}{p+1} \left[ \varsigma^{(p+1)}(0,t) - \varsigma^{(p+1)}(0) \right]$$

and hence integrating (2.18.1) results in

$$\int_1^t \sum_{r=0}^{q-1} \gamma_p\left(\frac{r+x}{q}\right) dx$$

$$= q(-1)^p \frac{\log^{p+1} q}{p+1}(t-1) + q \sum_{k=0}^p \binom{p}{k}(-1)^k \frac{(-1)^{p-k+1}}{p-k+1} \left[ \varsigma^{(p-k+1)}(0,t) - \varsigma^{(p-k+1)}(0) \right] \log^k q$$

We note that

$$\int_1^t \gamma_p\left(\frac{r+x}{q}\right) dx = q \int_{(r+1)/q}^{(r+t)/q} \gamma_p(u) du = q \int_1^{(r+t)/q} \gamma_p(u) du - q \int_1^{(r+1)/q} \gamma_p(u) du$$

$$= \frac{(-1)^{p+1} q}{p+1} \left[ \varsigma^{(p+1)}\left(0, \frac{r+t}{q}\right) - \varsigma^{(p+1)}\left(0, \frac{r+1}{q}\right) \right]$$

and therefore determine that

$$\sum_{r=0}^{q-1} \gamma_p \frac{(-1)^{p+1} q}{p+1} \left[ \varsigma^{(p+1)}\left(0, \frac{r+t}{q}\right) - \varsigma^{(p+1)}\left(0, \frac{r+1}{q}\right) \right]$$

$$= q(-1)^p \frac{\log^{p+1} q}{p+1}(t-1) + q \sum_{k=0}^{p} \binom{p}{k}(-1)^k \frac{(-1)^{p-k+1}}{p-k+1}\left[\varsigma^{(p-k+1)}(0,t) - \varsigma^{(p-k+1)}(0)\right]\log^k q$$

□

In passing, using the following expression

$$\varsigma_a(s) = \frac{1-2^{1-s}}{s-1} + \left(1-2^{1-s}\right)\sum_{p=0}^{\infty} \frac{(-1)^p}{p!}\gamma_p(s-1)^p$$

we see that

$$\varsigma_a(1) = \lim_{s \to 1} \frac{1-2^{1-s}}{s-1}$$

and using L'Hôpital's rule this becomes the well-known result

$$\varsigma_a(1) = \lim_{s \to 1}\left[2^{1-s}\log 2\right] = \log 2$$

**Lemma:**

From Gradshteyn and Ryzhik [28, p.35] we have the well-known formula (which is easily derived by using the complex number representation $e^{ijx}$ for the trigonometric functions)

(2.42) $$\sum_{j=1}^{n} \sin jx = \sin\frac{(n+1)x}{2}\sin\frac{nx}{2}\operatorname{cosec}\frac{x}{2}$$

and we immediately see that

(2.43) $$\sum_{j=1}^{q} \sin\left(\frac{2\pi jp}{q}\right) = 0$$

as required in (2.6) above.

From [28, p.35] we also have

$$\sum_{j=1}^{n} \cos jx = \cos\frac{(n+1)x}{2}\sin\frac{nx}{2}\operatorname{cosec}\frac{x}{2}$$

and differentiation results in

$$\sum_{j=1}^{n} j \sin jx = \cos \frac{(n+1)x}{2} \sin \frac{nx}{2} \operatorname{cosec} \frac{x}{2} \left[ \frac{(n+1)}{2} \tan \frac{(n+1)x}{2} - \frac{n}{2} \cot \frac{nx}{2} + \frac{1}{2} \cot \frac{x}{2} \right]$$

$$= \frac{(n+1)}{2} \sin \frac{(n+1)x}{2} \sin \frac{nx}{2} \operatorname{cosec} \frac{x}{2} - \frac{n}{2} \cos \frac{(n+1)x}{2} \cos \frac{nx}{2} \operatorname{cosec} \frac{x}{2}$$

$$+ \frac{1}{2} \cot \frac{x}{2} \operatorname{cosec} \frac{x}{2} \cos \frac{(n+1)x}{2} \sin \frac{nx}{2}$$

Since the terms involving $\sin \frac{nx}{2}$ will vanish when $\frac{nx}{2} = p\pi$ we therefore obtain

$$\sum_{j=1}^{q} j \sin \left( \frac{2\pi jp}{q} \right) = -\frac{q}{2} \cos \frac{(q+1)p\pi}{q} \cos p\pi \operatorname{cosec} \frac{p\pi}{q}$$

$$= -\frac{q}{2} \left[ \cos p\pi \cos \frac{p\pi}{q} - \sin p\pi \sin \frac{p\pi}{q} \right] \cos p\pi \operatorname{cosec} \frac{p\pi}{q}$$

and we get the result required in (2.11) above

(2.44) $$\sum_{j=1}^{q} j \sin \left( \frac{2\pi jp}{q} \right) = -\frac{q}{2} \cot \left( \frac{\pi p}{q} \right)$$

It is easily seen that differentiating (2.42) twice will provide an expression for $\sum_{j=1}^{q} j^2 \sin \left( \frac{2\pi jp}{q} \right)$.

The following lemma is required in (5.1) below. We showed in (2.7) that

$$\psi(1-s) + \frac{\pi}{2} \tan \left( \frac{\pi s}{2} \right) = \psi(s) + \frac{\pi}{2} \cot \left( \frac{\pi s}{2} \right)$$

and differentiation results in

$$-\psi'(1-s) + \frac{\pi^2}{4} \sec^2 \left( \frac{\pi s}{2} \right) = \psi'(s) - \frac{\pi^2}{4} \operatorname{cosec}^2 \left( \frac{\pi s}{2} \right)$$

and as $s \to 1$ this becomes

$$= \psi'(1) - \frac{\pi^2}{4}$$

Since $\psi'(1) = \varsigma(2) = \frac{\pi^2}{6}$ we have the limit

(2.45) $$\lim_{s \to 1}\left[-\psi'(1-s) + \frac{\pi^2}{4}\sec^2\left(\frac{\pi s}{2}\right)\right] = -\frac{\pi^2}{12}$$

$\square$

Referring back to (2.1.1)

$$\varsigma\left(s,\frac{p}{q}\right) - \varsigma\left(s,1-\frac{p}{q}\right) = 4\Gamma(1-s)(2\pi q)^{s-1}\cos\left(\frac{\pi s}{2}\right)\sum_{j=1}^{q}\sin\left(\frac{2\pi jp}{q}\right)\varsigma\left(1-s,\frac{j}{q}\right)$$

we see that

$$\varsigma\left(1,\frac{p}{q}\right) - \varsigma\left(1,1-\frac{p}{q}\right) = 2\pi\sum_{j=1}^{q}\sin\left(\frac{2\pi jp}{q}\right)\varsigma\left(0,\frac{j}{q}\right)$$

$$= 2\pi\sum_{j=1}^{q}\left(\frac{1}{2} - \frac{j}{q}\right)\sin\left(\frac{2\pi jp}{q}\right)$$

and using (2.6) this becomes

$$= -\frac{2\pi}{q}\sum_{j=1}^{q}j\sin\left(\frac{2\pi jp}{q}\right)$$

Now using (2.11)

$$\sum_{j=1}^{q}j\sin\left(\frac{2\pi jp}{q}\right) = -\frac{q}{2}\cot\left(\frac{\pi p}{q}\right) \quad \text{for } p < q$$

we obtain

$$\varsigma\left(1,\frac{p}{q}\right) - \varsigma\left(1,1-\frac{p}{q}\right) = \pi\cot\left(\frac{\pi p}{q}\right)$$

This then becomes the well-known identity [42, p.14]

$$\psi\left(\frac{p}{q}\right) - \psi\left(1-\frac{p}{q}\right) = -\pi\cot\left(\frac{\pi p}{q}\right)$$

$\square$

From (2.17) we see that

$$\gamma_0(x) - \gamma_0(y) = \lim_{s \to 1}[\varsigma(s,x) - \varsigma(s,y)]$$

and using the Hasse formula (1.1) this becomes

$$\gamma_0(x) - \gamma_0(y) = \lim_{s \to 1} \frac{1}{s-1}\left[\sum_{n=0}^{\infty}\frac{1}{n+1}\sum_{k=0}^{n}\binom{n}{k}\frac{(-1)^k}{(x+k)^{s-1}} - \sum_{n=0}^{\infty}\frac{1}{n+1}\sum_{k=0}^{n}\binom{n}{k}\frac{(-1)^k}{(y+k)^{s-1}}\right]$$

Applying L'Hôpital's rule we obtain

$$\gamma_0(x) - \gamma_0(y) = \sum_{n=0}^{\infty}\frac{1}{n+1}\sum_{k=0}^{n}\binom{n}{k}(-1)^k \log(y+k) - \sum_{n=0}^{\infty}\frac{1}{n+1}\sum_{k=0}^{n}\binom{n}{k}(-1)^k \log(x+k)$$

which concurs with (1.4). Differentiating the Hasse formula (1.1) gives us

$$(s-1)\varsigma'(s,x) + \varsigma(s,x) = \sum_{n=0}^{\infty}\frac{1}{n+1}\sum_{k=0}^{n}\binom{n}{k}\frac{(-1)^k \log(x+k)}{(x+k)^{s-1}}$$

and hence we have

$$\varsigma'(s,x) - \varsigma'(s,y) = -\frac{1}{s-1}\left[\sum_{n=0}^{\infty}\frac{1}{n+1}\sum_{k=0}^{n}\binom{n}{k}\frac{(-1)^k \log(x+k)}{(x+k)^{s-1}} - \sum_{n=0}^{\infty}\frac{1}{n+1}\sum_{k=0}^{n}\binom{n}{k}\frac{(-1)^k \log(y+k)}{(y+k)^{s-1}}\right]$$

$$-\frac{1}{s-1}[\varsigma(s,x) - \varsigma(s,y)]$$

Applying L'Hôpital's rule again we obtain

$$\lim_{s \to 1}[\varsigma'(s,x) - \varsigma'(s,y)] = \sum_{n=0}^{\infty}\frac{1}{n+1}\sum_{k=0}^{n}\binom{n}{k}(-1)^k \log^2(x+k) - \sum_{n=0}^{\infty}\frac{1}{n+1}\sum_{k=0}^{n}\binom{n}{k}(-1)^k \log^2(y+k)$$

$$-\lim_{s \to 1}[\varsigma'(s,x) - \varsigma'(s,y)]$$

and thus

$$\lim_{s \to 1}[\varsigma'(s,x) - \varsigma'(s,y)] = \frac{1}{2}\left[\sum_{n=0}^{\infty}\frac{1}{n+1}\sum_{k=0}^{n}\binom{n}{k}(-1)^k \log^2(x+k) - \sum_{n=0}^{\infty}\frac{1}{n+1}\sum_{k=0}^{n}\binom{n}{k}(-1)^k \log^2(y+k)\right]$$

Reference to (2.17) then results in

$$\gamma_1(x) - \gamma_1(y) = -\frac{1}{2} \left[ \sum_{n=0}^{\infty} \frac{1}{n+1} \sum_{k=0}^{n} \binom{n}{k} (-1)^k \log^2(x+k) - \sum_{n=0}^{\infty} \frac{1}{n+1} \sum_{k=0}^{n} \binom{n}{k} (-1)^k \log^2(y+k) \right]$$

which concurs with (1.14). In a similar way it is easily shown that

$$\gamma_p(x) - \gamma_p(y) = -\frac{1}{p+1} \left[ \sum_{n=0}^{\infty} \frac{1}{n+1} \sum_{k=0}^{n} \binom{n}{k} (-1)^k \log^{p+1}(x+k) - \sum_{n=0}^{\infty} \frac{1}{n+1} \sum_{k=0}^{n} \binom{n}{k} (-1)^k \log^{p+1}(y+k) \right]$$

which concurs with (1.16).

## 3. A CONNECTION WITH ADAMCHIK'S LOGARITHMIC INTEGRALS

Another derivation of (2.27) is set out below.

Using the integral definition of the gamma function for $s > 0$

$$\Gamma(s) = \int_0^{\infty} e^{-x} x^{s-1} \, dt$$

and making the substitution $x = ut$ (where $u > 0$) we obtain for $u, s > 0$

$$\int_0^{\infty} e^{-ut} t^{s-1} \, dt = \frac{\Gamma(s)}{u^s}$$

and using the binomial theorem we get for $n \geq 0$

(3.1) $$\int_0^{\infty} e^{-ut} (1-e^{-t})^n t^{s-1} \, dt = \Gamma(s) \sum_{k=0}^{n} \binom{n}{k} \frac{(-1)^k}{(k+u)^s}$$

Making the summation gives us

$$\sum_{n=0}^{\infty} \frac{1}{n+1} \int_0^{\infty} e^{-ut} (1-e^{-t})^n t^{s-1} \, dt = \sum_{n=0}^{\infty} \int_0^{\infty} e^{-(n+1)x} dx \int_0^{\infty} e^{-ut} (1-e^{-t})^n t^{s-1} \, dt$$

and this may be expressed as

$$\sum_{n=0}^{\infty} \int_0^{\infty} e^{-(n+1)x} dx \int_0^{\infty} e^{-ut} (1-e^{-t})^n t^{s-1} \, dt = \sum_{n=0}^{\infty} \int_0^{\infty} \int_0^{\infty} e^{-ut} e^{-x} \left[ e^{-x} (1-e^{-t}) \right]^n t^{s-1} \, dt \, dx$$

$$= \int_0^\infty \int_0^\infty \frac{e^{-ut}e^{-x}t^{s-1}}{1-e^{-x}[1-e^{-t}]} dt\, dx$$

We have

$$\int_0^\infty \frac{e^{-x}}{1-e^{-x}[1-e^{-t}]} dx = \frac{1}{1-e^{-t}} \log(1-e^{-x}[1-e^{-t}])\Big|_0^\infty = \frac{t}{1-e^{-t}}$$

and we therefore see that

$$\int_0^\infty \int_0^\infty \frac{e^{-ut}e^{-x}t^{s-1}}{1-e^{-x}[1-e^{-t}]} dt\, dx = \int_0^\infty \frac{e^{-ut}t^s}{1-e^{-t}} dt = \int_0^\infty \frac{e^{-(u-1)t}t^s}{e^t-1} dt$$

We note from [42, p.92] that for $\text{Re}(s)>0$ and $\text{Re}(u)>0$

(3.2) $$\int_0^\infty \frac{e^{-(u-1)t}t^s}{e^t-1} dt = \Gamma(s+1)\varsigma(s+1,u)$$

and with $u=1$ we obtain the familiar formula

(3.3) $$\int_0^\infty \frac{t^{s-1}}{e^t-1} dt = \Gamma(s)\varsigma(s,1) = \Gamma(s)\varsigma(s)$$

We therefore have

(3.4) $$\int_0^\infty \frac{e^{-(u-1)t}t^s}{e^t-1} dt = \Gamma(s)\sum_{n=0}^\infty \frac{1}{n+1}\sum_{k=0}^n \binom{n}{k}\frac{(-1)^k}{(k+u)^s}$$

and comparing this with (3.2) gives us another proof of the Hasse identity (1.1) for $\text{Re}(s)>0$

$$\varsigma(s+1,u) = \frac{1}{s}\sum_{n=0}^\infty \frac{1}{n+1}\sum_{k=0}^n \binom{n}{k}\frac{(-1)^k}{(k+u)^s}$$

Integrating (3.2) with respect to $u$ we see that

$$\int_1^x du \int_0^\infty \frac{e^{-(u-1)t}t^s}{e^t-1} dt = \int_0^\infty \frac{[1-e^{-(x-1)t}]t^{s-1}}{e^t-1} = \Gamma(s)[\varsigma(s)-\varsigma(s,x)]$$

We also note that

$$\int_0^\infty \frac{e^{-(u-1)t} t^s}{e^t - 1} dt = \Gamma(s+1)\varsigma(s+1,u) = \Gamma(s)s\varsigma(s+1,u) = -\Gamma(s)\frac{\partial}{\partial u}\varsigma(s,u)$$

$$= \Gamma(s)\sum_{n=0}^\infty \frac{1}{n+1}\sum_{k=0}^n \binom{n}{k}\frac{(-1)^k}{(k+u)^s}$$

We have

$$\int_0^\infty \frac{e^{-(u-1)t} t^s}{\Gamma(s)[e^t - 1]} dt = \sum_{n=0}^\infty \frac{1}{n+1}\sum_{k=0}^n \binom{n}{k}\frac{(-1)^k}{(k+u)^s}$$

and differentiating this with respect to $s$ results in

$$\int_0^\infty \frac{e^{-(u-1)t} t^s [\log t - \psi(s)]}{\Gamma(s)[e^t - 1]} dt = -\sum_{n=0}^\infty \frac{1}{n+1}\sum_{k=0}^n \binom{n}{k}\frac{(-1)^k \log(k+u)}{(k+u)^s}$$

Hence we have

(3.5) $$\int_0^\infty \frac{e^{-(u-1)t} t^s \log t}{e^t - 1} dt = \psi(s)\Gamma(s+1)\varsigma(s+1,u) - \Gamma(s)\sum_{n=0}^\infty \frac{1}{n+1}\sum_{k=0}^n \binom{n}{k}\frac{(-1)^k \log(k+u)}{(k+u)^s}$$

With $u = 1$ in (3.5) we get

$$\int_0^\infty \frac{t^s \log t}{e^t - 1} dt = \psi(s)\Gamma(s+1)\varsigma(s+1) - \Gamma(s)\sum_{n=0}^\infty \frac{1}{n+1}\sum_{k=0}^n \binom{n}{k}\frac{(-1)^k \log(k+1)}{(k+1)^s}$$

which is equivalent to

$$\int_0^\infty \frac{t^s \log t}{e^t - 1} dt = \frac{d}{ds}[\Gamma(s+1)\varsigma(s+1)]$$

With $s = 1$ in (3.5) we get

$$\int_0^\infty \frac{e^{-(u-1)t} t \log t}{e^t - 1} dt = -\gamma\varsigma(2,u) - \sum_{n=0}^\infty \frac{1}{n+1}\sum_{k=0}^n \binom{n}{k}\frac{(-1)^k \log(k+u)}{k+u}$$

and using (1.18) this becomes

$$\int_0^\infty \frac{e^{-(u-1)t} t \log t}{e^t - 1} dt = -\gamma\psi'(u) - \sum_{n=0}^\infty \frac{1}{n+1}\sum_{k=0}^n \binom{n}{k}\frac{(-1)^k \log(k+u)}{k+u}$$

Reference to (1.14) then shows that

$$(3.6) \quad \int_0^\infty \frac{e^{-(u-1)t} t \log t}{e^t - 1} dt = \gamma_1'(u) - \gamma \psi'(u) = \gamma_1'(u) + \gamma \gamma_0'(u)$$

Integrating this with respect to $u$ we get

$$\int_1^x du \int_0^\infty \frac{e^{-(u-1)t} t \log t}{e^t - 1} dt = \int_0^\infty \frac{[1 - e^{-(x-1)t}] \log t}{e^t - 1} dt = \gamma_1(x) - \gamma_1 - \gamma \psi(x) - \gamma^2$$

The substitution $y = e^{-t}$ gives us

$$(3.7) \quad \int_0^1 \frac{1 - y^{x-1}}{1 - y} \log \log(1/y) \, dy = \gamma_1(x) - \gamma_1 - \gamma \psi(x) - \gamma^2$$

and with $x = 1/2$ this becomes

$$(3.8) \quad \int_0^1 \frac{1 - 1/\sqrt{y}}{1 - y} \log \log(1/y) \, dy = -\int_0^1 \frac{1/\sqrt{y}}{1 + \sqrt{y}} \log \log(1/y) \, dy = \gamma_1\left(\frac{1}{2}\right) - \gamma_1 - \gamma \psi\left(\frac{1}{2}\right) - \gamma^2$$

We now recall Adamchik's result (see reference [2] and also Roach's paper [39])

(3.9)

$$\int_0^1 \frac{x^{p-1}}{1 + x^n} \log \log\left(\frac{1}{x}\right) dx = \frac{\gamma + \log(2n)}{2n}\left[\psi\left(\frac{p}{2n}\right) - \psi\left(\frac{n+p}{2n}\right)\right] + \frac{1}{2n}\left[\varsigma'\left(1, \frac{p}{2n}\right) - \varsigma'\left(1, \frac{n+p}{2n}\right)\right]$$

where, notwithstanding the suggestive notation, neither $p$ nor $n$ need be integers.

With $p = n$ we get

$$(3.10) \quad \int_0^1 \frac{x^{n-1}}{1 + x^n} \log \log\left(\frac{1}{x}\right) dx = \frac{\gamma + \log(2n)}{2n}\left[\psi\left(\frac{1}{2}\right) - \psi(1)\right] + \frac{1}{2n}\left[\varsigma'\left(1, \frac{1}{2}\right) - \varsigma'(1,1)\right]$$

From Adamchik's paper [2], we see that

$$(3.11) \quad \varsigma'\left(1, \frac{1}{2}\right) - \varsigma'(1,1) = \log^2 2 + 2\gamma \log 2$$

and hence, using $\psi\left(\frac{1}{2}\right) = -\gamma - 2\log 2$, we have

(3.12) $$\int_0^1 \frac{x^{n-1}}{1+x^n} \log\log\left(\frac{1}{x}\right) dx = -\frac{1}{2n}\left[\log^2 2 + 2\log 2 \log n\right]$$

With $p = n = 1/2$ we get

(3.13) $$\int_0^1 \frac{1/\sqrt{x}}{1+\sqrt{x}} \log\log\left(\frac{1}{x}\right) dx = \log^2 2$$

We then see from (3.8) and (3.13) that

(3.14) $$\gamma_1\left(\frac{1}{2}\right) = \gamma_1 - \log^2 2 - 2\gamma \log 2$$

and this is the same as the result obtained by Coffey [11] in (2.22) above.

From (3.7) we have

(3.15) $$\int_0^1 \frac{1-x^{-3/4}}{1-x} \log\log(1/x) dx = \gamma_1\left(\frac{1}{4}\right) - \gamma_1 - \gamma\psi\left(\frac{1}{4}\right) - \gamma^2$$

and using [42, p.22]

$$\psi\left(\frac{1}{4}\right) = -\gamma - \frac{1}{2}\pi - 3\log 2$$

we get

(3.16) $$\int_0^1 \frac{1-x^{-3/4}}{1-x} \log\log(1/x) dx = \gamma_1\left(\frac{1}{4}\right) - \gamma_1 + \frac{1}{2}\gamma\pi + 3\gamma \log 2$$

We have

$$\frac{1-x^{-3/4}}{1-x} = -\frac{x^{-3/4}}{2(x^{1/4}+1)} - \frac{x^{-1/2}}{2(x^{1/2}+1)} - \frac{x^{-3/4}}{2(x^{1/2}+1)}$$

Letting $p = 1/4$ and $n = 1/2$ in (3.9) we get

$$\int_0^1 \frac{x^{-3/4}}{1+x^{1/2}} \log\log\left(\frac{1}{x}\right) dx = \gamma\left[\psi\left(\frac{1}{4}\right) - \psi\left(\frac{3}{4}\right)\right] + \left[\varsigma'\left(1,\frac{1}{4}\right) - \varsigma'\left(1,\frac{3}{4}\right)\right]$$

and with $p = 1/2$ and $n = 1/2$ in (3.9) we see that

$$\int_0^1 \frac{x^{-1/2}}{1+x^{1/2}} \log\log\left(\frac{1}{x}\right) dx = \log^2 2$$

With $p = n = 1/4$ we get

$$\int_0^1 \frac{x^{-3/4}}{1+x^{1/4}} \log\log\left(\frac{1}{x}\right) dx = 6\log^2 2$$

We therefore obtain

$$\int_0^1 \frac{1-x^{-3/4}}{1-x} \log\log(1/x)\, dx = -\frac{7}{2}\log^2 2 - \frac{1}{2}\gamma\left[\psi\left(\frac{1}{4}\right) - \psi\left(\frac{3}{4}\right)\right] - \frac{1}{2}\left[\varsigma'\left(1,\frac{1}{4}\right) - \varsigma'\left(1,\frac{3}{4}\right)\right]$$

We then refer to (2.17) and see that

$$\varsigma'\left(1,\frac{1}{4}\right) - \varsigma'\left(1,\frac{3}{4}\right) = -\left[\gamma_1\left(\frac{1}{4}\right) - \gamma_1\left(\frac{3}{4}\right)\right]$$

and using [42, p.22]

$$\psi\left(\frac{1}{4}\right) = -\gamma - \frac{1}{2}\pi - 3\log 2 \qquad \psi\left(\frac{3}{4}\right) = -\gamma + \frac{1}{2}\pi - 3\log 2$$

we get

(3.17) $$\int_0^1 \frac{1-x^{-3/4}}{1-x} \log\log(1/x)\, dx = -\frac{7}{2}\log^2 2 + \frac{1}{2}\gamma\pi + \frac{1}{2}\left[\gamma_1\left(\frac{1}{4}\right) - \gamma_1\left(\frac{3}{4}\right)\right]$$

Using (2.18)

$$\gamma_1\left(\frac{1}{4}\right) - \gamma_1\left(\frac{3}{4}\right) = -\pi\left[\gamma + 4\log 2 + 3\log \pi - 4\log \Gamma\left(\frac{1}{4}\right)\right]$$

this then becomes

$$\int_0^1 \frac{1-x^{-3/4}}{1-x}\log\log(1/x)\,dx = -\frac{7}{2}\log^2 2 + \frac{1}{2}\gamma\pi - \frac{1}{2}\pi\left[\gamma + 4\log 2 + 3\log\pi - 4\log\Gamma\left(\frac{1}{4}\right)\right]$$

and equating this with (3.16) we get

(3.18)

$$\gamma_1\left(\frac{1}{4}\right) - \gamma_1 + \frac{1}{2}\gamma\pi + 3\gamma\log 2 = -\frac{7}{2}\log^2 2 + \frac{1}{2}\gamma\pi - \frac{1}{2}\pi\left[\gamma + 4\log 2 + 3\log\pi - 4\log\Gamma\left(\frac{1}{4}\right)\right]$$

This is in fact another proof of equation (2.27).

With $x = 2$ in (3.7) we obtain

$$\int_0^1 \log\log(1/y)\,dy = \gamma_1(2) - \gamma_1 - \gamma\psi(2) - \gamma^2$$

and since

$$\int_0^1 \log\log(1/y)\,dy = \Gamma'(1) = -\gamma \qquad \psi(2) = 1 - \gamma$$

we see that

(3.19) $\qquad \gamma_1(2) = \gamma_1$

which is a particular case of (2.17.1).

## 4. YET ANOTHER DERIVATION

In 1999, Shail [41] showed that

(4.1) $\qquad \sum_{n=0}^{\infty}\left[\frac{\log(3n+1)}{3n+1} - \frac{\log(3n+2)}{3n+2}\right] = \frac{\pi}{\sqrt{3}}\left[\log\frac{\Gamma(1/3)}{\Gamma(2/3)} - \frac{1}{3}(\gamma + \log(2\pi))\right]$

Using Euler's reflection formula for the gamma function we see that

$$\Gamma(2/3) = \frac{2\pi}{\sqrt{3}\,\Gamma(1/3)}$$

and hence Shail's formula may be written as

(4.2) $$\sum_{n=0}^{\infty}\left[\frac{\log(3n+1)}{3n+1}-\frac{\log(3n+2)}{3n+2}\right]=\frac{\pi}{\sqrt{3}}\left[2\log\Gamma(1/3)+\frac{1}{2}\log 3-\frac{1}{3}(\gamma+4\log(2\pi))\right]$$

We have from (1.18)

$$\gamma_1(x)-\gamma_1(y)=\sum_{n=0}^{\infty}\left[\frac{\log(n+x)}{n+x}-\frac{\log(n+y)}{n+y}\right]$$

and therefore we get

$$\gamma_1\left(\frac{1}{3}\right)-\gamma_1\left(\frac{2}{3}\right)=3\sum_{n=0}^{\infty}\left[\frac{\log(3n+1)}{3n+1}-\frac{\log(3n+2)}{3n+2}\right]-3\log 3\sum_{n=0}^{\infty}\left[\frac{1}{3n+1}-\frac{1}{3n+2}\right]$$

Using the following expression for the digamma function (1.24)

$$\psi(x)=-\gamma-\frac{1}{x}+\sum_{n=1}^{\infty}\left[\frac{1}{n}-\frac{1}{n+x}\right]$$

and we may deduce that

$$\psi(x)-\psi(y)=-\left[\frac{1}{x}-\frac{1}{y}\right]-\sum_{n=1}^{\infty}\left[\frac{1}{n+x}-\frac{1}{n+y}\right]$$

Hence we have

$$\psi\left(\frac{1}{3}\right)-\psi\left(\frac{2}{3}\right)=-\frac{3}{2}-3\sum_{n=1}^{\infty}\left[\frac{1}{3n+1}-\frac{1}{3n+2}\right]=-3\sum_{n=0}^{\infty}\left[\frac{1}{3n+1}-\frac{1}{3n+2}\right]$$

It is well known that [42, p.14]

$$\psi(x)-\psi(1-x)=-\pi\cot\pi x$$

and thus we have

$$\psi\left(\frac{1}{3}\right)-\psi\left(\frac{2}{3}\right)=-\frac{\pi}{\sqrt{3}}$$

This then gives us

$$\gamma_1\left(\frac{1}{3}\right)-\gamma_1\left(\frac{2}{3}\right)=3\sum_{n=0}^{\infty}\left[\frac{\log(3n+1)}{3n+1}-\frac{\log(3n+2)}{3n+2}\right]-\frac{\pi}{\sqrt{3}}\log 3$$

and using (4.2) this becomes

$$(4.3) \quad \gamma_1\left(\frac{1}{3}\right) - \gamma_1\left(\frac{2}{3}\right) = \frac{\pi}{\sqrt{3}}\left[6\log\Gamma\left(\frac{1}{3}\right) + \frac{1}{2}\log 3 - (\gamma + 4\log(2\pi))\right]$$

and we then see that this is equivalent to (2.30).

## 5. HIGHER ORDER STIELTJES CONSTANTS

We refer back to (2.12)

$$\varsigma'\left(s,\frac{p}{q}\right) - \varsigma'\left(s,1-\frac{p}{q}\right) = -4\Gamma(1-s)(2\pi q)^{s-1}\cos\left(\frac{\pi s}{2}\right)\sum_{j=1}^{q}\sin\left(\frac{2\pi jp}{q}\right)\varsigma'\left(1-s,\frac{j}{q}\right)$$

$$-4\Gamma(1-s)(2\pi q)^{s-1}\cos\left(\frac{\pi s}{2}\right)\left[\psi(1-s)-\log(2\pi q)+\frac{\pi}{2}\tan\left(\frac{\pi s}{2}\right)\right]\sum_{j=1}^{q}\sin\left(\frac{2\pi jp}{q}\right)\varsigma\left(1-s,\frac{j}{q}\right)$$

where differentiation results in

$$\varsigma''\left(s,\frac{p}{q}\right) - \varsigma''\left(s,1-\frac{p}{q}\right) = 4\Gamma(1-s)(2\pi q)^{s-1}\cos\left(\frac{\pi s}{2}\right)\sum_{j=1}^{q}\sin\left(\frac{2\pi jp}{q}\right)\varsigma''\left(1-s,\frac{j}{q}\right)$$

$$+8\Gamma(1-s)(2\pi q)^{s-1}\cos\left(\frac{\pi s}{2}\right)\left[\psi(1-s)-\log(2\pi q)+\frac{\pi}{2}\tan\left(\frac{\pi s}{2}\right)\right]\sum_{j=1}^{q}\sin\left(\frac{2\pi jp}{q}\right)\varsigma'\left(1-s,\frac{j}{q}\right)$$

$$+4\Gamma(1-s)(2\pi q)^{s-1}\cos\left(\frac{\pi s}{2}\right)\left[\psi(1-s)-\log(2\pi q)+\frac{\pi}{2}\tan\left(\frac{\pi s}{2}\right)\right]^2\sum_{j=1}^{q}\sin\left(\frac{2\pi jp}{q}\right)\varsigma\left(1-s,\frac{j}{q}\right)$$

$$-4\Gamma(1-s)(2\pi q)^{s-1}\cos\left(\frac{\pi s}{2}\right)\left[-\psi'(1-s)+\frac{\pi^2}{4}\sec^2\left(\frac{\pi s}{2}\right)\right]\sum_{j=1}^{q}\sin\left(\frac{2\pi jp}{q}\right)\varsigma\left(1-s,\frac{j}{q}\right)$$

Using (2.4), (2.8) and (2.45) we have for $s=1$

$$\varsigma''\left(1,\frac{p}{q}\right) - \varsigma''\left(1,1-\frac{p}{q}\right) = 2\pi\sum_{j=1}^{q}\sin\left(\frac{2\pi jp}{q}\right)\varsigma''\left(0,\frac{j}{q}\right) - 4\pi[\gamma+\log(2\pi q)]\sum_{j=1}^{q}\sin\left(\frac{2\pi jp}{q}\right)\varsigma'\left(0,\frac{j}{q}\right)$$

$$+2\pi[\gamma+\log(2\pi q)]^2\sum_{j=1}^{q}\sin\left(\frac{2\pi jp}{q}\right)\varsigma\left(0,\frac{j}{q}\right)+\frac{\pi^3}{6}\sum_{j=1}^{q}\sin\left(\frac{2\pi jp}{q}\right)\varsigma\left(0,\frac{j}{q}\right)$$

Using (2.6) this may be simplified to

$$\varsigma''\left(1,\frac{p}{q}\right)-\varsigma''\left(1,1-\frac{p}{q}\right)=2\pi\sum_{j=1}^{q}\sin\left(\frac{2\pi jp}{q}\right)\varsigma''\left(0,\frac{j}{q}\right)-4\pi[\gamma+\log(2\pi q)]\sum_{j=1}^{q}\sin\left(\frac{2\pi jp}{q}\right)\log\Gamma\left(\frac{j}{q}\right)$$

$$-\frac{\pi}{q}\left(2[\gamma+\log(2\pi q)]^2+\frac{\pi^2}{6}\right)\sum_{j=1}^{q}j\sin\left(\frac{2\pi jp}{q}\right)$$

and using (2.11) this becomes

(5.1)

$$\varsigma''\left(1,\frac{p}{q}\right)-\varsigma''\left(1,1-\frac{p}{q}\right)=2\pi\sum_{j=1}^{q}\sin\left(\frac{2\pi jp}{q}\right)\varsigma''\left(0,\frac{j}{q}\right)-4\pi[\gamma+\log(2\pi q)]\sum_{j=1}^{q}\sin\left(\frac{2\pi jp}{q}\right)\log\Gamma\left(\frac{j}{q}\right)$$

$$+\frac{\pi}{2}\left(2[\gamma+\log(2\pi q)]^2+\frac{\pi^2}{6}\right)\cot\left(\frac{\pi p}{q}\right)$$

For example we have

(5.2) $\varsigma''\left(1,\frac{1}{4}\right)-\varsigma''\left(1,\frac{3}{4}\right)$

$$=2\pi\left[\varsigma''\left(0,\frac{1}{4}\right)-\varsigma''\left(0,\frac{3}{4}\right)\right]-4\pi[\gamma+\log(8\pi)]\left[\log\Gamma\left(\frac{1}{4}\right)-\log\Gamma\left(\frac{3}{4}\right)\right]+\frac{\pi}{2}\left(2[\gamma+\log(8\pi)]^2+\frac{\pi^2}{6}\right)$$

Using Euler's reflection formula we see that

$$\Gamma\left(\frac{3}{4}\right)=\pi\sqrt{2}/\Gamma\left(\frac{1}{4}\right)$$

and we therefore get

$$\varsigma''\left(1,\frac{1}{4}\right)-\varsigma''\left(1,\frac{3}{4}\right)=$$

$$2\pi\left[\varsigma''\left(0,\frac{1}{4}\right)-\varsigma''\left(0,\frac{3}{4}\right)\right]-4\pi[\gamma+\log(8\pi)]\left[2\log\Gamma\left(\frac{1}{4}\right)-\log\pi-\frac{1}{2}\log 2\right]+\frac{\pi}{2}\left(2[\gamma+\log(8\pi)]^2+\frac{\pi^2}{6}\right)$$

Referring back to (2.17) we note that

$$\gamma_2(x)-\gamma_2(y)=\lim_{s\to 1}\frac{\partial^2}{\partial s^2}[\varsigma(s,x)-\varsigma(s,y)]$$

and hence we obtain

(5.3) $\quad \gamma_2\left(\dfrac{1}{4}\right)-\gamma_2\left(\dfrac{3}{4}\right)=$

$$2\pi\left[\varsigma''\left(0,\frac{1}{4}\right)-\varsigma''\left(0,\frac{3}{4}\right)\right]-4\pi[\gamma+\log(8\pi)]\left[2\log\Gamma\left(\frac{1}{4}\right)-\log\pi-\frac{1}{2}\log 2\right]+\frac{\pi}{2}\left(2[\gamma+\log(8\pi)]^2+\frac{\pi^2}{6}\right)$$

We recall (2.20)

$$\sum_{r=1}^{q-1}\varsigma\left(s,\frac{r}{q}\right)=(q^s-1)\varsigma(s)$$

and differentiation results in

$$\frac{d}{ds}\sum_{r=1}^{q-1}\varsigma\left(s,\frac{r}{q}\right)=(q^s-1)\varsigma'(s)+\varsigma(s)q^s\log q$$

$$\frac{d^2}{ds^2}\sum_{r=1}^{q-1}\varsigma\left(s,\frac{r}{q}\right)=(q^s-1)\varsigma''(s)+2\varsigma'(s)q^s\log q+\varsigma(s)q^s\log^2 q$$

$$\left.\frac{d^2}{ds^2}\sum_{r=1}^{q-1}\varsigma\left(s,\frac{r}{q}\right)\right|_{s=0}=2\varsigma'(0)\log q+\varsigma(0)\log^2 q$$

It is well known that [42, p.92]

$$\varsigma'(0)=-\frac{1}{2}\log(2\pi) \text{ and } \varsigma(0)=-\frac{1}{2}$$

and we then obtain

$$\left.\frac{d^2}{ds^2}\sum_{r=1}^{q-1}\varsigma\left(s,\frac{r}{q}\right)\right|_{s=0}=-\log(2\pi)\log q-\frac{1}{2}\log^2 q$$

In particular we have for $q = 4$

$$\frac{d^2}{ds^2} \sum_{r=1}^{3} \varsigma\left(s, \frac{r}{4}\right)\bigg|_{s=0} = -\log(2\pi)\log 4 - \frac{1}{2}\log^2 4 = -2\log(2\pi)\log 2 - 2\log^2 2$$

and the left-hand side may be written out explicitly as

$$\varsigma''\left(0, \frac{1}{4}\right) + \varsigma''\left(0, \frac{2}{4}\right) + \varsigma''\left(0, \frac{3}{4}\right) = \varsigma''\left(0, \frac{1}{4}\right) + \varsigma''\left(0, \frac{1}{2}\right) + \varsigma''\left(0, \frac{3}{4}\right)$$

We recall the identity [30] (which may be obtained from (2.20) with $q = 2$)

$$\varsigma\left(s, \frac{1}{2}\right) = [2^s - 1]\varsigma(s)$$

and differentiation results in

$$\varsigma'\left(s, \frac{1}{2}\right) = [2^s - 1]\varsigma'(s) + \varsigma(s)2^s \log 2$$

$$\varsigma''\left(s, \frac{1}{2}\right) = [2^s - 1]\varsigma''(s) + 2\varsigma'(s)2^s \log 2 + \varsigma(s)2^s \log^2 2$$

Specifically we get

(5.4)  $\varsigma''\left(0, \frac{1}{2}\right) = 2\varsigma'(0)\log 2 + \varsigma(0)\log^2 2 = -\log(2\pi)\log 2 - \frac{1}{2}\log^2 2$

We therefore have

(5.5)  $\varsigma''\left(0, \frac{1}{4}\right) + \varsigma''\left(0, \frac{3}{4}\right) = -\log(2\pi)\log 2 - \frac{3}{2}\log^2 2$

and hence we have obtained two simultaneous equations for $\varsigma''\left(0, \frac{1}{4}\right) - \varsigma''\left(0, \frac{3}{4}\right)$ and $\varsigma''\left(0, \frac{1}{4}\right) + \varsigma''\left(0, \frac{3}{4}\right)$ respectively.

□

The Dirichlet beta function is defined by

$$\beta(s) = \sum_{k=0}^{\infty} \frac{(-1)^k}{(2k+1)^s} = \frac{1}{4^s}\left[\varsigma\left(s,\frac{1}{4}\right) - \varsigma\left(s,\frac{3}{4}\right)\right]$$

and differentiation results in

$$\beta'(s) = \sum_{k=0}^{\infty} \frac{(-1)^{k+1}\log(2k+1)}{(2k+1)^s} = \frac{1}{4^s}\left[\varsigma'\left(s,\frac{1}{4}\right) - \varsigma'\left(s,\frac{3}{4}\right)\right] - \frac{\log 4}{4^s}\left[\varsigma\left(s,\frac{1}{4}\right) - \varsigma\left(s,\frac{3}{4}\right)\right]$$

In particular we have

$$\beta'(1) = \sum_{k=0}^{\infty} \frac{(-1)^{k+1}\log(2k+1)}{2k+1} = \frac{1}{4}\left[\varsigma'\left(1,\frac{1}{4}\right) - \varsigma'\left(1,\frac{3}{4}\right)\right] - \frac{1}{2}\log 2\left[\varsigma\left(1,\frac{1}{4}\right) - \varsigma\left(1,\frac{3}{4}\right)\right]$$

From (2.17) and [42, p.20] we have

$$\varsigma\left(1,\frac{1}{4}\right) - \varsigma\left(1,\frac{3}{4}\right) = \gamma_0\left(\frac{1}{4}\right) - \gamma_0\left(\frac{3}{4}\right) = \psi\left(\frac{3}{4}\right) - \psi\left(\frac{1}{4}\right) = \pi$$

Using (2.13) this becomes

(5.6) $$\beta'(1) = \sum_{k=0}^{\infty} \frac{(-1)^{k+1}\log(2k+1)}{2k+1} = \frac{\pi}{4}\left[\gamma + 2\log 2 + 3\log \pi - 4\log \Gamma\left(\frac{1}{4}\right)\right]$$

which was previously obtained in a different manner by Shail [41]. This was also obtained in Berndt's book [7a, Part I, p.201] by letting $x = 3/4$ in Kummer's formula (5.13). A different derivation was given by Malmstén [36a, p.20] in 1849.

Another differentiation gives us

$$\beta''(s) = \sum_{k=0}^{\infty} \frac{(-1)^k \log^2(2k+1)}{(2k+1)^s} = \frac{1}{4^s}\left[\varsigma''\left(s,\frac{1}{4}\right) - \varsigma''\left(s,\frac{3}{4}\right)\right] - 2\frac{\log 4}{4^s}\left[\varsigma'\left(s,\frac{1}{4}\right) - \varsigma'\left(s,\frac{3}{4}\right)\right]$$

$$+ \frac{\log^2 4}{4^s}\left[\varsigma\left(s,\frac{1}{4}\right) - \varsigma\left(s,\frac{3}{4}\right)\right]$$

In particular we have

$$\beta''(1) = \sum_{k=0}^{\infty} \frac{(-1)^k \log^2(2k+1)}{2k+1} = \frac{1}{4}\left[\varsigma''\left(1,\frac{1}{4}\right) - \varsigma''\left(1,\frac{3}{4}\right)\right] - 2\log 2\left[\varsigma'\left(1,\frac{1}{4}\right) - \varsigma'\left(1,\frac{3}{4}\right)\right]$$

$$+ \pi \log^2 2$$

and using (5.2) we may write this as

$$(5.7) \quad \sum_{k=0}^{\infty} \frac{(-1)^k \log^2(2k+1)}{2k+1}$$

$$= \frac{\pi}{4}\left\{2\varsigma''\left(0,\frac{1}{4}\right) - 2\varsigma''\left(0,\frac{3}{4}\right) - 4[\gamma + \log(2\pi)]\log\frac{\Gamma(1/4)}{\Gamma(3/4)} + [\gamma + \log(2\pi)]^2 + \frac{\pi^2}{12}\right\}$$

which was previously obtained in a different manner by Shail [41].

□

We have the well-known Hurwitz's formula for the Fourier expansion of the Riemann zeta function $\varsigma(s,t)$ as reported in Titchmarsh's treatise [44, p.37]

$$(5.8) \quad \varsigma(s,t) = 2\Gamma(1-s)\left[\sin\left(\frac{\pi s}{2}\right)\sum_{n=1}^{\infty}\frac{\cos 2n\pi t}{(2\pi n)^{1-s}} + \cos\left(\frac{\pi s}{2}\right)\sum_{n=1}^{\infty}\frac{\sin 2n\pi t}{(2\pi n)^{1-s}}\right]$$

where $\text{Re}(s) < 0$ and $0 < t \leq 1$. In 2000, Boudjelkha [8] showed that this formula also applies in the region $\text{Re}(s) < 1$. It may be noted that when $t = 1$ this reduces to Riemann's functional equation for $\varsigma(s)$. Letting $s \to 1-s$ we may write this as

$$(5.9) \quad \varsigma(1-s,t) = 2\Gamma(s)\left[\cos\left(\frac{\pi s}{2}\right)\sum_{n=1}^{\infty}\frac{\cos 2n\pi t}{(2\pi n)^s} + \sin\left(\frac{\pi s}{2}\right)\sum_{n=1}^{\infty}\frac{\sin 2n\pi t}{(2\pi n)^s}\right]$$

Letting $s = 0$ in (5.8) gives us

$$(5.10) \quad \varsigma(0,t) = \sum_{n=1}^{\infty}\frac{\sin 2n\pi t}{\pi n}$$

We have the well-known Fourier series [10] valid for $0 < t < 1$

$$\frac{1}{2} - t = \sum_{n=1}^{\infty}\frac{\sin 2n\pi t}{\pi n}$$

and validity is confirmed by noting that $\varsigma(0,t) = \frac{1}{2} - t$

Differentiating (5.8) gives us

$$\varsigma'(s,t) = 2\Gamma(1-s)\sin\left(\frac{\pi s}{2}\right)\sum_{n=1}^{\infty}\frac{\cos 2n\pi t \log(2\pi n)}{(2\pi n)^{1-s}} + \pi\Gamma(1-s)\cos\left(\frac{\pi s}{2}\right)\sum_{n=1}^{\infty}\frac{\cos 2n\pi t}{(2\pi n)^{1-s}}$$

$$-2\Gamma'(1-s)\sin\left(\frac{\pi s}{2}\right)\sum_{n=1}^{\infty}\frac{\cos 2n\pi t}{(2\pi n)^{1-s}}$$

$$+2\Gamma(1-s)\cos\left(\frac{\pi s}{2}\right)\sum_{n=1}^{\infty}\frac{\sin 2n\pi t \log(2\pi n)}{(2\pi n)^{1-s}} - \pi\Gamma(1-s)\sin\left(\frac{\pi s}{2}\right)\sum_{n=1}^{\infty}\frac{\sin 2n\pi t}{(2\pi n)^{1-s}}$$

$$-2\Gamma'(1-s)\cos\left(\frac{\pi s}{2}\right)\sum_{n=1}^{\infty}\frac{\sin 2n\pi t}{(2\pi n)^{1-s}}$$

which may be written as

$$\varsigma'(s,t) = 2\sum_{n=1}^{\infty}\left[\log(2\pi n) - \psi(1-s) + \frac{\pi}{2}\cot\left(\frac{\pi s}{2}\right)\right]\frac{\Gamma(1-s)}{(2\pi n)^{1-s}}\sin\left(\frac{\pi s}{2}\right)\cos 2n\pi t$$

$$+2\sum_{n=1}^{\infty}\left[\log(2\pi n) - \psi(1-s) - \frac{\pi}{2}\tan\left(\frac{\pi s}{2}\right)\right]\frac{\Gamma(1-s)}{(2\pi n)^{1-s}}\cos\left(\frac{\pi s}{2}\right)\sin 2n\pi t$$

With $s = 0$ we have

$$\varsigma'(0,t) = \pi\sum_{n=1}^{\infty}\frac{\cos 2n\pi t}{2\pi n} + 2\sum_{n=1}^{\infty}\frac{\sin 2n\pi t \log(2\pi n)}{2\pi n} - 2\Gamma'(1)\sum_{n=1}^{\infty}\frac{\sin 2n\pi t}{2\pi n}$$

$$= \frac{1}{2}\sum_{n=1}^{\infty}\frac{\cos 2n\pi t}{n} + \log(2\pi)\sum_{n=1}^{\infty}\frac{\sin 2n\pi t}{\pi n} + \frac{1}{\pi}\sum_{n=1}^{\infty}\frac{\log n}{n}\sin 2\pi nt + \gamma\sum_{n=1}^{\infty}\frac{\sin 2n\pi t}{\pi n}$$

Substituting the Fourier series shown in Carslaw's book [41, p.241]

(5.11) $$-\log(2\sin\pi t) = \sum_{n=1}^{\infty}\frac{\cos 2n\pi t}{n}$$

(5.12) $$\frac{1}{2}\pi(1-2t) = \sum_{n=1}^{\infty}\frac{\sin 2n\pi t}{n}$$

we obtain

$$\varsigma'(0,t) = -\frac{1}{2}\log(2\sin \pi t) + \frac{1}{2}(1-2t)[\gamma + \log(2\pi)] + \frac{1}{\pi}\sum_{n=1}^{\infty}\frac{\log n}{n}\sin 2\pi nt$$

and using Lerch's identity (2.5) we deduce Kummer's Fourier series for $\log \Gamma(t)$

(5.13) $$\log \Gamma(t) = \frac{1}{2}\log\frac{\pi}{\sin \pi t} + \frac{1}{2}(1-2t)[\gamma + \log(2\pi)] + \frac{1}{\pi}\sum_{n=1}^{\infty}\frac{\log n}{n}\sin 2\pi nt$$

□

Letting $t \to 1-t$ in (5.8) we get

$$\varsigma(s,1-t) = 2\Gamma(1-s)\left[\sin\left(\frac{\pi s}{2}\right)\sum_{n=1}^{\infty}\frac{\cos 2n\pi t}{(2\pi n)^{1-s}} - \cos\left(\frac{\pi s}{2}\right)\sum_{n=1}^{\infty}\frac{\sin 2n\pi t}{(2\pi n)^{1-s}}\right]$$

and we therefore see that

(5.14) $$\varsigma(s,t) + \varsigma(s,1-t) = 4\Gamma(1-s)\sin\left(\frac{\pi s}{2}\right)\sum_{n=1}^{\infty}\frac{\cos 2n\pi t}{(2\pi n)^{1-s}}$$

(5.15) $$\varsigma(s,t) - \varsigma(s,1-t) = 4\Gamma(1-s)\cos\left(\frac{\pi s}{2}\right)\sum_{n=1}^{\infty}\frac{\sin 2n\pi t}{(2\pi n)^{1-s}}$$

Differentiation gives us

$$\varsigma'(s,t) + \varsigma'(s,1-t)$$

$$= 4\Gamma(1-s)\sin\left(\frac{\pi s}{2}\right)\sum_{n=1}^{\infty}\frac{\cos 2n\pi t \log(2\pi n)}{(2\pi n)^{1-s}} + 2\pi\Gamma(1-s)\cos\left(\frac{\pi s}{2}\right)\sum_{n=1}^{\infty}\frac{\cos 2n\pi t}{(2\pi n)^{1-s}}$$

$$-4\Gamma'(1-s)\sin\left(\frac{\pi s}{2}\right)\sum_{n=1}^{\infty}\frac{\cos 2n\pi t}{(2\pi n)^{1-s}}$$

With $s = 0$ we have

$$\varsigma'(0,t) + \varsigma'(0,1-t) = \sum_{n=1}^{\infty}\frac{\cos 2n\pi t}{n}$$

and using (5.11) we obtain

$$\varsigma'(0,t) + \varsigma'(0,1-t) = -\log[2\sin \pi t]$$

This formula may also be obtained directly from Lerch's identity (2.5).

A further differentiation results in

$$\varsigma''(s,t)+\varsigma''(s,1-t)$$

$$=4\Gamma(1-s)\sin\left(\frac{\pi s}{2}\right)\sum_{n=1}^{\infty}\frac{\cos 2n\pi t \log^2(2\pi n)}{(2\pi n)^{1-s}}$$

$$+\left[2\pi\,\Gamma(1-s)\cos\left(\frac{\pi s}{2}\right)-4\Gamma'(1-s)\sin\left(\frac{\pi s}{2}\right)\right]\sum_{n=1}^{\infty}\frac{\cos 2n\pi t \log(2\pi n)}{(2\pi n)^{1-s}}$$

$$+2\pi\,\Gamma(1-s)\cos\left(\frac{\pi s}{2}\right)\sum_{n=1}^{\infty}\frac{\cos 2n\pi t \log(2\pi n)}{(2\pi n)^{1-s}}$$

$$-2\pi\left[\Gamma(1-s)\frac{\pi}{2}\sin\left(\frac{\pi s}{2}\right)+\Gamma'(1-s)\cos\left(\frac{\pi s}{2}\right)\right]\sum_{n=1}^{\infty}\frac{\cos 2n\pi t}{(2\pi n)^{1-s}}$$

$$-4\Gamma'(1-s)\sin\left(\frac{\pi s}{2}\right)\sum_{n=1}^{\infty}\frac{\cos 2n\pi t \log(2\pi n)}{(2\pi n)^{1-s}}$$

$$-4\left[\Gamma'(1-s)\frac{\pi}{2}\cos\left(\frac{\pi s}{2}\right)-\Gamma''(1-s)\sin\left(\frac{\pi s}{2}\right)\right]\sum_{n=1}^{\infty}\frac{\cos 2n\pi t}{(2\pi n)^{1-s}}$$

With $s=0$ we have

$$\varsigma''(0,t)+\varsigma''(0,1-t)=2\sum_{n=1}^{\infty}\frac{\cos 2n\pi t \log(2\pi n)}{n}+2\gamma\sum_{n=1}^{\infty}\frac{\cos 2n\pi t}{n}$$

With $t=1/2$ we have

$$\varsigma''\left(0,\frac{1}{2}\right)=\sum_{n=1}^{\infty}\frac{(-1)^n \log(2\pi n)}{n}+\gamma\sum_{n=1}^{\infty}\frac{(-1)^n}{n}$$

$$=-\log(2\pi)\log 2+\sum_{n=1}^{\infty}\frac{(-1)^n \log n}{n}-\gamma\log 2$$

$$=-\log(2\pi)\log 2+\varsigma'_a(1)-\gamma\log 2$$

and using (7.13) we rediscover (5.4)

$$\varsigma''\left(0, \frac{1}{2}\right) = -\log(2\pi)\log 2 - \frac{1}{2}\log^2 2$$

Letting $t = 1/4$ we have

$$\varsigma''\left(0, \frac{1}{4}\right) + \varsigma''\left(0, \frac{3}{4}\right) = 2\sum_{n=1}^{\infty}\frac{\cos(n\pi/2)\log n}{n} + 2[\gamma + \log(2\pi)]\sum_{n=1}^{\infty}\frac{\cos(n\pi/2)}{n}$$

We note that

$$\sum_{n=1}^{\infty}\frac{\cos(n\pi/2)\log n}{n} = -\frac{\log 2}{2} + \frac{\log 4}{4} - \frac{\log 6}{6} \cdots = \sum_{n=1}^{\infty}\frac{(-1)^n \log(2n)}{2n}$$

$$= \frac{1}{2}\log 2 \sum_{n=1}^{\infty}\frac{(-1)^n}{n} + \frac{1}{2}\sum_{n=1}^{\infty}\frac{(-1)^n \log n}{n}$$

$$= -\frac{1}{2}\log^2 2 + \frac{1}{2}\varsigma_a^{(1)}(1)$$

and we recall from (2.34) that

$$\varsigma_a^{(1)}(1) = \gamma \log 2 - \frac{1}{2}\log^2 2$$

From (5.11) we see that

$$\sum_{n=1}^{\infty}\frac{\cos(n\pi/2)}{n} = -\log[2\sin(\pi/4)] = -\frac{1}{2}\log 2$$

and hence we have

$$\varsigma''\left(0, \frac{1}{4}\right) + \varsigma''\left(0, \frac{3}{4}\right) = -\log 2 \log(2\pi) - \frac{3}{2}\log 2$$

in agreement with (5.5).

We also see from (5.15) that

$$\varsigma'(s,t) - \varsigma'(s, 1-t)$$

$$= -4\Gamma(1-s)\cos\left(\frac{\pi s}{2}\right)\sum_{n=1}^{\infty}\frac{\sin 2n\pi t \log(2\pi n)}{(2\pi n)^{1-s}} - 2\pi\Gamma(1-s)\sin\left(\frac{\pi s}{2}\right)\sum_{n=1}^{\infty}\frac{\sin 2n\pi t}{(2\pi n)^{1-s}}$$

$$-4\Gamma'(1-s)\cos\left(\frac{\pi s}{2}\right)\sum_{n=1}^{\infty}\frac{\sin 2n\pi t}{(2\pi n)^{1-s}}$$

so that

$$\varsigma'(0,t) - \varsigma'(0,1-t) = -2\sum_{n=1}^{\infty}\frac{\sin 2n\pi t \log(2\pi n)}{\pi n} + 2\gamma\sum_{n=1}^{\infty}\frac{\sin 2n\pi t}{\pi n}$$

$$= -2\sum_{n=1}^{\infty}\frac{\sin 2n\pi t \log n}{\pi n} + 2[\gamma - \log(2\pi)]\sum_{n=1}^{\infty}\frac{\sin 2n\pi t}{\pi n}$$

$$= -2\sum_{n=1}^{\infty}\frac{\sin 2n\pi t \log n}{\pi n} + 2[\gamma - \log(2\pi)]\left(\frac{1}{2} - t\right)$$

## 6. INTEGRAL REPRESENTATIONS OF THE CLAUSEN FUNCTIONS

The Clausen functions $\text{Cl}_N(\theta)$ are defined by [42, p.115]

$$\text{Cl}_{2N}(\theta) = \sum_{n=1}^{\infty}\frac{\sin n\theta}{n^{2N}}$$

$$\text{Cl}_{2N+1}(\theta) = \sum_{n=1}^{\infty}\frac{\cos n\theta}{n^{2N+1}}$$

and may be represented by the following integrals [24]

(6.1a) $\quad \text{Cl}_{2n}(\theta) = -\frac{\sin\theta}{(2n-1)!}\int_0^1 \frac{\log^{2n-1} x}{1 - 2x\cos\theta + x^2}dx$

(6.1b) $\quad \text{Cl}_{2n+1}(\theta) = -\frac{1}{(2n)!}\int_0^1 \frac{(x - \cos\theta)\log^{2n} x}{1 - 2x\cos\theta + x^2}dx$

A novel way of deriving the above two identities has been recently given by Efthimiou [26]. His method is shown below.

Since, for $\text{Re}(s) > 0$, $\frac{1}{n^s} = \frac{1}{\Gamma(s)}\int_0^{\infty} e^{-nt} t^{s-1} dt$ we may write

$$\sum_{n=1}^{\infty} \frac{\cos nx}{n^s} = \frac{1}{\Gamma(s)} \sum_{n=1}^{\infty} \cos nx \int_0^{\infty} e^{-nt} t^{s-1} dt$$

$$= \frac{1}{\Gamma(s)} \int_0^{\infty} \left( \sum_{n=1}^{\infty} \cos nx \, e^{-nt} \right) t^{s-1} dt$$

$$= \frac{1}{\Gamma(s)} \int_0^{\infty} \frac{e^{-t}(\cos x - e^{-t})}{1 - 2\cos x \, e^{-t} + e^{-2t}} t^{s-1} dt$$

Letting $u = e^{-t}$ we obtain

(6.2a) $$\sum_{n=1}^{\infty} \frac{\cos nx}{n^s} = \frac{(-1)^s}{\Gamma(s)} \int_0^1 \frac{(u - \cos x) \log^{s-1} u}{1 - 2u \cos x + u^2} du$$

Similarly we obtain

(6.2b) $$\sum_{n=1}^{\infty} \frac{\sin nx}{n^s} = \frac{(-1)^{s-1} \sin x}{\Gamma(s)} \int_0^1 \frac{\log^{s-1} u}{1 - 2u \cos x + u^2} du$$

We may also write (6.2b) as

(6.3) $$\sum_{n=1}^{\infty} \frac{\sin nx}{n^s} = \frac{\sin x}{\Gamma(s)} \int_0^1 \frac{\log^{s-1}(1/u)}{1 - 2u \cos x + u^2} du$$

and then differentiate with respect to $s$ to obtain

(6.4) $$-\sum_{n=1}^{\infty} \frac{\log n}{n^s} \sin nx = \frac{\sin x}{\Gamma(s)} \int_0^1 \frac{\log^{s-1}(1/u) \log \log(1/u)}{1 - 2u \cos x + u^2} du$$

$$- \frac{\psi(s) \sin x}{\Gamma(s)} \int_0^1 \frac{\log^{s-1}(1/u)}{1 - 2u \cos x + u^2} du$$

or alternatively using (6.3)

(6.5) $$-\sum_{n=1}^{\infty} \frac{\log n}{n^s} \sin nx = \frac{\sin x}{\Gamma(s)} \int_0^1 \frac{\log^{s-1}(1/u) \log \log(1/u)}{1 - 2u \cos x + u^2} du - \psi(s) \sum_{n=1}^{\infty} \frac{\sin nx}{n^s}$$

With $s = 1$ we have

$$-\sum_{n=1}^{\infty}\frac{\log n}{n}\sin nx=\sin x\int_{0}^{1}\frac{\log\log(1/u)}{1-2u\cos x+u^{2}}du+\gamma\sum_{n=1}^{\infty}\frac{\sin nx}{n}$$

and with $x\to 2\pi x$ this becomes

$$-\sum_{n=1}^{\infty}\frac{\log n}{n}\sin 2\pi nx=\sin 2\pi x\int_{0}^{1}\frac{\log\log(1/u)}{1-2u\cos 2\pi x+u^{2}}du+\gamma\sum_{n=1}^{\infty}\frac{\sin 2\pi nx}{n}$$

Then referring to Kummer's Fourier series (5.13)

$$\log\Gamma(x)=\frac{1}{2}\log\frac{\pi}{\sin\pi x}+\frac{1}{2}(1-2x)[\gamma+\log(2\pi)]+\frac{1}{\pi}\sum_{n=1}^{\infty}\frac{\log n}{n}\sin 2\pi nx$$

we see that

$$\sin 2\pi x\int_{0}^{1}\frac{\log\log(1/u)}{1-2u\cos 2\pi x+u^{2}}du=$$

$$-\frac{1}{2}\gamma\pi(1-2x)+\frac{1}{2}\pi\log\frac{\pi}{\sin\pi x}+\frac{1}{2}\pi(1-2x)[\gamma+\log(2\pi)]-\pi\log\Gamma(x)$$

We therefore obtain for $0<x<1$

(6.6)

$$\sin 2\pi x\int_{0}^{1}\frac{\log\log(1/u)}{1-2u\cos 2\pi x+u^{2}}du=\frac{1}{2}\pi\log\frac{\pi}{\sin\pi x}+\frac{1}{2}\pi(1-2x)\log(2\pi)-\pi\log\Gamma(x)$$

It is easily seen that both sides of the above equation vanish at $x=1/2$. With $x=1/4$ we have

(6.7) $$\int_{0}^{1}\frac{\log\log(1/u)}{1+u^{2}}du=\frac{3}{4}\pi\log\pi+\frac{1}{2}\pi\log 2-\pi\log\Gamma\left(\frac{1}{4}\right)$$

One is then immediately reminded of Adamchik's result (3.9)

$$\int_{0}^{1}\frac{u^{p-1}}{1+u^{n}}\log\log\left(\frac{1}{u}\right)du=\frac{\gamma+\log(2n)}{2n}\left[\psi\left(\frac{p}{2n}\right)-\psi\left(\frac{n+p}{2n}\right)\right]+\frac{1}{2n}\left[\varsigma'\left(1,\frac{p}{2n}\right)-\varsigma'\left(1,\frac{n+p}{2n}\right)\right]$$

where with $p=1$ and $n=2$ we get

$$\int_0^1 \frac{\log\log(1/u)}{1+u^2} du = \frac{1}{4}[\gamma + 2\log 2]\left[\psi\left(\frac{1}{4}\right) - \psi\left(\frac{3}{4}\right)\right] + \frac{1}{4}\left[\varsigma'\left(1,\frac{1}{4}\right) - \varsigma'\left(1,\frac{3}{4}\right)\right]$$

We have previously seen in (2.13) that

$$\varsigma'\left(1,\frac{1}{4}\right) - \varsigma'\left(1,\frac{3}{4}\right) = \pi\left[\gamma + 4\log 2 + 3\log \pi - 4\log\Gamma\left(\frac{1}{4}\right)\right]$$

and it is well known that [42, p.20]

$$\psi\left(\frac{1}{4}\right) - \psi\left(\frac{3}{4}\right) = -\pi$$

and hence we obtain

$$\int_0^1 \frac{\log\log(1/u)}{1+u^2} du = -\frac{\pi}{4}[\gamma + 2\log 2] + \frac{1}{4}\pi\left[\gamma + 4\log 2 + 3\log\pi - 4\log\Gamma\left(\frac{1}{4}\right)\right]$$

which is the same as (6.7).

With $x = 1/6$ in (6.6) we get

$$\frac{\sqrt{3}}{2}\int_0^1 \frac{\log\log(1/u)}{1-u+u^2} du = \frac{1}{2}\pi\log(2\pi) + \frac{\pi}{3}\log(2\pi) - \pi\log\Gamma\left(\frac{1}{6}\right)$$

$$= \frac{5}{6}\pi\log(2\pi) - \pi\log\Gamma\left(\frac{1}{6}\right)$$

which is in agreement with Adamchik's paper [2].

We may also write (6.2b) as

$$\sum_{n=1}^{\infty} \frac{\cos nx}{n^s} = -\frac{1}{\Gamma(s)} \int_0^1 \frac{(u - \cos x)\log^{s-1}(1/u)}{1 - 2u\cos x + u^2} du$$

and then differentiate with respect to $s$ to obtain

$$-\sum_{n=1}^{\infty} \frac{\log n}{n^s} \cos nx = -\frac{1}{\Gamma(s)} \int_0^1 \frac{(u - \cos x)\log^{s-1}(1/u)\log\log(1/u)}{1 - 2u\cos x + u^2} du$$

$$+\frac{\psi(s)}{\Gamma(s)}\int_0^1 \frac{(u-\cos x)\log^{s-1}(1/u)}{1-2u\cos x+u^2}du$$

and this may be written as

(6.8)

$$\frac{1}{\Gamma(s)}\int_0^1 \frac{(u-\cos x)\log^{s-1}(1/u)\log\log(1/u)}{1-2u\cos x+u^2}du = \sum_{n=1}^{\infty}\frac{\log n}{n^s}\cos nx - \psi(s)\sum_{n=1}^{\infty}\frac{\cos nx}{n^s}$$

For example, with $s=2$ and $x=\pi/2$ we may note that

$$\sum_{n=1}^{\infty}\frac{\cos(n\pi/2)}{n^2} = -\frac{1}{2^2}+\frac{1}{4^2}-\frac{1}{6^2}\cdots$$

$$= -\frac{1}{2^2}\left[\frac{1}{1^2}-\frac{1}{2^2}+\frac{1}{3^2}\cdots\right]$$

$$= -\frac{1}{4}\varsigma_a(2) = -\frac{1}{8}\varsigma(2)$$

and it is known that

$$\sum_{n=1}^{\infty}\frac{\log n\cos(n\pi/2)}{n^2} = \frac{\pi^2}{48}[\log(2\pi)+\gamma-1]+\frac{\pi^2}{4}\varsigma'(-1)$$

Therefore we have

(6.9)

$$\int_0^1 \frac{u\log(1/u)\log\log(1/u)}{1+u^2}du = \frac{\pi^2}{48}[\log(2\pi)+\gamma-1]+\frac{\pi^2}{4}\varsigma'(-1)+\frac{1}{8}(1-\gamma)\varsigma(2)$$

Reference should also be made to the 2002 paper by Koyama and Kurokawa, "Kummer's formula for the multiple gamma functions" [34] where they show that

(6.10) $\quad \log\Gamma_2^*(x) = -\frac{1}{2\pi^2}\sum_{n=1}^{\infty}\frac{\log n}{n^2}\cos 2\pi nx - \frac{\log(2\pi)+\gamma-1}{2\pi^2}\sum_{n=1}^{\infty}\frac{\cos 2\pi nx}{n^2}$

$$+\frac{1}{4\pi}\sum_{n=1}^{\infty}\frac{\sin 2\pi nx}{n^2}+(1-x)\log\Gamma_1(x)$$

(6.11) $$\log \Gamma_3^*(x) = -\frac{1}{4\pi^3} \sum_{n=1}^{\infty} \frac{\log n}{n^3} \sin 2\pi nx - \frac{2\log(2\pi) + 2\gamma - 3}{8\pi^3} \sum_{n=1}^{\infty} \frac{\sin 2\pi nx}{n^3}$$

$$+ \frac{1}{8\pi^2} \sum_{n=1}^{\infty} \frac{\cos 2\pi nx}{n^3} + \left(\frac{3}{2} - x\right) \log \Gamma_2^*(x) - \frac{1}{2}(1-x)^2 \log \Gamma_1^*(x)$$

where $\Gamma_1^*(x) = \frac{\Gamma(x)}{\sqrt{2\pi}}$. These equations will then enable us to evaluate the integral in (6.8) when $s = 2$ and the integral in (6.5) when $s = 3$. It should however be noted that the multiple gamma functions $\Gamma_n^*(x)$ considered by Koyama and Kurokawa [35] are not exactly the same as those traditionally employed by Barnes [6], Adamchik [3] and Srivastava and Choi [42] etc.

We could also differentiate (6.6) with respect to $x$ to give us the integral

(6.12)
$$-4\pi \sin^2 2\pi x \int_0^1 \frac{u \log \log(1/u)}{[1 - 2u \cos 2\pi x + u^2]^2} du + 2\pi \cos 2\pi x \int_0^1 \frac{\log \log(1/u)}{1 - 2u \cos 2\pi x + u^2} du =$$

$$\frac{1}{2} \pi^2 \cot \pi x - \pi \log(2\pi) - \pi \psi(x)$$

and we note that Adamchik [2] has also evaluated integrals such as

$$\int_0^1 \frac{u \log \log(1/u)}{[1 + u^2]^2} du$$

which arise from (6.12) when $x = 1/4$.

We may write (6.6) as

$$\sin 2\pi x \int_0^1 \frac{\log \log(1/u)}{1 - 2u \cos 2\pi x + u^2} du = \frac{\pi}{2} \log \left[\frac{\pi}{\sin(\pi x)(2\pi)^{2x-1} \Gamma^2(x)}\right]$$

and with the substitution $u = e^{-t}$ this becomes

(6.13) $$\sin 2\pi x \int_0^1 \frac{e^{-t} \log t}{e^{-2t} - 2e^{-t} \cos 2\pi x + 1} dt = \frac{\pi}{2} \log \left[\frac{\Gamma(1-x)}{(2\pi)^{2x-1} \Gamma(x)}\right]$$

This integral was previously derived by Yue and Williams [48] in 1993. Reference should also be made to the recent paper by Medina and Moll [37] who considered integrals such as (6.6).

It is easily seen that

$$\int_0^1 \frac{e^{-t}\log t}{e^{-2t}-2e^{-t}\cos 2\pi x+1}\,dt = \int_0^1 \frac{\log t}{e^t+e^{-t}-2\cos 2\pi x}\,dt$$

$$= \frac{1}{2}\int_0^1 \frac{\log t}{\cosh t-\cos 2\pi x}\,dt$$

Accordingly (6.13) may be written as

$$\sin 2\pi x \int_0^1 \frac{\log t}{\cosh t-\cos 2\pi x}\,dt = \pi(1-2x)\log(2\pi)+\pi\log\left[\frac{\Gamma(1-x)}{\Gamma(x)}\right]$$

We also have

$$\int_0^1 \frac{\log\log(1/u)}{1-2u\cos 2\pi x+u^2}\,du = \frac{1}{2}\int_0^\infty \frac{\log t}{\cosh t-\cos 2\pi x}\,dt$$

Glasser [26a] reported in 1971 that

(6.14) $\sin 2\pi x \int_0^\infty \dfrac{\log(\gamma t)}{\cosh t-\cos 2\pi x}\,dt = \pi(1-2x)\log(2\pi\gamma)+\pi\log\left[\dfrac{\Gamma(1-x)}{\Gamma(x)}\right]$

and the Wolfram Mathematica Online Integrator gives us

$$\int \frac{1}{\cosh t-a}\,dt = \frac{2}{\sqrt{1-a^2}}\tan^{-1}\left[\frac{(1+a)\tanh(t/2)}{\sqrt{1-a^2}}\right]+c$$

and therefore we have the definite integral

$$\sin 2\pi x \int_0^\infty \frac{1}{\cosh t-\cos 2\pi x}\,dt = 2\tan^{-1}\left[\frac{(1+\cos 2\pi x)}{\sin 2\pi x}\right]$$

$$= 2\tan^{-1}\cot \pi x$$

Since for $z>0$ we have

$$\tan^{-1}\left(\frac{1}{z}\right) = \frac{\pi}{2} - \tan^{-1} z$$

and thus $2\tan^{-1}\cot \pi x = \pi(1-2x)$.

We according see that Glasser's result (6.14) is equivalent to the one subsequently obtained by Yue and Williams (6.13).

In fact, the history goes back even further. With regard to (6.6) we note that the Swedish mathematician Malmstén [36a, p.24] computed the following integral in 1849

$$(6.15) \qquad \sin 2\pi x \int_0^1 \frac{\log\log(1/u)}{1+2u\cos 2\pi x + u^2}\, du = \frac{1}{2}\pi \log \frac{(2\pi)^{2x}\Gamma\left(\frac{1}{2}+x\right)}{\Gamma\left(\frac{1}{2}-x\right)}$$

and this may be easily derived from (6.6) with the substitution $x \to x+1/2$ where we also employ the well-known relation

$$\Gamma\left(x+\frac{1}{2}\right)\Gamma\left(\frac{1}{2}-x\right) = \frac{\pi}{\cos \pi x}$$

which results from Euler's reflection formula with the substitution $x \to x+1/2$.

The integral (6.15) was also recently reported by Medina and Moll [37].

Integrating (6.6) may also give rise to interesting results.

## 7. SOME CONNECTIONS WITH THE ALTERNATING HURWITZ ZETA FUNCTION $\varsigma_a(s,t)$

Upon a separation of terms according to the parity of $n$ we see that for $\sigma > 1$

$$\varsigma_a(s,t) = \sum_{n=0}^{\infty} \frac{(-1)^n}{(n+t)^s} = \sum_{n=0}^{\infty} \frac{1}{(2n+t)^s} - \sum_{n=0}^{\infty} \frac{1}{(2n+1+t)^s}$$

$$= 2^{-s}\left[\sum_{n=0}^{\infty} \frac{1}{(n+t/2)^s} - \sum_{n=0}^{\infty} \frac{1}{(n+(t+1)/2)^s}\right]$$

and we therefore see that $\varsigma_a(s,t)$ is related to the Hurwitz zeta function by the formula

$$(7.1) \quad \varsigma_a(s,t) = 2^{-s}\left[\varsigma\left(s,\frac{t}{2}\right) - \varsigma\left(s,\frac{1+t}{2}\right)\right]$$

Hansen and Patrick [25] showed in 1962 that the Hurwitz zeta function could be written as

$$(7.2) \quad \varsigma(s,x) = 2^s \varsigma(s,2x) - \varsigma\left(s, x+\frac{1}{2}\right)$$

and, by analytic continuation, this holds for all $s$. With $x = t/2$ this becomes

$$(7.3) \quad \varsigma\left(s,\frac{t}{2}\right) = 2^s \varsigma(s,t) - \varsigma\left(s,\frac{1+t}{2}\right)$$

and hence we have for $\sigma > 1$

$$(7.4) \quad \varsigma_a(s,t) = \varsigma(s,t) - 2^{1-s}\varsigma\left(s,\frac{1+t}{2}\right)$$

and

$$(7.5) \quad \varsigma_a(s,t) = 2^{1-s}\varsigma\left(s,\frac{t}{2}\right) - \varsigma(s,t)$$

Since $\varsigma(s,t)$ can be continued analytically to the whole complex plane except for a simple pole at $s=1$, $\varsigma_a(s,t)$ can be continued analytically to become an entire function and (7.1), (7.4) and (7.5) therefore hold in the whole complex plane.

Referring back to (7.5) we see that

$$(7.6) \quad \varsigma_a(s,t) = 2^{1-s}\left[\varsigma\left(s,\frac{t}{2}\right) - \varsigma(s,t)\right] + [2^{1-s} - 1]\varsigma(s,t)$$

and we easily obtain the limit

$$(7.7) \quad \lim_{s \to 1}\varsigma_a(s,t) = \gamma_0\left(\frac{t}{2}\right) - \gamma_0(t) + \lim_{s \to 1}\left[\frac{2^{1-s}-1}{s-1}[(s-1)\varsigma(s,t)]\right]$$

Using L'Hôpital's rule we get

$$(7.8) \quad \lim_{s \to 1}\varsigma_a(s,t) = \sum_{n=0}^{\infty}\frac{(-1)^n}{n+t} = \gamma_0\left(\frac{t}{2}\right) - \gamma_0(t) - \log 2$$

or equivalently

(7.9a) $$\sum_{n=0}^{\infty}\frac{(-1)^n}{n+t}=\psi(t)-\psi\left(\frac{t}{2}\right)-\log 2$$

Similarly, using (7.1) and (7.4) we may obtain

(7.9b) $$\sum_{n=0}^{\infty}\frac{(-1)^n}{n+t}=\frac{1}{2}\left[\psi\left(\frac{1+t}{2}\right)-\psi\left(\frac{t}{2}\right)\right]$$

which appears in [46, p.262]

(7.9c) $$\sum_{n=0}^{\infty}\frac{(-1)^n}{n+t}=\psi\left(\frac{1+t}{2}\right)-\psi(t)+\log 2$$

The above identities (7.9a), (7.9b) and (7.9c) were also derived in [47].

Equation (7.6) may be written as

(7.10) $$\varsigma_a(s,t)=2^{1-s}\left[\varsigma\left(s,\frac{t}{2}\right)-\varsigma(s,t)\right]+\frac{2^{1-s}-1}{s-1}[(s-1)\varsigma(s,t)]$$

and differentiating (7.10) gives us

$$\varsigma_a'(s,t)=2^{1-s}\left[\varsigma'\left(s,\frac{t}{2}\right)-\varsigma'(s,t)\right]-2^{1-s}\log 2\left[\varsigma\left(s,\frac{t}{2}\right)-\varsigma(s,t)\right]$$

$$+f(s)\frac{d}{ds}[(s-1)\varsigma(s,t)]+(s-1)\varsigma(s,t)f'(s)$$

where we define $f(s)$ by

$$f(s)=\frac{2^{1-s}-1}{s-1}$$

We see that

$$2^{1-s}-1=\exp[-(s-1)\log 2]-1=\sum_{k=1}^{\infty}\frac{(-1)^k(s-1)^k\log^k 2}{k!}$$

and differentiating

$$\frac{2^{1-s}-1}{s-1}=\sum_{k=1}^{\infty}\frac{(-1)^k(s-1)^{k-1}\log^k 2}{k!}$$

results in

(7.11) $$f^{(n)}(1) = (-1)^{n+1} \frac{\log^{n+1} 2}{n+1}$$

Using (1.15) and (2.17) this gives us

$$\varsigma_a'(1,t) = -\left[\gamma_1\left(\frac{t}{2}\right) - \gamma_1(t)\right] - \log 2\left[\gamma_0\left(\frac{t}{2}\right) - \gamma_0(t)\right] - \gamma_0(t)\log 2 + \frac{1}{2}\log^2 2$$

$$= -\left[\gamma_1\left(\frac{t}{2}\right) - \gamma_1(t)\right] - \gamma_0\left(\frac{t}{2}\right)\log 2 + \frac{1}{2}\log^2 2$$

and with $t=1$ we have

(7.12) $$\varsigma_a'(1,1) = \varsigma_a'(1) = \gamma_1 - \gamma_1\left(\frac{1}{2}\right) - \gamma_0\left(\frac{1}{2}\right)\log 2 + \frac{1}{2}\log^2 2$$

Since $\varsigma_a(s) = (1-2^{1-s})\varsigma(s)$ we have

$$\varsigma_a(s) = -\frac{2^{1-s}-1}{s-1}[(s-1)\varsigma(s)]$$

and using (7.11) we obtain

(7.13) $$\varsigma_a'(1) = \gamma \log 2 - \frac{1}{2}\log^2 2$$

Using this and (2.22) we see from (7.12) that

(7.14) $$\gamma_1\left(\frac{1}{2}\right) = \gamma_1 - \log^2 2 - 2\gamma \log 2$$

which we also saw in equation (4.3.233) in [21].

We note from (2.27) that

$$\gamma_1\left(\frac{1}{4}\right) = \frac{1}{2}[2\gamma_1 - 15\log^2 2 - 6\gamma \log 2] - \frac{1}{2}\pi\left[\gamma + 4\log 2 + 3\log \pi - 4\log \Gamma\left(\frac{1}{4}\right)\right]$$

and hence we may compute $\varsigma'_a\left(1, \frac{1}{2}\right)$.

Differentiating the Hasse formula for the alternating Hurwitz zeta function gives us

$$\varsigma'_a(s,t) = \frac{\partial}{\partial s} \varsigma_a(s,t) = -\sum_{n=0}^{\infty} \frac{1}{2^{n+1}} \sum_{k=0}^{n} \binom{n}{k} \frac{(-1)^k \log(t+k)}{(t+k)^s}$$

and we therefore see how this connects with equation (4.16) of [23a].

## 8. AN APPLICATION OF THE ABEL-PLANA SUMMATION FORMULA

Adamchik [3a] noted that the Hermite integral for the Hurwitz zeta function may be derived from the Abel-Plana summation formula [46, p.90]

(8.1) $$\sum_{k=0}^{\infty} f(k) = \frac{1}{2} f(0) + \int_0^{\infty} f(x)\,dx + i\int_0^{\infty} \frac{f(ix) - f(-ix)}{e^{2\pi x} - 1}\,dx$$

which applies to functions which are analytic in the right-hand plane and satisfy the convergence condition

$$\lim_{y \to \infty} e^{-2\pi y} |f(x+iy)| = 0$$

uniformly on any finite interval of $x$. Derivations of the Abel-Plana summation formula are given in [49, p.145], [50, p.108] and [31a, p.338].

Letting $f(k) = (k+t)^{-s}$ we obtain

(8.2) $$\varsigma(s,t) = \sum_{k=0}^{\infty} \frac{1}{(k+t)^s} = \frac{t^{-s}}{2} + \frac{t^{1-s}}{s-1} + i\int_0^{\infty} \frac{(t+ix)^{-s} - (t-ix)^{-s}}{e^{2\pi x} - 1}\,dx$$

Then, noting that

$$(t+ix)^{-s} - (t-ix)^{-s} = (re^{i\theta})^{-s} - (re^{-i\theta})^{-s}$$

$$= r^{-s}[e^{-is\theta} - e^{is\theta}]$$

$$= \frac{2}{i(t^2+x^2)^{s/2}} \sin(s\tan^{-1}(x/t))$$

we may write (8.2) as Hermite's integral for $\varsigma(s,u)$

$$(8.3) \qquad \varsigma(s,t) = \frac{t^{-s}}{2} + \frac{t^{1-s}}{s-1} + 2\int_0^\infty \frac{\sin(s\tan^{-1}(x/t))}{(t^2+x^2)^{s/2}(e^{2\pi x}-1)}dx$$

Some applications of (8.2) were considered in [23a].

Here we now apply this methodology to the alternating Hurwitz zeta function. Using (8.2) we have

$$\varsigma\left(s,\frac{t}{2}\right) = 2^{s-1}t^{-s} + \frac{2^{s-1}t^{1-s}}{s-1} + i\int_0^\infty \frac{(t/2+ix)^{-s} - (t/2-ix)^{-s}}{e^{2\pi x}-1}dx$$

and with the substitution $x = y/2$ this becomes

$$= 2^{s-1}t^{-s} + \frac{2^{s-1}t^{1-s}}{s-1} + i2^{-s-1}\int_0^\infty \frac{(t+iy)^{-s} - (t-iy)^{-s}}{e^{\pi y}-1}dy$$

Then using (7.5) we have for all $s$ and $0 < t \leq 1$

$$(8.4) \qquad \varsigma_a(s,t) = \frac{t^{-s}}{2} + i\int_0^\infty \frac{[(t+ix)^{-s} - (t-ix)^{-s}]e^{\pi x}}{e^{2\pi x}-1}dx$$

or equivalently

$$(8.5) \qquad \varsigma_a(s,t) = \frac{t^{-s}}{2} + 2\int_0^\infty \frac{\sin(s\tan^{-1}(x/t))e^{\pi x}}{(t^2+x^2)^{s/2}(e^{2\pi x}-1)}dx$$

Differentiating (8.4) gives us

$$\varsigma_a'(s,t) = -\frac{t^{-s}\log t}{2} + i\int_0^\infty \frac{[(t-ix)^{-s}\log(t-ix) - (t+ix)^{-s}\log(t+ix)]e^{\pi x}}{e^{2\pi x}-1}dx$$

and in particular we have

$$\varsigma_a'(1,t) = -\frac{1}{2}\frac{\log t}{t} + i\int_0^\infty \frac{[(t+ix)\log(t-ix) - (t-ix)\log(t+ix)]e^{\pi x}}{(t^2+x^2)(e^{2\pi x}-1)}dx$$

or alternatively

$$\varsigma_a'(1,t) = -\frac{1}{2}\frac{\log t}{t} + 2t\int_0^\infty \frac{\tan^{-1}(x/t)e^{\pi x}}{(t^2+x^2)(e^{2\pi x}-1)}dx - \int_0^\infty \frac{x\log(t^2+x^2)e^{\pi x}}{(t^2+x^2)(e^{2\pi x}-1)}dx$$

More generally we have

$$(-1)^n \varsigma_a^{(n)}(s,t) = \frac{t^{-s} \log^n t}{2} + i\int_0^\infty \frac{[(t+ix)^{-s} \log^n(t+ix) - (t-ix)^{-s} \log^n(t-ix)]e^{\pi x}}{e^{2\pi x} - 1} dx$$

and with $s = 1$ we have

$$(-1)^n \varsigma_a^{(n)}(1,t) = \frac{1}{2} \frac{\log^n t}{t} + i\int_0^\infty \frac{[(t+ix)^{-1} \log^n(t+ix) - (t-ix)^{-1} \log^n(t-ix)]e^{\pi x}}{e^{2\pi x} - 1} dx$$

In particular we have

$$\varsigma_a'(1,t) = -\frac{1}{2}\frac{\log t}{t} + i\int_0^\infty \frac{[(t+ix)\log(t-ix) - (t-ix)\log(t+ix)]e^{\pi x}}{(t^2 + x^2)(e^{2\pi x} - 1)} dx$$

and with simple algebra this is equal to

$$= -\frac{1}{2}\frac{\log t}{t} - i\int_0^\infty \frac{[(t-ix)\log(t+ix) - (t+ix)\log(t-ix)]}{(t^2 + x^2)(e^{\pi x} - 1)} dx$$

$$+ i\int_0^\infty \frac{[(t-ix)\log(t+ix) - (t+ix)\log(t-ix)]}{(t^2 + x^2)(e^{2\pi x} - 1)} dx$$

Using the substitution $x = 2y$ we have

$$\int_0^\infty \frac{[(t-ix)\log(t+ix) - (t+ix)\log(t-ix)]}{(t^2 + x^2)(e^{\pi x} - 1)} dx$$

$$= 2\int_0^\infty \frac{[(t-2iy)\log(t+2iy) - (t+2iy)\log(t-2iy)]}{(t^2 + 4y^2)(e^{2\pi y} - 1)} dy$$

$$= 2\int_0^\infty \frac{[(2(t/2) - 2iy)\log(2(t/2) + 2iy) - (2(t/2) + 2iy)\log(2(t/2) - 2iy)]}{4[(t/2)^2 + y^2](e^{2\pi y} - 1)} dy$$

$$= \int_0^\infty \frac{[(t/2) - iy]\log[2(t/2) + 2iy] - [(t/2) + iy]\log[2(t/2) - 2iy]}{[(t/2)^2 + y^2](e^{2\pi y} - 1)} dy$$

$$= \log 2 \int_0^\infty \frac{[(t/2)-iy]-[(t/2)+iy]}{[(t/2)^2+y^2](e^{2\pi y}-1)} dy$$

$$+ \int_0^\infty \frac{[(t/2)-iy]\log[(t/2)+iy]-[(t/2)+iy]\log[(t/2)-iy]}{[(t/2)^2+y^2](e^{2\pi y}-1)} dy$$

$$= -2i \log 2 \int_0^\infty \frac{y}{[(t/2)^2+y^2](e^{2\pi y}-1)} dy - if_1(t/2)$$

where $f_n(t)$ is defined below in (8.7).

The following integral is well known (see for example [46, p.251] and [23a])

(8.6) $$\psi(u) = -\frac{1}{2u} + \log u - 2\int_0^\infty \frac{x}{(u^2+x^2)(e^{2\pi x}-1)} dx$$

and hence we have

$$2\int_0^\infty \frac{y}{[(t/2)^2+y^2](e^{2\pi y}-1)} dy = -\frac{1}{t} + \log(t/2) - \psi(t/2)$$

This results in

$$i\int_0^\infty \frac{[(t-ix)\log(t+ix)-(t+ix)\log(t-ix)]}{(t^2+x^2)(e^{\pi x}-1)} dx = \left[-\frac{1}{t} - \psi(t/2) + \log(t/2)\right]\log 2 + f_1(t/2)$$

We showed in (4.3) of [23a] (see also Coffey's paper [15])

$$\gamma_n(t) = \frac{1}{2}\frac{\log^n t}{t} - \frac{1}{n+1}\log^{n+1} t + i\int_0^\infty \frac{(t-ix)\log^n(t+ix)-(t+ix)\log^n(t-ix)}{(t^2+x^2)(e^{2\pi x}-1)} dx$$

and for convenience we define $f_n(t)$ by

(8.7)
$$f_n(t) = \gamma_n(t) - \frac{1}{2}\frac{\log^n t}{t} + \frac{1}{n+1}\log^{n+1} t = i\int_0^\infty \frac{(t-ix)\log^n(t+ix)-(t+ix)\log^n(t-ix)}{(t^2+x^2)(e^{2\pi x}-1)} dx$$

Therefore we obtain

$$\varsigma_a^{(1)}(1,t) = -\frac{1}{2}\frac{\log t}{t} - \left[-\frac{1}{t} - \psi(t/2) + \log(t/2)\right]\log 2 + f_1(t) - f_1(t/2)$$

or alternatively

(8.8) $\quad \varsigma_a^{(1)}(1,t) = -\frac{1}{2}\frac{\log t}{t} - \left[-\frac{1}{t} - \psi(t/2) + \log(t/2)\right]\log 2 + \gamma_1(t) - \frac{1}{2}\frac{\log t}{t} + \frac{1}{2}\log^2 t$

$$-\gamma_1(t/2) + \frac{\log(t/2)}{t} - \frac{1}{2}\log^2(t/2)$$

With $t=1$ and noting (2.34)

$$\varsigma_a^{(1)}(1) = \gamma \log 2 - \frac{1}{2}\log^2 2$$

we recover (2.23)

$$\gamma_1(1/2) = \gamma_1 - \log^2 2 - 2\gamma \log 2$$

## 9. A FORMULA CONNECTING THE BERNOULLI NUMBERS WITH THE STIELTJES CONSTANTS

Noting that

$$(s-1)\varsigma(s) = \left[\frac{s-1}{1-2^{1-s}}\right][(1-2^{1-s})\varsigma(s)]$$

we see that

$$\frac{d^{m+1}}{ds^{m+1}}[(s-1)\varsigma(s)]\bigg|_{s=1} = \frac{d^{m+1}}{ds^{m+1}}\left[\frac{s-1}{1-2^{1-s}}\right][(1-2^{1-s})\varsigma(s)]\bigg|_{s=1}$$

$$= \frac{d^{m+1}}{ds^{m+1}}\left[\frac{s-1}{1-2^{1-s}}\varsigma_a(s)\right]\bigg|_{s=1}$$

and we apply the Leibniz rule for differentiation to obtain

$$= \sum_{l=0}^{m+1}\binom{m+1}{l}\frac{d^{m+1-l}}{ds^{m+1-l}}\left[\frac{s-1}{1-2^{1-s}}\right]\bigg|_{s=1}\varsigma_a^{(l)}(1)$$

We have

$$\frac{d^{m+1}}{ds^{m+1}}[(s-1)\varsigma(s)]\bigg|_{s=1} = (-1)^m (m+1)\gamma_m$$

and therefore

$$(-1)^m (m+1)\gamma_m = \sum_{l=0}^{m+1} \binom{m+1}{l} \frac{d^{m+1-l}}{ds^{m+1-l}}\left[\frac{s-1}{1-2^{1-s}}\right]\bigg|_{s=1} \varsigma_a^{(l)}(1)$$

The Bernoulli numbers $B_n$ are given by the generating function

$$\frac{t}{e^t - 1} = \sum_{n=0}^{\infty} B_n \frac{t^n}{n!} \qquad , (|t| < 2\pi)$$

and we therefore have

$$\frac{1}{x^{1-s} - 1} = \frac{1}{\exp[-(s-1)\log x] - 1} = -\frac{1}{(s-1)\log x} + \sum_{k=0}^{\infty} \frac{(-1)^k B_{k+1} \log^k x}{(k+1)!} (s-1)^k$$

which holds for $|s-1| < \frac{2\pi}{\log x}$.

Hence we have

$$\frac{s-1}{1-x^{1-s}} = \frac{1}{\log x} + \sum_{k=0}^{\infty} \frac{(-1)^{k+1} B_{k+1} \log^k x}{(k+1)!} (s-1)^{k+1}$$

and differentiation results in

$$\frac{d^{m+1-l}}{ds^{m+1-l}}\left[\frac{s-1}{1-x^{1-s}}\right] = \sum_{k=0}^{\infty} \frac{(-1)^{k+1} B_{k+1} \log^k x}{(k+1)!} (k+1)k \cdots (k-m+l+1)(s-1)^{k-m+l}$$

so that

$$\frac{d^{m+1-l}}{ds^{m+1-l}}\left[\frac{s-1}{1-x^{1-s}}\right]\bigg|_{s=1} = (-1)^{m-l+1} B_{m-l+1} \log^{m-l} x$$

and in particular we have

$$\frac{d^{m+1-l}}{ds^{m+1-l}}\left[\frac{s-1}{1-2^{1-s}}\right]\bigg|_{s=1}=(-1)^{m-l+1}B_{m-l+1}\log^{m-l}2$$

We then see that

$$(9.1)\qquad (m+1)\gamma_m=-\sum_{l=0}^{m+1}\binom{m+1}{l}(-1)^l B_{m-l+1}\,\varsigma_a^{(l)}(1)\log^{m-l}2$$

In 2006 Coffey ([12] and [16a]) showed that

$$\gamma_m=-m!\sum_{l=1}^{m+1}\frac{B_{m-l+1}}{(m-l+1)!}\frac{\log^{m-l}2}{l!}(-1)^l\varsigma_a^{(l)}(1)-\frac{B_{m+1}}{m+1}\log^{m+1}2$$

and using some basic binomial number identities we see that

$$=-\sum_{l=1}^{m+1}\binom{m}{l}\frac{B_{m-l+1}}{m-l+1}\log^{m-l}2(-1)^l\varsigma_a^{(l)}(1)-\frac{B_{m+1}}{m+1}\log^{m+1}2$$

and we therefore obtain

$$(9.2)\qquad \gamma_m=-\frac{1}{m+1}\sum_{l=1}^{m+1}\binom{m+1}{l}B_{m-l+1}\log^{m-l}2(-1)^l\varsigma_a^{(l)}(1)-\frac{B_{m+1}}{m+1}\log^{m+1}2$$

which is a version of the expression originally derived by Liang and Todd [36] in 1972.

Since $\varsigma_a^{(0)}(1)=\log 2$, with the summation starting at $l=0$, we may write this as

$$(9.3)\qquad \gamma_m=-\frac{1}{m+1}\sum_{l=0}^{m+1}\binom{m+1}{l}B_{m-l+1}\log^{m-l}2.(-1)^l\varsigma_a^{(l)}(1)$$

and this is the formula reported by Zhang and Williams [49] in 1994.

We have

$$(9.4)\qquad (-1)^l\varsigma_a^{(l)}(1)=\frac{1}{l+1}\log^{l+1}2-\sum_{k=0}^{l-1}\binom{l}{k}\gamma_k\log^{l-k}2$$

as shown, for example, by Dilcher in [25]. Substituting this in (9.3) however does not appear to produce any additional information

$$\gamma_m = -\frac{1}{m+1}\sum_{l=0}^{m+1}\binom{m+1}{l}B_{m-l+1}\log^{m-l}2\cdot\left[\frac{1}{l+1}\log^{l+1}2-\sum_{k=0}^{l-1}\binom{l}{k}\gamma_k\log^{l-k}2\right]$$

$$=-\frac{1}{m+1}\sum_{l=0}^{m+1}\binom{m+1}{l}\frac{B_{m-l+1}}{l+1}\log^{m+1}2+\frac{1}{m+1}\sum_{l=0}^{m+1}\binom{m+1}{l}B_{m-l+1}\log^m 2\cdot\sum_{k=0}^{l-1}\binom{l}{k}\gamma_k\log^{-k}2$$

For example, with $m=1$ we simply find that $\gamma_1 = \gamma_1$!

## 10. SOME INFINITE SERIES INVOLVING THE STIELTJES CONSTANTS

It is easily seen from the definition of the Hurwitz zeta function that for $\mathrm{Re}(s) > 1$

$$\frac{\partial^n}{\partial x^n}\varsigma(s,1-x) = s(s+1)\cdots(s+n-1)\varsigma(s+n,1-x)$$

and hence, as noted by Ramanujan [38a], the Maclaurin expansion of the Hurwitz zeta function may be written for $|x| < 1$ as

$$\varsigma(s,1-x) = \sum_{n=0}^{\infty}\binom{s+n-1}{n}\varsigma(s+n)x^n$$

Letting $x \to -x$ we have

$$\varsigma(s,1+x) = \sum_{n=0}^{\infty}(-1)^n\binom{s+n-1}{n}\varsigma(s+n)x^n$$

It is easily seen from the definition of the Hurwitz zeta function that

$$\varsigma(s,x+n) = \varsigma(s,x) - \sum_{k=0}^{n-1}\frac{1}{(k+x)^s}$$

and in particular

$$\varsigma(s,x+1) = \varsigma(s,x) - \frac{1}{x^s}$$

Hence we have

$$\varsigma(s,x) = \frac{1}{x^s} + \sum_{n=0}^{\infty}(-1)^n\binom{s+n-1}{n}\varsigma(s+n)x^n$$

Since $\binom{s+n-1}{n} = \frac{\Gamma(s+n)}{\Gamma(n+1)\Gamma(s)}$ this may be written as

$$\varsigma(s,x) = \frac{1}{x^s} + \sum_{n=0}^{\infty}(-1)^n \frac{\Gamma(s+n)}{\Gamma(n+1)\Gamma(s)} \varsigma(s+n) x^n$$

A number of similar series are reported in Chapter 3 of the treatise by Srivastava and Choi [42] where they are primarily expressed in terms of the Pochhammer symbol $(s)_n$ which is defined by

$$(s)_n = s(s+1)\cdots(s+n-1) = \frac{\Gamma(s+n)}{\Gamma(s)}$$

Starting the summation at $n=1$ we have

(10.1) $\quad \varsigma(s,x) - \varsigma(s) = \frac{1}{x^s} + \sum_{n=1}^{\infty}(-1)^n \frac{\Gamma(s+n)}{\Gamma(n+1)\Gamma(s)} \varsigma(s+n) x^n$

and we have the limit

$$\lim_{s \to 1}[\varsigma(s,x) - \varsigma(s)] = \frac{1}{x} + \sum_{n=1}^{\infty}(-1)^n \varsigma(n+1) x^n$$

Reference to (2.17) gives us

$$\lim_{s \to 1}[\varsigma(s,x) - \varsigma(s)] = \gamma_0(x) - \gamma$$

and we then obtain

$$\gamma_0(x) - \gamma = \frac{1}{x} + \sum_{n=1}^{\infty}(-1)^n \varsigma(n+1) x^n$$

Since $\gamma_0(x) = -\psi(x)$ and $\psi(x+1) - \psi(x) = \frac{1}{x}$, the above equation may be written as the well-known Maclaurin expansion of the digamma function [42, p.160]

$$\psi(x+1) = -\gamma + \sum_{k=2}^{\infty}(-1)^k \varsigma(k) x^{k-1}$$

□

From the treatise by Srivastava and Choi [42, p.146] we find a different expression for the Hurwitz zeta function

$$\varsigma(s,x) = \frac{x^{1-s}}{s-1} - \sum_{n=1}^{\infty}(-1)^n \frac{\Gamma(s+n)}{(n+1)\Gamma(n+1)\Gamma(s)} \varsigma(s+n,x)$$

and subtracting $\varsigma(s)$ from both sides gives us

(10.2) $$\varsigma(s,x) - \varsigma(s) = \frac{x^{1-s}}{s-1} - \varsigma(s) - \sum_{n=1}^{\infty}(-1)^n \frac{\Gamma(s+n)}{(n+1)\Gamma(n+1)\Gamma(s)} \varsigma(s+n,x)$$

$$= f(s) - g(s) - \sum_{n=1}^{\infty}(-1)^n \frac{\Gamma(s+n)}{(n+1)\Gamma(n+1)\Gamma(s)} \varsigma(s+n,x)$$

where, for convenience, we have designated $f(s) = \dfrac{x^{1-s}-1}{s-1}$ and $g(s) = \varsigma(s) - \dfrac{1}{s-1}$.

Using L'Hôpital's rule we have the limit as $s \to 1$

$$\lim_{s \to 1} f(s) = -\log x$$

and it is well known that $\lim_{s \to 1} g(s) = \gamma$. Hence we have as $s \to 1$

$$\gamma_0(x) = -\log x - \sum_{n=1}^{\infty} \frac{(-1)^n}{n+1} \varsigma(n+1,x)$$

This concurs with the equation in [42, p.159, Eq (2)] with $t = 1$

$$\sum_{k=2}^{\infty} \frac{(-1)^k t^k}{k} \varsigma(k,x) = \log \Gamma(x+t) - \log \Gamma(x) - t\psi(x)$$

$\square$

Differentiating (10.1) results in

(10.3) $$\varsigma'(s,x) - \varsigma'(s) = -\frac{\log x}{x^s} + \sum_{n=1}^{\infty}(-1)^n \frac{\Gamma(s+n)}{\Gamma(n+1)\Gamma(s)} \varsigma'(s+n)x^n$$

$$+ \sum_{n=1}^{\infty}(-1)^n \frac{\Gamma(s+n)}{\Gamma(n+1)\Gamma(s)} [\psi(s+n) - \psi(s)]\varsigma(s+n)x^n$$

and we have the limit

$$\lim_{s\to 1}[\varsigma'(s,x)-\varsigma'(s)] = -\frac{\log x}{x} + \sum_{n=1}^{\infty}(-1)^n \varsigma'(n+1)x^n + \sum_{n=1}^{\infty}(-1)^n[\psi(n+1)-\psi(1)]\varsigma(n+1)x^n$$

Hence using (2.17) we obtain for $0 < x < 1$

(10.4) $$\gamma_1(x) - \gamma_1 = \frac{\log x}{x} - \sum_{n=1}^{\infty}(-1)^n \varsigma'(n+1)x^n - \sum_{n=1}^{\infty}(-1)^n H_n \varsigma(n+1)x^n$$

where $H_n^{(r)}$ are the harmonic numbers $H_n^{(r)} = \sum_{k=1}^{n}\frac{1}{k^r}$ and $H_n^{(1)} = H_n$.

By differentiating (10.2) we obtain

(10.5) $$\varsigma'(s,x) - \varsigma'(s) = f'(s) - g'(s) - \sum_{n=1}^{\infty}(-1)^n \frac{\Gamma(s+n)}{(n+1)\Gamma(n+1)\Gamma(s)}\varsigma'(s+n,x)$$

$$-\sum_{n=1}^{\infty}(-1)^n \frac{\Gamma(s+n)}{(n+1)\Gamma(n+1)\Gamma(s)}[\psi(s+n)-\psi(s)]\varsigma(s+n,x)$$

We can represent $f(s)$ by the following integral

$$f(s) = \frac{x^{1-s}-1}{s-1} = -\int_{1}^{x} t^{-s} dt$$

so that

$$f^{(p)}(s) = (-1)^{p+1}\int_{1}^{x} t^{-s} \log^p t \, dt$$

and thus

(10.5.1) $$f^{(p)}(1) = (-1)^{p+1}\int_{1}^{x} \frac{\log^p t}{t} dt = (-1)^{p+1}\frac{\log^{p+1} x}{p+1}$$

We have from (1.10)

$$\lim_{s\to 1} g^{(p)}(s) = (-1)^p \gamma_p$$

and therefore, using (2.17), we obtain

(10.6) $$\gamma_1(x) = -\frac{1}{2}\log^2 x + \sum_{n=1}^{\infty}(-1)^n \frac{\varsigma'(n+1,x)}{n+1} + \sum_{n=1}^{\infty}\frac{(-1)^n H_n \varsigma(n+1,x)}{n+1}$$

Coffey [17] has recently shown by a different route that

(10.7)
$$\gamma_1(x) = \log x \log\left(1+\frac{1}{x}\right) - \frac{1}{2}\log^2(1+x) + \sum_{n=1}^{\infty}\frac{(-1)^n}{n+1}\left(\left[\varsigma(n+1,x) - x^{-(n+1)}\right]H_n + \varsigma'(n+1,x)\right)$$

and by integrating the generating function

$$\sum_{n=1}^{\infty}(-1)^n H_n x^n = -\frac{\log(1+x)}{1+x}$$

it is easily seen that this is equivalent to (10.6).

We see from (10.6) that

(10.8)
$$\gamma_1 = \sum_{n=1}^{\infty}(-1)^n \frac{\varsigma'(s+n) + H_n \varsigma(n+1)}{n+1}$$

Equating (10.4) with (10.6) gives us

$$\gamma_1 + \frac{\log x}{x} - \sum_{n=1}^{\infty}(-1)^n \varsigma'(n+1)x^n - \sum_{n=1}^{\infty}(-1)^n H_n \varsigma(n+1)x^n$$

$$= -\frac{1}{2}\log^2 x + \sum_{n=1}^{\infty}(-1)^n \frac{\varsigma'(n+1,x)}{n+1} + \sum_{n=1}^{\infty}\frac{(-1)^n H_n \varsigma(n+1,x)}{n+1}$$

□

Differentiating (10.3) results in

(10.9)
$$\varsigma''(s,x) - \varsigma''(s) = \frac{\log^2 x}{x^s} + \sum_{n=1}^{\infty}(-1)^n \frac{\Gamma(s+n)}{\Gamma(n+1)\Gamma(s)}\varsigma''(s+n)x^n$$

$$+2\sum_{n=1}^{\infty}(-1)^n \frac{\Gamma(s+n)}{\Gamma(n+1)\Gamma(s)}[\psi(s+n) - \psi(s)]\varsigma'(s+n)x^n$$

$$+\sum_{n=1}^{\infty}(-1)^n \frac{\Gamma(s+n)}{\Gamma(n+1)\Gamma(s)}[\psi(s+n) - \psi(s)]^2\varsigma(s+n)x^n$$

$$+\sum_{n=1}^{\infty}(-1)^n \frac{\Gamma(s+n)}{\Gamma(n+1)\Gamma(s)}[\psi'(s+n) - \psi'(s)]\varsigma(s+n)x^n$$

and, as before, we have the limit

$$\lim_{s \to 1}[\varsigma''(s,x) - \varsigma''(s)] = \gamma_2(x) - \gamma_2$$

which results in

(10.10) $\gamma_2(x) - \gamma_2 = \dfrac{\log^2 x}{x} + \sum_{n=1}^{\infty}(-1)^n \varsigma''(n+1) x^n + 2\sum_{n=1}^{\infty}(-1)^n H_n^{(1)} \varsigma'(n+1) x^n$

$$+ \sum_{n=1}^{\infty}(-1)^n \left(\left[H_n^{(1)}\right]^2 - H_n^{(2)}\right) \varsigma(n+1) x^n$$

since $\psi'(n+1) - \psi'(1) = -H_n^{(2)}$.

By differentiating (10.5) we obtain

(10.11) $\varsigma''(s,x) - \varsigma''(s) = f''(s) - g''(s) - \sum_{n=1}^{\infty}(-1)^n \dfrac{\Gamma(s+n)}{(n+1)\Gamma(n+1)\Gamma(s)} \varsigma''(s+n, x)$

$$-2\sum_{n=1}^{\infty}(-1)^n \dfrac{\Gamma(s+n)}{(n+1)\Gamma(n+1)\Gamma(s)} [\psi(s+n) - \psi(s)] \varsigma'(s+n, x)$$

$$-\sum_{n=1}^{\infty}(-1)^n \dfrac{\Gamma(s+n)}{(n+1)\Gamma(n+1)\Gamma(s)} [\psi'(s+n) - \psi'(s)] \varsigma(s+n, x)$$

$$-\sum_{n=1}^{\infty}(-1)^n \dfrac{\Gamma(s+n)}{(n+1)\Gamma(n+1)\Gamma(s)} [\psi(s+n) - \psi(s)]^2 \varsigma(s+n, x)$$

and, as before, we have the limit

(10.12) $\gamma_2(x) = -\dfrac{1}{3}\log^3 x - \sum_{n=1}^{\infty}\dfrac{(-1)^n}{n+1} \varsigma''(n+1, x) - 2\sum_{n=1}^{\infty}\dfrac{(-1)^n}{n+1} H_n^{(1)} \varsigma'(n+1, x)$

$$-\sum_{n=1}^{\infty}\dfrac{(-1)^n}{n+1}\left(\left[H_n^{(1)}\right]^2 - H_n^{(2)}\right) \varsigma(n+1, x)$$

The corresponding result derived by Coffey [17] is (where I have corrected a misprint)

(10.13)
$\gamma_2(x) = \dfrac{\log^2 x}{x} - \dfrac{1}{3}\log^3(1+x) - \sum_{n=1}^{\infty}\dfrac{(-1)^n}{n+1} \varsigma''(n+1, x+1) + 2\sum_{n=1}^{\infty}\dfrac{s(n+1, 2)}{(n+1)n!} \varsigma'(n+1, x+1)$

$$-2\sum_{n=1}^{\infty}\frac{s(n+1,3)}{(n+1)n!}\varsigma(n+1,x+1)$$

where $s(n,k)$ are the Stirling numbers of the first kind [42, p.56]. In particular, we have

$$s(n+1,2) = (-1)^{n+1} n! H_n^{(1)}$$

$$s(n+1,3) = \frac{1}{2}(-1)^n n!\left(\left[H_n^{(1)}\right]^2 - H_n^{(2)}\right)$$

and we therefore have

(10.14)

$$\gamma_2(x) = \frac{\log^2 x}{x} - \frac{1}{3}\log^3(1+x) - \sum_{n=1}^{\infty}\frac{(-1)^n}{n+1}\varsigma''(n+1,x+1) - 2\sum_{n=1}^{\infty}\frac{(-1)^n}{n+1}H_n^{(1)}\varsigma'(n+1,x+1)$$

$$-\sum_{n=1}^{\infty}\frac{(-1)^n}{n+1}\left(\left[H_n^{(1)}\right]^2 - H_n^{(2)}\right)\varsigma(n+1,x+1)$$

With $x=1$ in (10.14) we obtain

$$\gamma_2 = -\frac{1}{3}\log^3 2 - \sum_{n=1}^{\infty}\frac{(-1)^n}{n+1}\varsigma''(n+1) - 2\sum_{n=1}^{\infty}\frac{(-1)^n}{n+1}H_n^{(1)}\varsigma'(n+1)$$

$$-\sum_{n=1}^{\infty}\frac{(-1)^n}{n+1}\left(\left[H_n^{(1)}\right]^2 - H_n^{(2)}\right)\varsigma(n+1) + \sum_{n=1}^{\infty}\frac{(-1)^n}{n+1}\left(\left[H_n^{(1)}\right]^2 - H_n^{(2)}\right)$$

where we have used

$$\varsigma(s,x+1) = \varsigma(s,x) - \frac{1}{x^s} \qquad \varsigma^{(p)}(s,2) = \varsigma^{(p)}(s,1) = \varsigma^{(p)}(s)$$

The final summation is

$$\sum_{n=1}^{\infty}\frac{(-1)^n}{n+1}\left(\left[H_n^{(1)}\right]^2 - H_n^{(2)}\right) = 2\sum_{n=1}^{\infty}\frac{s(n+1,3)}{(n+1)!}$$

We employ the Maclaurin expansion due to Cauchy

$$\log^k(1+x) = k!\sum_{n=k}^{\infty} s(n,k)\frac{x^n}{n!}$$

and since $s(n,k) = \frac{1}{k!}\frac{d^n}{dx^n}\log^k(1+x)\Big|_{x=0}$ it is clear that $s(n,k) = 0 \ \forall \ n \leq k-1$.

Reindexing gives us

$$\log^k(1+x) = k! \sum_{m=k-1}^{\infty} s(m+1,k)\frac{x^{m+1}}{(m+1)!}$$

and with $x=1$ we see that

$$2\sum_{n=1}^{\infty} \frac{s(n+1,3)}{(n+1)!} = \frac{1}{3}\log^3 2$$

We then deduce that

$$\gamma_2 = -\sum_{n=1}^{\infty} \frac{(-1)^n}{n+1}\varsigma''(n+1) - 2\sum_{n=1}^{\infty} \frac{(-1)^n}{n+1} H_n^{(1)} \varsigma'(n+1) - \sum_{n=1}^{\infty} \frac{(-1)^n}{n+1}\left(\left[H_n^{(1)}\right]^2 - H_n^{(2)}\right)\varsigma(n+1)$$

which concurs with $x=1$ in (10.12).

Higher derivatives of (10.11) will result in a family of similar infinite series involving the (exponential) complete Bell polynomials (see [23c] for a brief review of these polynomials in a similar application).

□

From [42, p.96] we have

$$\varsigma(n+1) = (-1)^n \int_0^1 \frac{\log^n y}{1-y} dy$$

and the following integral representation of the harmonic numbers [23c] is well known

$$H_n = -n\int_0^1 (1-t)^{n-1} \log t \, dt$$

We then obtain

$$\sum_{n=1}^{\infty} (-1)^n H_n \varsigma(n+1) x^n = -\sum_{n=1}^{\infty} nx^n \int_0^1 (1-t)^{n-1} \log t \, dt \int_0^1 \frac{\log^n y}{1-y} dy$$

$$= -\int_0^1\int_0^1 \sum_{n=1}^{\infty} \frac{n[x(1-t)\log y]^n \log t}{(1-t)(1-y)} dt \, dy$$

Noting the derivative of the geometric series we have

$$\sum_{n=1}^{\infty}(-1)^n H_n \varsigma(n+1)x^n = -\int_0^1\int_0^1 \frac{\log t}{(1-t)(1-y)[1-x(1-t)\log y]^2} \, dt \, dy$$

but this double integral appears to be rather intractable.

□

Alternatively, we may employ Wilton's formula [42, p.144] for $|x| < |a|$ and $\text{Re}(a) > 0$

$$\varsigma(s, a-x) = \sum_{n=0}^{\infty} \frac{\Gamma(s+n)}{\Gamma(n+1)\Gamma(s)} \varsigma(s+n, a) x^n$$

which is valid for all values of $s \neq 1$. We may write this as

(10.15) $\qquad \varsigma(s, a-x) - \varsigma(s, a) = \sum_{n=1}^{\infty} \frac{\Gamma(s+n)}{\Gamma(n+1)\Gamma(s)} \varsigma(s+n, a) x^n$

and we have the limit

$$\gamma_0(a-x) - \gamma_0(a) = \sum_{n=1}^{\infty} \varsigma(n+1, a) x^n$$

or equivalently [42, p.144]

$$\psi(a) - \psi(a-x) = \sum_{n=1}^{\infty} \varsigma(n+1, a) x^n = \sum_{n=2}^{\infty} \varsigma(n, a) x^{n-1}$$

Differentiating (10.15) results in

$$\varsigma'(s, a-x) - \varsigma'(s, a) = \sum_{n=1}^{\infty} \frac{\Gamma(s+n)}{\Gamma(n+1)\Gamma(s)} \{\varsigma'(s+n, a) + [\psi(s+n) - \psi(s)]\varsigma(s+n, a)\} x^n$$

and we have the limit (as recently reported by Coffey [17a])

(10.16) $\qquad \gamma_1(a-x) - \gamma_1(a) = -\sum_{n=1}^{\infty}\left[\varsigma'(n+1, a) + H_n^{(1)}\varsigma(n+1, a)\right] x^n$

With $a = 1$ and $x \to -x$ this becomes

(10.17) $$\gamma_1(1+x) - \gamma_1 = -\sum_{n=1}^{\infty}(-1)^n\left[\varsigma'(n+1) + H_n^{(1)}\varsigma(n+1)\right]x^n$$

which, as shown below, is equivalent to (10.4).

Since $\varsigma(s, 1+x) - \varsigma(s,x) = -\dfrac{1}{x^s}$ we have

$$\varsigma^{(p)}(s,1+x) - \varsigma^{(p)}(s,x) = (-1)^{p+1}\dfrac{\log^p x}{x^s}$$

and therefore using (2.17) we obtain

$$\gamma_p(1+x) - \gamma_p(x) = -\dfrac{\log^p x}{x}$$

With $a \to 1+a$ and $x \to -x$ (10.16) becomes

$$\gamma_1(1+a+x) - \gamma_1(1+a) = -\sum_{n=1}^{\infty}(-1)^n\left[\varsigma'(n+1, 1+a) + H_n^{(1)}\varsigma(n+1, 1+a)\right]x^n$$

Wilton's formula was also used by Coffey in [12], [17] and [17a].

$\square$

From [42, p.147] we have for $|x| < 1$

(10.18) $$\varsigma(s-1, 1-x) - \varsigma(s-1) = \sum_{n=0}^{\infty}\binom{s+n-1}{n+1}\varsigma(s+n)x^{n+1}$$

and since $\binom{s+n-1}{n+1} = \dfrac{s-1}{n+1}\binom{s+n-1}{n}$ this may be written as

$$\varsigma(s-1, 1-x) - \varsigma(s-1) = (s-1)\sum_{n=0}^{\infty}\dfrac{1}{n+1}\binom{s+n-1}{n}\varsigma(s+n)x^{n+1}$$

or equivalently

(10.19) $$\varsigma(s-1, 1-x) - \varsigma(s-1) = (s-1)\sum_{n=0}^{\infty}\dfrac{1}{n+1}\dfrac{\Gamma(s+n)}{\Gamma(n+1)\Gamma(s)}\varsigma(s+n)x^{n+1}$$

Differentiation gives us

$$\varsigma'(s-1,1-x) - \varsigma'(s-1)$$

$$= (s-1)\sum_{n=0}^{\infty} \frac{1}{n+1} \frac{\Gamma(s+n)}{\Gamma(n+1)\Gamma(s)} \left[\varsigma'(s+n) + [\psi(s+n) - \psi(s)]\varsigma(s+n)\right] x^{n+1}$$

$$+ \sum_{n=0}^{\infty} \frac{1}{n+1} \frac{\Gamma(s+n)}{\Gamma(n+1)\Gamma(s)} \varsigma(s+n) x^{n+1}$$

which, upon letting $s = 2$, becomes

$$\varsigma'(1,1-x) - \varsigma'(1) = \sum_{n=0}^{\infty} \left[\varsigma'(n+2) + [\psi(n+2) - \psi(2) + 1]\varsigma(n+2)\right] x^{n+1}$$

$$= \sum_{n=0}^{\infty} \left[\varsigma'(n+2) + [\psi(n+2) - \psi(1) + \psi(1) - \psi(2) + 1]\varsigma(n+2)\right] x^{n+1}$$

$$= \sum_{n=0}^{\infty} \left[\varsigma'(n+2) + [\psi(n+2) - \psi(1)]\varsigma(n+2)\right] x^{n+1}$$

$$= \sum_{n=0}^{\infty} \left[\varsigma'(n+2) + H_{n+1}^{(1)} \varsigma(n+2)\right] x^{n+1}$$

and reindexing gives us

$$\varsigma'(1,1-x) - \varsigma'(1) = \sum_{n=1}^{\infty} \left[\varsigma'(n+1) + H_n^{(1)} \varsigma(n+1)\right] x^n$$

or equivalently

$$\gamma_1(1-x) - \gamma_1 = -\sum_{n=1}^{\infty} \left[\varsigma'(n+1) + H_n^{(1)} \varsigma(n+1)\right] x^n$$

It is easily seen that this corresponds with (10.17). Using (2.33) for $\gamma_1(1/2)$ we find that

(10.20) $$\sum_{n=1}^{\infty} \frac{1}{2^n} \left[\varsigma'(n+1) + H_n^{(1)} \varsigma(n+1)\right] = \log^2 2 + 2\gamma \log 2$$

□

Equation (10.19) may be written as (where we have started the summation at $n = 1$)

$$\varsigma(s-1,1-x)-\varsigma(s-1)-(s-1)\varsigma(s)x = (s-1)\sum_{n=1}^{\infty}\frac{1}{n+1}\frac{\Gamma(s+n)}{\Gamma(n+1)\Gamma(s)}\varsigma(s+n)x^{n+1}$$

and taking the limit as $s \to 1$ we have

$$\varsigma(0,1-x)-\varsigma(0)-x = 0$$

since $\lim_{s \to 1}[(s-1)\varsigma(s)] = 1$. This is in accordance with the well-known formulae:

$$\varsigma(0,t) = \frac{1}{2} - t$$

$$\varsigma(0,1) = \varsigma(0) = -\frac{1}{2}$$

We also have

$$\varsigma'(s-1,1-x)-\varsigma'(s-1)-x\frac{d}{ds}[(s-1)\varsigma(s)]$$

$$= (s-1)\sum_{n=1}^{\infty}\frac{1}{n+1}\frac{\Gamma(s+n)}{\Gamma(n+1)\Gamma(s)}\left[\varsigma'(s+n)+[\psi(s+n)-\psi(s)]\varsigma(s+n)\right]x^{n+1}$$

$$+\sum_{n=1}^{\infty}\frac{1}{n+1}\frac{\Gamma(s+n)}{\Gamma(n+1)\Gamma(s)}\varsigma(s+n)x^{n+1}$$

which becomes upon letting $s = 1$

$$\varsigma'(0,1-x)-\varsigma'(0)-\gamma x = \sum_{n=1}^{\infty}\frac{\varsigma(n+1)}{n+1}x^{n+1}$$

Using Lerch's identity (2.5) this may be expressed as the familiar Maclaurin series which is valid for $-1 \leq x < 1$

$$\log\Gamma(1-x)-\gamma x = \sum_{n=1}^{\infty}\frac{\varsigma(n+1)}{n+1}x^{n+1}$$

The second derivative gives us

$$\varsigma''(s-1,1-x)-\varsigma''(s-1)-x\frac{d^2}{ds^2}[(s-1)\varsigma(s)]$$

$$= (s-1)\sum_{n=1}^{\infty} \frac{1}{n+1} \frac{\Gamma(s+n)}{\Gamma(n+1)\Gamma(s)} \varsigma''(s+n)x^{n+1}$$

$$+2(s-1)\sum_{n=1}^{\infty} \frac{1}{n+1} \frac{\Gamma(s+n)}{\Gamma(n+1)\Gamma(s)} [\psi(s+n)-\psi(s)]\varsigma'(s+n)x^{n+1}$$

$$+2\sum_{n=1}^{\infty} \frac{1}{n+1} \frac{\Gamma(s+n)}{\Gamma(n+1)\Gamma(s)} \varsigma'(s+n)x^{n+1}$$

$$+(s-1)\sum_{n=1}^{\infty} \frac{1}{n+1} \frac{\Gamma(s+n)}{\Gamma(n+1)\Gamma(s)} [\psi'(s+n)-\psi'(s)]\varsigma(s+n)x^{n+1}$$

$$+(s-1)\sum_{n=1}^{\infty} \frac{1}{n+1} \frac{\Gamma(s+n)}{\Gamma(n+1)\Gamma(s)} [\psi(s+n)-\psi(s)]^2 \varsigma(s+n)x^{n+1}$$

$$+2\sum_{n=1}^{\infty} \frac{1}{n+1} \frac{\Gamma(s+n)}{\Gamma(n+1)\Gamma(s)} [\psi(s+n)-\psi(s)]\varsigma(s+n)x^{n+1}$$

which becomes upon letting $s=1$

$$\varsigma''(0,1-x)-\varsigma''(0)+2\gamma_2 x = 2\sum_{n=1}^{\infty} \frac{1}{n+1}\left[\varsigma'(n+1)+H_n^{(1)}\varsigma(n+1)\right]x^{n+1}$$

Letting $x \to -x$ we have

$$\varsigma''(0,1+x)-\varsigma''(0)-2\gamma_2 x = -2\sum_{n=1}^{\infty} \frac{(-1)^n}{n+1}\left[\varsigma'(n+1)+H_n^{(1)}\varsigma(n+1)\right]x^{n+1}$$

or alternatively

(10.21)    $$\varsigma''(0,x)-\varsigma''(0)-\log^2 x - 2\gamma_2 x = -2\sum_{n=1}^{\infty} \frac{(-1)^n}{n+1}\left[\varsigma'(n+1)+H_n^{(1)}\varsigma(n+1)\right]x^{n+1}$$

□

Integrating (10.4) results in

$$\frac{1}{2}[\varsigma''(0,t)-\varsigma''(0)]-\gamma_1(t-1) = \frac{1}{2}\log^2 t - \sum_{n=1}^{\infty} \frac{(-1)^n}{n+1}\left(\varsigma'(n+1)+H_n\varsigma(n+1)\right)[t^{n+1}-1]$$

where we have used [21]

$$\int_1^t \gamma_n(x)dx = \frac{(-1)^{n+1}}{n+1}\left[\varsigma^{(n+1)}(0,t)-\varsigma^{(n+1)}(0)\right]$$

Using (10.8) this may be written as

(10.22) $\quad \dfrac{1}{2}[\varsigma''(0,t)-\varsigma''(0)]-\gamma_1 t = \dfrac{1}{2}\log^2 t - \sum_{n=1}^{\infty}\dfrac{(-1)^n}{n+1}[\varsigma'(n+1)+H_n\varsigma(n+1)]t^{n+1}$

and we see that this simply corresponds with (10.21).

Referring back to (5.4)

$$\varsigma''\left(0,\frac{1}{2}\right) = -\log(2\pi)\log 2 - \frac{1}{2}\log^2 2$$

and using the result previously obtained by Ramanujan [7b] and Apostol [5a] (see also the recent paper [23d])

$$\varsigma''(0) = \gamma_1 + \frac{1}{2}\gamma^2 - \frac{1}{24}\pi^2 - \frac{1}{2}\log^2(2\pi)$$

upon letting $t=1/2$ in (10.22) we obtain

(10.23)

$$\sum_{n=1}^{\infty}\frac{(-1)^n}{(n+1)2^n}[\varsigma'(n+1)+H_n\varsigma(n+1)] = 2\gamma_1 + \frac{1}{2}\gamma^2 - \frac{1}{24}\pi^2 - \frac{1}{2}\log^2(2\pi) + \log(2\pi)\log 2 + \frac{3}{2}\log^2 2$$

## 11. SOME CONNECTIONS WITH RAMANUJAN'S WORK

We now recall (1.18)

$$\gamma_1(x) - \gamma_1(1) = \sum_{n=0}^{\infty}\left[\frac{\log(n+x)}{n+x} - \frac{\log(n+1)}{n+1}\right]$$

and noting that Ramanujan [7a, p.197] defined $\varphi(x)$ by

(11.1) $\qquad\qquad \varphi(x) = \sum_{n=1}^{\infty}\left[\frac{\log n}{n} - \frac{\log(n+x)}{n+x}\right]$

we see that

$$\gamma_1(x) - \gamma_1(1) = -\sum_{n=0}^{\infty} \left[ \frac{\log(n+1)}{n+1} - \frac{\log(n+x)}{n+x} \right]$$

$$= -\sum_{n=1}^{\infty} \left[ \frac{\log(n+1)}{n+1} - \frac{\log(n+x)}{n+x} \right] + \frac{\log x}{x}$$

$$= \frac{\log x}{x} - \sum_{n=1}^{\infty} \left[ \frac{\log n}{n} - \frac{\log(n+x)}{n+x} + \frac{\log(n+1)}{n+1} - \frac{\log n}{n} \right]$$

We then have

(11.2) $$\gamma_1(x) - \gamma_1(1) = \frac{\log x}{x} - \sum_{n=1}^{\infty} \left[ \frac{\log n}{n} - \frac{\log(n+x)}{n+x} \right]$$

and thus

(11.3) $$\gamma_1(x) - \gamma_1 = \frac{\log x}{x} - \varphi(x)$$

We note that $\varphi(1) = 0$. It is not immediately obvious how this representation may be used to define $\varphi(x)$ for $x \leq 0$ but referring to (10.4)

$$\gamma_1(x) - \gamma_1 = \frac{\log x}{x} - \sum_{n=1}^{\infty} (-1)^n \varsigma'(n+1) x^n - \sum_{n=1}^{\infty} (-1)^n H_n \varsigma(n+1) x^n$$

it seems natural to write

(11.4) $$\varphi(x) = \sum_{n=1}^{\infty} (-1)^n \varsigma'(n+1) x^n + \sum_{n=1}^{\infty} (-1)^n H_n \varsigma(n+1) x^n$$

whereupon we have for negative values

$$\varphi(-x) = \sum_{n=1}^{\infty} \varsigma'(n+1) x^n + \sum_{n=1}^{\infty} H_n \varsigma(n+1) x^n$$

Letting $x \to -x$ in (10.17) gives us

$$\gamma_1(1-x) - \gamma_1 = -\left[ \sum_{n=1}^{\infty} \varsigma'(n+1) x^n + \sum_{n=1}^{\infty} H_n \varsigma(n+1) x^n \right]$$

and we therefore have

(11.5) $$\gamma_1(1-x) - \gamma_1 = -\varphi(-x)$$

With $x = 1/2$ we have

$$\gamma_1\left(\frac{1}{2}\right) - \gamma_1 = -\left[\sum_{n=1}^{\infty}\frac{\varsigma'(n+1)}{2^2} + \sum_{n=1}^{\infty}\frac{H_n\varsigma(n+1)}{2^2}\right]$$

Using (2.23)

$$\gamma_1\left(\frac{1}{2}\right) = \gamma_1 - \log^2 2 - 2\gamma \log 2$$

we see that

(11.6) $$\log^2 2 + 2\gamma \log 2 = \left[\sum_{n=1}^{\infty}\frac{\varsigma'(n+1)}{2^2} + \sum_{n=1}^{\infty}\frac{H_n\varsigma(n+1)}{2^2}\right]$$

which agrees with Ramanujan's result [7a, p.199]

$$\varphi\left(-\frac{1}{2}\right) = \log^2 2 + 2\gamma \log 2$$

$\square$

Letting $-x \to x - 1$ in (11.5) we obtain

$$\gamma_1(x) - \gamma_1 = -\varphi(x-1)$$

and we then deduce that (see [7a, p.200])

(11.7) $$\varphi(x-1) - \varphi(-x) = \gamma_1(1-x) - \gamma_1(x)$$

(11.8) $$\varphi(x-1) + \varphi(-x) = -[\gamma_1(1-x) + \gamma_1(x)] + 2\gamma_1$$

For example, we have with $x = 1/4$

$$\varphi(-3/4) + \varphi(-1/4) = -[\gamma_1(3/4) + \gamma_1(1/4)] + 2\gamma_1$$

and referring to (2.27) and (2.28) we see that

$$\varphi(-3/4) + \varphi(-1/4) = 7\log^2 2 + 6\gamma \log 2$$

which concurs with Ramanujan's result [7a, p.199].

We also have

$$\varphi(-1/4) - \varphi(-3/4) = \gamma_1(1/4) - \gamma_1(3/4)$$

and referring to (2.27) and (2.28) we obtain

(11.9) $\quad \varphi(-1/4) - \varphi(-3/4) = -\pi\left[\gamma + 4\log 2 + 3\log \pi - 4\log \Gamma\left(\frac{1}{4}\right)\right]$

Using the definition (11.1) Ramanujan [7a, p.198] showed that

$$\frac{1}{4}[\varphi(-1/4) - \varphi(-3/4)] = \sum_{n=1}^{\infty}\left[\frac{\log(n-3/4)}{4n-3} - \frac{\log(n-1/4)}{4n-1}\right]$$

$$= \sum_{n=0}^{\infty}\frac{(-1)^n \log(2n+1)}{2n+1} - \log 4 \sum_{n=0}^{\infty}\frac{(-1)^n}{2n+1}$$

$$= \sum_{n=0}^{\infty}\frac{(-1)^n \log(2n+1)}{2n+1} - \frac{1}{2}\pi \log 2$$

and using (5.6)

$$\sum_{n=0}^{\infty}\frac{(-1)^n \log(2n+1)}{2n+1} = -\left[\gamma + 2\log 2 + 3\log \pi - 4\log \Gamma\left(\frac{1}{4}\right)\right]$$

this confirms (11.9) above.

$\square$

Ramanujan [7a, p.220] also considered a generalised function $\varphi_r(x)$ defined by

(11.10) $\quad \varphi_r(x) = \sum_{k=1}^{\infty}\left[\frac{\log k}{k^r} - \frac{\log(k+x)}{(k+x)^r}\right]$

and we have

$$\varphi_r(x) = \sum_{k=1}^{\infty}\log k\left[\frac{1}{k^r} - \frac{1}{(k+x)^r}\right] - \sum_{k=1}^{\infty}\frac{\log(1+x/k)}{(k+x)^r}$$

$$= S_1 + S_2$$

For $|x| < 1$ we have

$$S_1 = \sum_{k=1}^{\infty} \log k \left[ \frac{1}{k^r} - \frac{1}{k^r} \sum_{n=0}^{\infty} \binom{-r}{n} \left(\frac{x}{k}\right)^n \right]$$

$$= -\sum_{k=1}^{\infty} \frac{\log k}{k^r} \sum_{n=1}^{\infty} \binom{-r}{n} \left(\frac{x}{k}\right)^n$$

$$= -\sum_{k=1}^{\infty} \frac{\log k}{k^{r+n}} \sum_{n=1}^{\infty} \binom{-r}{n} x^n$$

$$= \sum_{n=1}^{\infty} \binom{-r}{n} \varsigma'(n+r) x^n$$

(in the second line we have inserted a minus sign which has inadvertently been missed in Berndt's book [7a, p.220]).

Secondly, for $|x| < 1$ we have

$$S_2 = -\sum_{k=1}^{\infty} \frac{\log(1+x/k)}{(k+x)^r}$$

$$= -\sum_{k=1}^{\infty} \frac{1}{k^r} \frac{\log(1+x/k)}{(1+x/k)^r}$$

$$= -\sum_{k=1}^{\infty} \frac{1}{k^r} \sum_{j=0}^{\infty} \binom{-r}{j} \left(\frac{x}{k}\right)^j \sum_{m=1}^{\infty} \frac{(-1)^{m+1}}{m} \left(\frac{x}{k}\right)^m$$

$$= \sum_{k=1}^{\infty} \frac{1}{k^r} \sum_{n=1}^{\infty} (-1)^n \sum_{j=0}^{n-1} \frac{(r)_j}{j!(n-j)!} \left(\frac{x}{k}\right)^n$$

$$= \sum_{k=1}^{\infty} \frac{1}{k^r} \sum_{n=1}^{\infty} (-1)^n \sum_{j=0}^{n-1} \binom{n}{j} \frac{(r)_j}{n!} \left(\frac{x}{k}\right)^n$$

Using the following identity from Hansen's tables [31, p.126] (this corrects the misprint in [7a, p.281])

$$\sum_{j=0}^{n-1} \frac{(r)_j}{j!(n-j)!} = \frac{(r)_n}{n!} \sum_{j=0}^{n-1} \frac{1}{k+r}$$

or equivalently

$$(11.10.1) \qquad \sum_{j=0}^{n-1} \binom{n}{j} \frac{(r)_j}{n!} = \frac{(r)_n}{n!} H_n^{(1)}(r)$$

we see that

$$S_2 = \sum_{k=1}^{\infty} \frac{1}{k^r} \sum_{n=1}^{\infty} \frac{(-1)^n H_n(r)}{n!} \left(\frac{x}{k}\right)^n$$

$$= \sum_{k=1}^{\infty} \frac{1}{k^r} \binom{-r}{n} H_n(r) \varsigma(n+r) x^n$$

Hence we have for $r > 0$ and $|x| < 1$

$$(11.11) \qquad \varphi_r(x) = \sum_{n=1}^{\infty} \binom{-r}{n} [H_n(r) \varsigma(n+r) + \varsigma'(n+r)] x^n$$

where we have removed the minus sign which appears in Berndt's book [7a, p.220].

Here $H_n(r)$ is defined by

$$H_n^{(m)}(r) = \sum_{k=0}^{n-1} \frac{1}{(k+r)^m} \qquad H_n^{(1)}(r) = H_n(r) = \sum_{k=0}^{n-1} \frac{1}{k+r}$$

so that $H_n(1) = H_n$ corresponds with the familiar harmonic numbers (it should be noted that the notation employed here slightly differs from that employed by Berndt and Ramanujan).

Equation (11.11) may be written as

$$(11.12) \qquad \varphi_r(x) = \sum_{n=1}^{\infty} (-1)^n \frac{\Gamma(r+n)}{\Gamma(n+1)\Gamma(r)} [H_n(r) \varsigma(n+r) + \varsigma'(n+r)] x^n$$

and in particular for $r = 1$ we have

$$(11.13) \qquad \varphi_1(x) = \varphi(x) = \sum_{n=1}^{\infty} (-1)^n [H_n \varsigma(n+1) + \varsigma'(n+1)] x^n$$

which concurs with (11.4) above.

Differentiation of (11.12) results in

$$\frac{\partial}{\partial r}\varphi_r(x) = \sum_{n=1}^{\infty}(-1)^n \frac{\Gamma(r+n)}{\Gamma(n+1)\Gamma(r)}[H_n(r)\varsigma'(n+r) - H_n^{(2)}(r)\varsigma(n+r) + \varsigma''(n+r)]x^n$$

$$+ \sum_{n=1}^{\infty}(-1)^n \frac{\Gamma(r+n)}{\Gamma(n+1)\Gamma(r)}[\psi(r+n) - \psi(r)][H_n(r)\varsigma(n+r) + \varsigma'(n+r)]x^n$$

and with $r = 1$ we have

$$\left.\frac{\partial}{\partial r}\varphi_r(x)\right|_{r=1} = \sum_{n=1}^{\infty}(-1)^n[H_n^{(1)}\varsigma'(n+1) - H_n^{(2)}\varsigma(n+1) + \varsigma''(n+1)]x^n$$

$$+ \sum_{n=1}^{\infty}(-1)^n H_n^{(1)}[H_n^{(1)}\varsigma(n+1) + \varsigma'(n+1)]x^n$$

(11.14) $$= \sum_{n=1}^{\infty}(-1)^n[2H_n^{(1)}\varsigma'(n+1) + \left(\left[H_n^{(1)}\right]^2 - H_n^{(2)}\right)\varsigma(n+1) + \varsigma''(n+1)]x^n$$

From its definition (11.1) we see that

$$\frac{\partial}{\partial r}\varphi_r(x) = -\sum_{n=1}^{\infty}\left[\frac{\log^2 n}{n^r} - \frac{\log^2(n+x)}{(n+x)^r}\right]$$

and in particular

$$\left.\frac{\partial}{\partial r}\varphi_r(x)\right|_{r=1} = -\sum_{n=1}^{\infty}\left[\frac{\log^2 n}{n} - \frac{\log^2(n+x)}{n+x}\right]$$

We see that

$$\gamma_r(x) - \gamma_r(1) = -\sum_{n=0}^{\infty}\left[\frac{\log^r(n+1)}{n+1} - \frac{\log^r(n+x)}{n+x}\right]$$

$$= -\sum_{n=1}^{\infty}\left[\frac{\log^r(n+1)}{n+1} - \frac{\log^r(n+x)}{n+x}\right] + \frac{\log^r x}{x}$$

$$= \frac{\log^r x}{x} - \sum_{n=1}^{\infty}\left[\frac{\log^r n}{n} - \frac{\log^r(n+x)}{n+x} + \frac{\log^r(n+1)}{n+1} - \frac{\log^r n}{n}\right]$$

We then have

(11.15) $$\gamma_r(x) - \gamma_r(1) = \frac{\log^r x}{x} - \sum_{n=1}^{\infty}\left[\frac{\log^r n}{n} - \frac{\log^r(n+x)}{n+x}\right]$$

and in particular

$$\gamma_2(x) - \gamma_2(1) = \frac{\log^2 x}{x} - \sum_{n=1}^{\infty}\left[\frac{\log^2 n}{n} - \frac{\log^2(n+x)}{n+x}\right]$$

We then have

$$\gamma_2(x) - \gamma_2(1) = \frac{\log^2 x}{x} + \frac{\partial}{\partial r}\varphi_r(x)\bigg|_{r=1}$$

and referring to (11.14) gives us

(11.16)
$$\gamma_2(x) - \gamma_2(1) = \frac{\log^2 x}{x} + \sum_{n=1}^{\infty}(-1)^n[2H_n^{(1)}\varsigma'(n+1) + \left(\left[H_n^{(1)}\right]^2 - H_n^{(2)}\right)\varsigma(n+1) + \varsigma''(n+1)]x^n$$

which corresponds with (10.10).

Upon differentiating (11.15) we obtain

$$\gamma_r'(x) = \frac{d}{dx}\frac{\log^r x}{x} + \sum_{n=1}^{\infty}\left[\frac{r\log^{r-1}(n+x)}{(n+x)^2} - \frac{\log^r(n+x)}{(n+x)^2}\right]$$

and in particular we have

$$\gamma_1'(x) = \frac{d}{dx}\frac{\log x}{x} + \sum_{n=1}^{\infty}\left[\frac{1}{(n+x)^2} - \frac{\log(n+x)}{(n+x)^2}\right]$$

We have from (11.1)

$$\varphi_1(x) = \sum_{n=1}^{\infty}\left[\frac{\log n}{n} - \frac{\log(n+x)}{n+x}\right]$$

and therefore

$$\varphi_1'(x) = -\sum_{n=1}^{\infty}\left[\frac{1 - \log(n+x)}{(n+x)^2}\right]$$

and hence we have

$$\gamma_1'(x) = \frac{d}{dx}\frac{\log x}{x} - \varphi_1'(x)$$

Integration results in

$$\gamma_1(x) - \gamma_1(1) = \frac{\log x}{x} - \varphi_1(x) + \varphi_1(1)$$

and since $\varphi_1(1) = 0$ we have

$$\gamma_1(x) - \gamma_1(1) = \frac{\log x}{x} - \varphi_1(x)$$

Substituting (11.13) we rediscover (10.4).

Integrating (11.1) gives us

$$\int_0^t \varphi_1(x)dx = \sum_{n=1}^{\infty}\left[t\frac{\log n}{n} + \frac{1}{2}\log^2 n - \frac{1}{2}\log^2(n+t)\right]$$

and therefore integrating (11.13) shows that

(11.17) $$\sum_{n=1}^{\infty}\frac{(-1)^n}{n+1}[H_n\varsigma(n+1) + \varsigma'(n+1)]t^{n+1} = \sum_{n=1}^{\infty}\left[t\frac{\log n}{n} + \frac{1}{2}\log^2 n - \frac{1}{2}\log^2(n+t)\right]$$

With $t = 1/2$ we have using (10.23)

(11.18)

$$\sum_{n=1}^{\infty}\left[\frac{\log n}{n} + \log^2 n - \log^2(n+1/2)\right] = 2\gamma_1 + \frac{1}{2}\gamma^2 - \frac{1}{24}\pi^2 - \frac{1}{2}\log^2(2\pi) + \log(2\pi)\log 2 + \frac{3}{2}\log^2 2$$

Differentiating (11.17) gives us

$$\sum_{n=1}^{\infty}(-1)^n[H_n\varsigma(n+1) + \varsigma'(n+1)]t^n = \sum_{n=1}^{\infty}\left[\frac{\log n}{n} - \frac{\log(n+t)}{n+t}\right]$$

Dilcher [25a] has defined generalised gamma functions by

$$\Gamma_k(x) = \lim_{n\to\infty} \frac{\exp\left[\dfrac{x}{k+1}\log^{k+1} n\right]\prod_{j=1}^{n}\exp\left[\dfrac{1}{k+1}\log^{k+1} j\right]}{\prod_{j=0}^{n}\exp\left[\dfrac{1}{k+1}\log^{k+1}(j+x)\right]}$$

where we note that $\Gamma_0(x)$ corresponds with the classical gamma function

$$\Gamma_0(x) = \lim_{n\to\infty}\frac{n^x n!}{x(x+1)\cdots(x+n)} = \Gamma(x)$$

and we have

$$\Gamma_k(1) = 1$$

$$\Gamma_k(x+1) = \exp\left[\frac{1}{k+1}\log^{k+1} x\right]\Gamma_k(x)$$

Dilcher [25a] showed that

(11.19) $\qquad \log\Gamma_1(x+1) + \gamma_1 x = \sum_{n=1}^{\infty}\dfrac{(-1)^n}{n+1}[H_n\varsigma(n+1) + \varsigma'(n+1)]x^{n+1}$

By definition we have

$$\log\Gamma_1(x) = \lim_{n\to\infty}\left[\frac{x}{2}\log^2 n + \frac{1}{2}\sum_{j=1}^{n}\log^2 j - \frac{1}{2}\sum_{j=0}^{n}\log^2(j+x)\right]$$

and we have

$$\log\Gamma_1(x+1) = \frac{1}{2}\log^2 x + \log\Gamma_1(x)$$

Therefore we have

$$\log\Gamma_1(x+1) + \gamma_1 x = \lim_{n\to\infty}\left[\frac{x}{2}\log^2 n + \frac{1}{2}\sum_{j=1}^{n}\log^2 j - \frac{1}{2}\sum_{j=1}^{n}\log^2(j+x)\right] + \lim_{n\to\infty}\sum_{j=1}^{n}\left[x\frac{\log j}{j} - \frac{x}{2}\log^2 n\right]$$

$$= \lim_{n\to\infty}\left[x\frac{\log j}{j} + \frac{1}{2}\sum_{j=1}^{n}\log^2 j - \frac{1}{2}\sum_{j=1}^{n}\log^2(j+x)\right]$$

and hence we see that (11.17) follows from (11.19).

It was also shown by Dilcher [25a] that

$$\frac{1}{\Gamma_k(x)} = e^{\gamma_k x} \exp\left[\frac{x}{k+1}\log^{k+1} x\right] \prod_{n=1}^{\infty} \exp\left[-\frac{x}{n}\log^k n\right] \exp\left[\frac{1}{k+1}[\log^{k+1}(n+x) - \log^{k+1} n]\right]$$

whereupon it directly follows that

$$\log \Gamma_k(x+1) + \gamma_k x = \lim_{n\to\infty}\left[x\frac{\log^k j}{j} + \frac{1}{k+1}\sum_{j=1}^n \log^k j - \frac{1}{k+1}\sum_{j=1}^n \log^k(j+x)\right]$$

□

Using the Euler-Maclaurin summation formula [5b], Hardy [31a, p.333] showed that the Riemann zeta function could be expressed as follows

(11.20) $\qquad \varsigma(s) = \lim_{n\to\infty}\left[\sum_{k=1}^n \frac{1}{k^s} + \frac{n^{1-s}}{s-1} - \frac{1}{2}n^{-s}\right] \qquad \text{Re}(s) > -1$

and we have

(11.21) $\qquad \varsigma(s) - \frac{1}{s-1} = \lim_{n\to\infty}\left[\sum_{k=1}^n \frac{1}{k^s} + f(s) - \frac{1}{2}n^{-s}\right] \qquad \text{Re}(s) > -1$

where $f(s) = \frac{n^{1-s}-1}{s-1}$. We saw in (10.5.1) that

$$f^{(p)}(1) = (-1)^{p+1}\frac{\log^{p+1} n}{p+1}$$

and we therefore have the limit

$$\lim_{s\to 1}\left[\varsigma(s) - \frac{1}{s-1}\right] = \lim_{n\to\infty}\left[\sum_{k=1}^n \frac{1}{k} - \log n - \frac{1}{2}n^{-1}\right]$$

We thereby recover the familiar limit definition of Euler's constant

$$\gamma = \lim_{n\to\infty}[H_n - \log n]$$

Differentiation of (11.21) formally gives us

$$\frac{d^p}{ds^p}\left[\varsigma(s)-\frac{1}{s-1}\right]=(-1)^p\lim_{n\to\infty}\left[\sum_{k=1}^{n}\frac{\log^p k}{k^s}+(-1)^p f^{(p)}(s)-\frac{1}{2}n^{-s}\log^p n\right]$$

and we have

$$\frac{d^p}{ds^p}\left[\varsigma(s)-\frac{1}{s-1}\right]\bigg|_{s=1}=(-1)^p\lim_{n\to\infty}\left[\sum_{k=1}^{n}\frac{\log^p k}{k}-\frac{\log^{p+1} n}{p+1}-\frac{1}{2}n^{-1}\log^p n\right]$$

Hence we obtain another derivation of (1.22)

$$\gamma_p=\lim_{N\to\infty}\left(\sum_{k=1}^{N}\frac{\log^p k}{k}-\frac{\log^{p+1} N}{p+1}\right)$$

since as shown in (1.22.1) we have $\lim_{n\to\infty}[n^{-1}\log^p n]=0$.

$\square$

Ramanujan [7a, p.198] also showed that

(11.18) $$\phi(x)-\frac{1}{n}\sum_{k=0}^{n-1}\phi\left(\frac{x-k}{n}\right)=\psi(x+1)\log n-\frac{1}{2}\log^2 n$$

and with $x=0$ we have

(11.19) $$\sum_{k=0}^{n-1}\phi\left(-\frac{k}{n}\right)=\gamma n\log n+\frac{1}{2}n\log^2 n$$

Using (11.5) we may express (11.19) as

(11.20) $$n\gamma_1-\sum_{k=0}^{n-1}\gamma_1\left(1-\frac{k}{n}\right)=\gamma n\log n+\frac{1}{2}n\log^2 n$$

Reindexing so that $k=n-j$ we obtain

(11.21) $$n\gamma_1-\sum_{j=1}^{n}\gamma_1\left(\frac{k}{n}\right)=\gamma n\log n+\frac{1}{2}n\log^2 n$$

and it is easily seen that this is equivalent to Coffey's formula (2.19).

$\square$

Since (10.4) is valid in the limit as $x\to 1$, we then have

$$\sum_{n=1}^{\infty}(-1)^n[\varsigma'(n+1)+H_n\varsigma(n+1)]=0$$

and convergence of this series implies that

$$\lim_{n\to\infty}[\varsigma'(n+1)+H_n\varsigma(n+1)]=0$$

We have

$$\varsigma'(n+1)+H_n\varsigma(n+1)=\varsigma'(n+1)+\log n\varsigma(n+1)+[H_n-\log n]\varsigma(n+1)$$

and we then obtain

$$\lim_{n\to\infty}[\varsigma'(n+1)+\log n\varsigma(n+1)]=-\gamma$$

or equivalently

$$\lim_{n\to\infty}\sum_{k=1}^{\infty}\frac{\log n-\log k}{k^{n+1}}=-\gamma$$

## 12. SOME IDENTITIES INVOLVING THE HARMONIC NUMBERS

We recall (11.10.1)

$$\sum_{j=0}^{n-1}\binom{n}{j}\frac{(x)_j}{n!}=\frac{(x)_n}{n!}H_n^{(1)}(x)$$

and substituting

$$(x)_n=\frac{\Gamma(x+n)}{\Gamma(x)}$$

we have

(12.1) $$\sum_{j=0}^{n-1}\binom{n}{j}\frac{\Gamma(x+j)}{\Gamma(x+n)}=H_n^{(1)}(x)$$

With $x=1$ we have

(12.2) $$\sum_{j=0}^{n-1}\frac{1}{(n-j)!}=H_n^{(1)}$$

Differentiating (12.1) gives us

(12.3) $$\sum_{j=0}^{n-1}\binom{n}{j}\frac{\Gamma(x+j)}{\Gamma(x+n)}[\psi(x+j)-\psi(x+n)]=-H_n^{(2)}(x)$$

and substituting

$$H_n^{(1)}(x)=\psi(n+x)-\psi(x)$$

gives us

(12.4) $$\sum_{j=0}^{n-1}\binom{n}{j}\frac{\Gamma(x+j)}{\Gamma(x+n)}H_j^{(1)}(x)=\left[H_n^{(1)}(x)\right]^2-H_n^{(2)}(x)$$

With $x=1$ we have

(12.5) $$\sum_{j=0}^{n-1}\frac{H_j^{(1)}}{(n-j)!}=\left[H_n^{(1)}\right]^2-H_n^{(2)}$$

and substituting (12.2) gives us

$$\sum_{j=0}^{n-1}\frac{1}{(n-j)!}\sum_{k=0}^{j-1}\frac{1}{(j-k)!}=\left[H_n^{(1)}\right]^2-H_n^{(2)}$$

Higher derivatives will obviously result in additional identities.

If $h'(x)=h(x)g(x)$ then we have

$$\frac{d^m}{dx^m}h(x)=\frac{d^m}{dx^m}e^{\log h(x)}=h(x)Y_m\left(g(x),g^{(1)}(x),...,g^{(m-1)}(x)\right)$$

where $Y_m\left(g(x),g^{(1)}(x),...,g^{(m-1)}(x)\right)$ are the (exponential) complete Bell polynomials defined in (2.38) above. Hence, from (12.1) we have

$$\sum_{j=0}^{n-1}\binom{n}{j}\frac{\Gamma(x+j)}{\Gamma(x+n)}Y_m\left(g(x),g^{(1)}(x),...,g^{(m-1)}(x)\right)=(-1)^m m!H_n^{(m+1)}(x)$$

where $g(x)=\psi(x+j)-\psi(x+n)$.

□

Further information regarding the Stieltjes constants is contained in the papers by Coffey (see [11] to [17a] and the papers referred to therein). In [21] we highlighted a connection with the Fresnel integral defined in [28, p.880]

$$C(x) = \frac{2}{\sqrt{2\pi}} \int_0^x \cos(t^2)\,dt$$

and also with the cosine integral $Ci(x)$ defined in [28, p.878] and [1, p.231] as

$$Ci(x) = \gamma + \log x + \int_0^x \frac{\cos t - 1}{t}\,dt = \gamma + \log x + \sum_{n=1}^{\infty} \frac{(-1)^n x^{2n}}{2n(2n)!}$$

This paper contains extracts from two earlier papers ([21] and [22]) with the addition of some new material.

Introduction to the General Theory of Infinite Processes and of Analytic Functions; With an Account of the Principal Transcendental Functions. Fourth Ed., Cambridge University Press, Cambridge, London and New York, 1963.

Donal F. Connon
Elmhurst
Dundle Road
Matfield
Kent TN12 7HD
dconnon@btopenworld.com